\documentclass[12pt,reqno]{article}

\usepackage[utf8]{inputenc}
\usepackage[T1]{fontenc}
\usepackage{lmodern}
\usepackage[english]{babel}
\usepackage{csquotes}
\usepackage{microtype}
\usepackage{enumerate}

\usepackage[a4paper,hmargin=1.75cm,bottom=3cm,top=2.5cm,footskip=3\baselineskip]{geometry}

\usepackage{tikz}
\definecolor{midnightblue}{rgb}{0.1, 0.1, 0.44}
\usepackage[pdfencoding=auto,psdextra,colorlinks,allcolors=midnightblue]{hyperref}
\AtBeginDocument{}
\AtBeginDocument{}

\usepackage[backend=biber,hyperref=true,giveninits=true]{biblatex}
\usepackage{fsmath}
\usepackage{cancel}

\usepackage{bm}
\newcommand{\x}{{\bm{x}}}
\newcommand{\y}{{\bm{y}}}

\newcommand{\ep}{\eps}
\newcommand{\R}{\setR}

\DeclareMathOperator{\graph}{graph}
\DeclareMathOperator{\erf}{erf}
\DeclareMathOperator{\erfc}{erfc}

\newcommand{\ext}{\mathrm{ext}}
\newcommand{\ee}{\mathrm{ee}}

\DeclareMathOperator{\KE}{KE}

\newcommand{\note}[1]{\footnote{\color{red}#1}}
\newcommand{\alert}[1]{{\color{red}#1}}
\newcommand{\old}[1]{{\color{orange!70!red}#1}}
\renewcommand{\note}[1]{\errmessage{ALERT}}
\renewcommand{\alert}[1]{\errmessage{ALERT}}
\renewcommand{\old}[1]{\errmessage{OLD}}

\begin{document}

%\title{First order expansion of a penalized optimal transport functional from DFT}
\title{\vspace{-1cm}First order expansion in the semiclassical limit of \\ the Levy-Lieb functional}
\author{
	Maria Colombo\thanks{
	EPFL, AMCV, Lausanne,
	\url{maria.colombo@epfl.ch}}
	\and
	Simone Di Marino\thanks{
	Università di Genova (DIMA), MaLGa, Genova,
	\url{simone.dimarino@unige.it}}
	\and
	Federico Stra\thanks{
	EPFL, AMCV, Lausanne,
	\url{federico.stra@epfl.ch}}
}

\maketitle

\begin{abstract}
We prove the conjectured first order expansion of  the Levy-Lieb functional in the semiclassical limit, arising from Density Functional Theory (DFT). This is accomplished by interpreting the problem as the singular perturbation of an Optimal Transport problem via a Dirichlet penalization.% and obtain the first order expansion. 

%We study the singular perturbation of an Optimal Transport problem via a Dirichlet penalization and obtain the first order expansion. The functional originates from Density Functional Theory (DFT) and we provide the first rigorous proof of the conjectured first order term in the expansion of the semiclassical limit.
\end{abstract}

\setcounter{tocdepth}{1}
\tableofcontents

\paragraph*{Acknowledgements.}
%The authors thank the anonymous referee for pointing out a mistake in a previous version of the manuscript, which led to a substantial improvement in the paper.
M.C.\ and F.S.\ have been supported by the SNSF Grant 182565.

\clearpage

%!TEX root = ../semiclassical.tex
\section{Introduction}

%\alert{\textbf{References}: non mi è chiaro se con [Bindini-DePascale] si intende sempre \cite{BDeP} (praticamente sempre) o \cite{BDeP2} (mi sembra quasi mai); non mi è chiaro chi sono [Friesecke1]--> \cite{CFK} e [Friesecke2]--> \cite{CFK2}.}

A revolutionary approach to finding the ground state energy of a many-body electron system was developed in the 60s by Hohenberg and Kohn \cite{HK64}: Density Functional Theory (DFT). Their idea, subsequently  formalized mathematically by Lieb and Levy \cite{LL83}, can be seen as breaking up the minimization over all wave functions into a first minimization over those wavefunctions having a given one-particle density and then minimizing the resulting function over the one-particle density. In particular we define 
\begin{equation}
\label{eqn:Fhbar}
F^\eps(\rho) = \inf\set*{ \int_{\setR^{dN}} \left(\frac\eps 2 \abs{\nabla\psi}^2
	+ V_{\ee} \abs{\psi}^2 \right) \d x_1 \dots \d x_N}
	{\psi\in H^1(\R^{dN}),\ \psi\mapsto\rho},
\end{equation}
where $\psi \mapsto \rho$ means that the one electron density of $\psi$ is $\rho$, that is, for all $i=1, \ldots, N$, we require $\rho(x_i)= \iint | \psi|^2(x_1, \ldots, x_n) \, dx_1 \, \ldots \, \widehat{ d x_i } \, \ldots d x_N$, and moreover
\[
V_{\ext}(x_1,\dots,x_n) = \sum_{i<j} \frac{1}{\abs{x_i-x_j}}.
\]
%(i.e.  $(e_i)_{\sharp} ( | \psi|^2 \mathcal{L}^{dN} )= \rho \mathcal{L}^d$ where $e_i(x_1, \ldots, x_N)=x_i$
Then, for any external potential $V_{ext}$ the corresponding energy of the ground state will be equal to 
\[
E(V_{\ext}) = \inf\set*{ F^{\hbar^2}(\rho) + \int_{\setR^d} V_{\ext} \d\rho}
	{\rho\in\Meas_+(\setR^d),\ \rho(\setR^d)=N}.
\]
In particular $F^{\hbar^2}(\rho)$ is a \emph{universal functional} in the sense that depends only on the number of electrons, whereas the dependence on the external field appears only in the term $\int V_\ext \d\rho$ in the outer minimization. It becomes then fundamental for applications to approximate $F^\hbar(\rho)$ and compute this value efficiently.
In recent years a new approach gained importance, relying on the analysis of the Strictly Correlated Electrons (SCE) case. The interest in this approach resides in two main reasons: it is mathematically rigorous since it is a limiting procedure starting from the exact functional, and it is highly non-local in nature, thus it can be thought as complementary in some sense to the more classical Local Density Approximation (LDA) approach \cite{LewinLDA}. The idea is to have a parameter $\beta>0$ which tunes the strength of the interaction between the electrons and the functional
\[
F^{\hbar^2}_\beta(\rho) = \inf\set*{ \int_{\setR^{dN}} \left(\frac{\hbar^2}2 \abs{\nabla\psi}^2
	+ \beta V_{\ee} \abs{\psi}^2 \right) \d x_1 \dots \d x_N}
	{\psi\in H^1(\R^{dN}),\ \psi\mapsto\rho}.
\]
The SCE limit considers the case $\beta \to \infty$, the asymptotic expansion of $F^{\hbar^2}_{\beta} (\rho)$ in $\beta$, to then deduce information about $\beta<\infty$.
This limiting procedure has more equivalent formulations. Using the homogeneity of $V_{ee}$, $F^{\hbar^2}_{\beta} $ is computed from $F^{\hbar^2}$ by scaling, namely $F^{\hbar^2}_{\beta}(\rho)= \beta^2 F^{\hbar^2} (\rho_\beta)$ where $\rho_{\beta}(x)= \beta^{-d} \rho(x/\beta)$. %: this is called sometimes the slow-varying regime, since the variations of the density $\rho_{\beta}$ are at scale $\beta$.
Another equivalent approach considers a varying  kinetic energy coefficient: %we remark that this procedure is not physical since $\hbar$ is a physical constant, but it is a useful mathematical simplification. 
since $F^{\hbar^2}_{\beta} (\rho) = \beta F^{\hbar^2/\beta}(\rho)$, the asymptotic expansion of $F^{\hbar^2}_{\beta}$ as $\beta \to +\infty$ can be studied by means of the asymptotic expansion of $F^{\ep}(\rho)$ as $\ep\to 0$, which is mathematically convenient since the zeroth-order term of the expansion does not need to be renormalized.
This is formally the same as changing the value of $\hbar$ in the original definition of the Levy-Lieb functional, but the rigorous justification of this non-physical procedure (since $\hbar$ is a physical constant) relies on the previous explanation.
%{\color{red} DISCORSO SU SCALING E STRONG INTERACTION CHE SONO LA STESSA COSA DI $\hbar \to 0$}

%Looking at the explicit dependence of the DFT functional from the parameter $\hbar$, we can think of $F^{\hbar}_{HK-LL}(\rho)$. What is usually employed as an approximation is the regime $\hbar \gg 1$, while we are interested in the case $\hbar \ll 1$: the idea is that in order to get the intermediate values one must interpolate. Of course more is known about the two extreme regimes the better approximation we can find. While the case $\hbar \gg 1$ is rather well known (citazioni), 

Gori-Giorgi, Savin and Seidl were the first to conjecture the SCE limit case for $\ep \to 0$ in \cite{GGSV}: this new functional found its use for example in \cite{depa2018dissociating, CA2, CA1, CA3}. 

$$ \lim_{\ep \to 0} F^{\ep}(\rho) = F_{OT}(\rho) := \inf \left\{ \int_{ \R^{3N} } V_{\ee} \, d \gamma \; : \; \gamma \in \Pi_N(\rho) \right\},$$
where $\Pi_N(\rho)$ is the set of probabilities in $\R^{dN}$ which induce $\rho$, namely such that its push-forward through any projection $e_i: (x_1, \ldots, x_N) \mapsto x_i$ is $\rho$. We will denote in the sequel with $\Pi_0(\rho)$ the set of $\gamma \in \Pi_N(\rho)$ which are minimizers for $F_{OT}(\rho)$. Moreover, in the physics literature, ansatz of minimizers in the 1D case and radial case were conjectured in \cite{S99, SGSmaps07}: in the 1D case the conjecture was proven to hold in \cite{CDD}, while in the radial case various counterexamples were found in \cite{CS, Seidletc, BDK}.
The functional $F_{OT}$ has been studied in the last years also with regard to the limit of infinitely many particles: the zeroth order expansion for $N \to \infty$ was investigated in \cite{CFP}, proving the mean field limit, while the first order was proven independently in \cite{LewinUniformGas} and \cite{CP2,CP}, with ideas coming from the seminal papers \cite{PetracheSerfaty,SandierSerfaty}, and in connection with the Lieb-Oxford inequality \cite{LIEB1979444,LiebOxford}.

Shortly after \cite{GGSV}, the next order in the asymptotic expansion of $F^\eps$ for $\eps\to0$ was conjectured in \cite{ZPO}:
\begin{equation} \label{eqn:expansion} F^{\ep}(\rho) = F_{OT}(\rho) + \sqrt{ \ep} F_{ZPO}(\rho) + O(\ep^{3/4} ); \end{equation}
together with an explicit conjecture for the \emph{zero point oscillation} functional $F_{ZPO}(\rho)$ in the case $d=1$, involving the eigenvalues of the Hessian of an effective potential described in the following. %In order to write explicitly the functional $F_{ZPO}$ we 
As in \cite{BDG}, we introduce the potential $u$ which solves the dual problem for $F_{OT}$; its existence and regularity have been studied in \cite{Auglog, DeP, BCD, CDMS}. It is known to be unique if the support of $\rho$ is connected, $u$ is Lipschitz, $F_{OT}(\rho)= N\int u \d \rho$,% and crucially we have 
\begin{equation}\label{eq:V}
V(x_1, \ldots, x_N) = V_{ee}(x_1, \ldots, x_N) - u(x_1) - \ldots - u(x_N) \geq 0 \qquad \forall x_1, \ldots, x_N \in \R^d;
\end{equation}
and the optimal plans are supported in the set $\{V=0\}$. In particular $D^2V$ is well defined $\gamma$-almost everywhere for every optimal plan $\gamma$ and so we can define
\[
F_{ZPO}(\rho):= \inf \set*{ \frac12\int_{\R^{dN}}\tr\oleft(\sqrt{D^2V}\right) \d\gamma }{ \gamma \in \Pi_0 (\rho)}.
\]
In the physics literature $V$ is called the effective potential; $u$ is the Lagrange multiplier associated with the constraint on the marginals, thought as a confining potential which keeps the optimal plan concentrated on minimal energy configurations with respect to the effective potential; $\frac12 \tr\oleft(\sqrt{D^2V}\right)$ is the lowest eigenvalue of the harmonic oscillator associated with the quadratic form given by $D^2V$.

From the mathematical  viewpoint, the $0$-th order approximation was obtained in \cite{CFK} for the bosonic case and for the fermionic case only when $N=2$. Then it was settled in \cite{BDeP} for $N=3$ with the help of spins; in \cite{Lewin} with a different approach involving mixed states (that is, considering the convex relaxation  $\bar{F}_{\ep}(\rho)$ of $F_{\ep}(\rho)$); finally, \cite{CFK2} provided a proof using only pure states for every $N$.

%even the proof for the $0$-th order approximation was closed after $10$ years by [Friesecke2] after a number of papers dealing with it [Friesecke1] \cite{BDeP,BDeP2,Lewin}. Notice that in [Friesecke1] the proof of the convergence was already complete for the bosonic case: the fermionic case was more deceptive and it was solved in [Friesecke1] only in the case $N=2$, then in \cite{BDeP,BDeP2} for $N=3$ with the help of spins, and in \cite{Lewin} using a slightly different approach with mixed states (that is, considering the convex relaxation of $\bar{F}_{\ep}(\rho)$ of $F_{\ep}(\rho)$); in the end in [Friesecke2] a complete proof was presented using only pure states for a generic $N$.

Regarding the $1$-st order, the best known bound in \cite{Lewin} proves that $\bar{F}_{\ep} (\rho) \leq  F_{OT}(\rho) + \sqrt{\ep} C_{\rho}$ for every admissible $\rho$, namely an upper bound on the conjectured order of the next term. Our main contribution can be summarised in the following

\begin{theorem}\label{thm:main_intro} Let $\rho \in \mathcal{P}(\R^d)$ be such that $F^{\ep}(\rho) < \infty$, $u$ the optimal Kantorovich potential associated to the problem $F_{OT}(\rho)$ and suppose that $V$ as in \eqref{eq:V} is such that $V \in C^2_{loc} ( \{ V=0 \})$. Then 
\begin{itemize}
\item for every $N \geq 2$ and $d \geq 1$ we have
\begin{equation}\label{eq:Eeps-liminf}
\liminf_{\ep \to 0} \frac { F^{\ep}( \rho) - F_{OT}(\rho) }{\sqrt{\ep}} \geq F_{ZPO} (\rho);
\end{equation}
\item if $d=1$, $N=2$ and $\rho \in C^1(\R)$ is a positive probability density, then
\begin{equation}\label{eq:Eeps-limsup}
\limsup_{\ep \to 0} \frac { F^{\ep}( \rho) - F_{OT}(\rho) }{\sqrt{\ep}} \leq F_{ZPO} (\rho). \end{equation}
\end{itemize}

\end{theorem}

Our approach is based on the variational formulation for the quantity $ \frac{ F^{\ep}(\rho) - F_{OT}(\rho)}{\sqrt{\ep}}$. 
Using that $F_{OT}(\rho)=N\int u \, d \rho =\int_{\R^{dN}} (u(x_1) + \ldots + u(x_N) ) | \psi |^2(x_1, \ldots, x_N) \,d x $ for every $\psi \mapsto \rho$, we obtain that 
\begin{align}
\frac{ F^{\ep}(\rho) - F_{OT}(\rho)}{\sqrt{\ep}}
&= \min \set*{ E_{\ep}(\psi)} {\psi\mapsto\rho}, \label{eq:Feps-FOT}\\
E_{\ep}(\psi) &=  \int_{ \R^{dN} }\ep^{1/2} \frac{ |\nabla \psi (x)|^2}2 + \ep^{-1/2} V(x) | \psi (x)|^2 \d x. \label{eq:Eeps}
\end{align}

%In order to study the limit of minima of $E_{\ep}(\psi)$ it is natural to consider its $\Gamma$-convergence. There are two kind of $\Gamma$-convergence which can be interesting: the first one is retaining the marginal constraint (i.e. the one of physical interest), while the second one, which is easier from the technical standpoint, is done without the marginal constraint along the sequence.

%Our first results concerns precisely the latter situation, which will be also propedeutic to the former, since it helps us understanding the structure of the minimizing sequence for $E_{\ep}$. We show that
The core of Theorem~\ref{thm:main_intro} relies in a precise convergence result on $E_\eps$. Firstly we study the convergence when we \emph{drop the marginals constraint}, which simplifies our analysis but is not the one of physical interest, and we show roughly speaking that
\begin{equation}\label{eqn:conv-summary}
E_{\ep}(\psi) \stackrel{\Gamma}{ \longrightarrow}
\begin{cases} \displaystyle \frac12\int_{\R^{dN}}\tr\oleft(\sqrt{D^2V}\right) \d\gamma \quad & \text{ if $\gamma ( \{ V>0 \}) =0$} \\
+\infty & \text{ otherwise.}
\end{cases}
\end{equation}
The notion of convergence in \eqref{eqn:conv-summary} is the so called $\Gamma$-convergence, the right variational tool to take the limit in sequences of minimizers, as it happens in Theorem~\ref{thm:main_intro}. For brevity, we don't introduce this abstract notion here but we just present its outcomes in \autoref{thm:nonfixed-marginals}.% in particular we get in the limit the conjectured $F_{ZPO}$ functional, albeit using wavefunctions which satisfy the marginal constraint only at the limit.
%The main result is naturally stated as $\Gamma$-convergence, since we are looking at the re-scaled behavior of the energy, as $\ep \to 0$, of any minimizing sequence for the non-scaled version: in this sense it is a \emph{second order} $\Gamma$-convergence.

\begin{theorem}[$\Gamma$-convergence with free maginals]\label{thm:nonfixed-marginals}
Let $V \in C^0(\R^d)$ with $V \geq 0$ and $V $ locally  $C^2$ in a neighborhood of $\{V=0\}$. %and $\liminf_{|x| \to \infty} V(x)>0$.
Then we have  \begin{itemize}
\item[(a)] for every sequence $\psi_{\ep} \in H^1(\R^d)$ with $\int | \psi_{\ep}|^2=1$ and $\liminf F_{\ep} (\psi_{\ep})  < \infty$, there exists $\gamma \in \mathcal{M}(\R^n)$ with  $\gamma(\{ V>0 \})=0$ and $\gamma(\R^d) \leq 1$ such that $ | \psi_{\ep}|^2 \weakto \gamma$ up to subsequences and
\[
\liminf_{n \to \infty} E_{\ep}(\psi_{\ep}) \geq \frac12\int_{\R^{dN}}\tr\oleft(\sqrt{D^2V}\right) \d\gamma;
\]
\item[(b)] for every $\gamma \in \mathcal{P}(\R^n)$ with   $\gamma(\{ V>0 \})=0$, there exists a \emph{recovery sequence} $\psi_{\ep} \in H^1(\R^d)$ with $\int | \psi_{\ep}|^2=1$ such that $ | \psi_{\ep}|^2 \weakto \gamma$ and
\[
\limsup_{n \to \infty} E_{\ep}(\psi_{\ep}) \leq \frac12\int_{\R^{dN}}\tr\oleft(\sqrt{D^2V}\right) \d\gamma.
\]
\end{itemize}
\end{theorem}
In order to deal with the actual first order expansion of $F_{\ep}(\rho)$ we have to enforce the marginal constraint: while for the lower bound we can use \autoref{thm:nonfixed-marginals} (a), in order to prove the upper bound we need to construct a recovery sequence $\bar{\psi}_{\ep}$ with an almost optimal energy, but with the additional constraint on the marginals. This is in general nontrivial, since the Lavrentiev gap phenomenon could appear as happens for instance in an example of \cite{CFK2}, where the constrained limit energy is smaller than the unconstrained one.  This result is established in \autoref{thm:rec}. %Notice that the fact that the limit energy is the same even if we put the maginal constraint is not at all evident: in [Firesecke 2018] the authors show in fact a functional for which the two limits are different, due to a Lavrentiev gap phenomenon.

\begin{theorem}\label{thm:rec}
Let $\rho\in C^1(\setR)$ be a strictly positive probability density with finite kinetic energy, let $\gamma$ be the optimal transport plan for the problem with Coulomb cost and marginals $\rho$, let $u$ be its Kantorovich potential and assume that $u'' \in  L^\infty(\R)$. Let $V(x,y) = |x-y|^{-1} - u(x) - u(y)$. %Let $V\in C^2\bigl([0,1]^2;[0,\infty)\bigr)$. 
 
 Then %for every $\gamma\in\Prob([0,1]^2)$ %with $\gamma(\{V>0\})=0$ 
 there exists a recovery sequence $\psi_\eps\in H^1(\setR^2)$ such that $\int\abs{\psi_\eps}^2=1$, $\abs{\psi_\eps}^2\weakto\gamma$, $\abs{\psi_\eps}^2$ has the same marginals as $\gamma$ and
\[
\limsup_{\eps\to0} \int_{\setR^2} \biggl(
	\eps^{1/2} \frac{\abs{\nabla\psi_\eps}^2}2
	+ \eps^{-1/2} V \abs{\psi_\eps}^2 \biggr) \d x
\leq \frac12 \int_{\setR^2} \tr\bigl(\sqrt{D^2V}\bigr) \d\gamma .
\]
\end{theorem}
Theorem~\ref{thm:rec} is stated with a rather implicit assumption, namely that $u'' \in  L^\infty(\R)$. This is essentially equivalent to the condition \eqref{hp:tails} on the growth of the tails of $\rho$. For instance, we show in \autoref{lem:good-rho} that \autoref{thm:rec} can be applied to any strictly positive $C^1$ probability density $\rho$ such that $\KE(\rho)<\infty$ and
\begin{equation}
\label{hp:tails}
\liminf_{\abs{x}\to\infty}{} \abs{x}^3 \rho(x) > 0.
\end{equation}
%We conclude the introduction with some comments and an outline of the paper. 

\paragraph{Singular perturbations of optimal transport.} Since without loss of generality we can assume $\psi>0$, we rewrite the functional $F^{\ep}(\rho)$ in terms of $\gamma= | \psi|^2 $
$$ F^{\ep}(\rho) = \inf \left\{ \int_{\R^{dN}} V_{ee} (x_1, \ldots, x_n) \, d \gamma  +  \ep \int_{\R^{dN}} \frac { |\nabla \gamma|^2}{\gamma} \, d x \, : \; \gamma \in \Pi_N(\rho) \right\}.$$
We think of $F^{\ep}$ as a penalized multi-marginal optimal transport problem, where the penalization on the plan $\gamma$ is the Fisher information with respect to the Lebesgue measure $I(\gamma)=\int \frac { |\nabla \gamma|^2}{\gamma} \, d x$. The ZPO conjecture is about the asymptotic expansion of this penalized OT problem. Other results which are similar in spirit regard only the case $N=2$ with the entropic penalization $Ent(\gamma)= \int \gamma \ln(\gamma) \, dx$ and either the Monge cost (\cite{DiMarinoLouet}, only in 1D), or the quadratic distance cost (\cite{ConfortiTamanini} in a general Riemannian manifold) or a smooth convex cost (\cite{Pal}, in $\R^d$). However, in \cite{ConfortiTamanini} and in \cite{Pal} the peculiar structure of the problem and an additivity property of the entropy are heavily used. The approach closest to ours is \cite{DiMarinoLouet}, but the construction of the recovery sequence is heavily simplified by taking advantage of working in the region $\{ T = {\rm id}\}$. To our knowledge, the present result is the first next-order asymptotics with marginal constraint proven with a rather general strategy.

The setting in \autoref{thm:nonfixed-marginals} is rather general and works for any nonnegative potential which is locally $C^2$ around its zero set (we recall that here $\rho$ does not play a role). Similar studies in the literature look at the spectrum of the operator $- \Delta + \ep V$ as $\ep \to 0$ (see for example \cite{SimonI,SimonII}), while our point of view is slightly different. Instead of focusing on the behavior of the eigenvalues we want to know what is the minimial asymptotic energy given the asymptotic behavior of the wavefunction (that is, their weak limit $\gamma$). This is of course in view of applying the result to a generic $\gamma \in \Pi_N(\rho)$ which might not come from a sequence of eigenvectors.

%The behavior of the ground state energy as $\ep \to 0$ has been already studied in the case $\argmin V$ consists in a finite number of points, with also more general behaviors of the potentials around the minimal points: in the literature however the focus is not on the actual ground state energy but rather on the behavior of the whole spectrum of the operator. However, since we are interested only in the minimal value, we can allow for a more general structure of $\argmin V$, but still retaining the $C^2$ behavior around $\{V=0\}$ which guarantees the correct scaling $\sqrt{\ep}$ for the expansion of $F_{\ep}(\rho)$.  (CITAZIONIIIII)
\autoref{thm:nonfixed-marginals} builds on the idea that, to have the correct minimal energy, the wavefunction must localize around $\{V=0\}$ and behave locally like a gaussian with variance depending on $D^2V$. Moreover this minimal energy can be achieved by a recovery sequence built superposing these gaussians, which will be elongated along the zero set of $V$. In the proof of \autoref{thm:nonfixed-marginals} we construct a wavefunction for $\ep>0$ which has  the proven almost correct energy and therefore can be of interest also for numerical purposes.

\paragraph{Ansatz on wavefunctions with asymptotically minimal energy.} The ansatz of a recovery sequence obtained as a superposition of Gaussians is already present in \cite{ZPO}, but then, instead of superposing elongated Gaussians, a local nonlinear change of variables is considered, as well as Gaussians in the variable perpendicuar to the manifold $\{V=0\}$. In our opinion, there are no guarantees on the smallness of the kinetic energy in the direction tangent to the manifold.
In \cite{ZPO} the constraint on the marginals is also diuscussed via the introduction of a Lagrange multiplier. However, the latter has no reason to exist for $\ep>0$ and it is not clear the claim that it intervenes in the asymptotic expansion only at higher order terms.
%
%\alert{Menzionare la costruzione dei fisici/chimici/numerici e giustificare che la nostra costruzione è un po' migliore in qualche senso perché ci consente di avere una remaining mass più regolare che possiamo deconvolvere con energia cinetica controllata.
%
%Nel conto dei numerici, i marginali sono sbagliati, quindi ci sta che gli venga con energia minore.
%}

\paragraph{Future perpectives.}
It would be interesting proving the ZPO conjectured asymptotics for every $\rho$ and in any dimension, in particular removing the heavy tails assumption \eqref{hp:tails} since one expects physical densities to have exponential tails. While the restriction to $N=2$ in \autoref{thm:rec} is just used to simplify technical details and not make the notation heavier, a different strategy would be needed in dimension $d \geq 2$. In fact we rely heavily on the one dimensional setting for the existence, regularity and structure of optimal maps. % In particular we know that if $\rho$ is smooth, the marginals of a localized part of the optimal plan are still smooth: this allow us to build locally our solution. 
This is an open problem if $d \geq 2$ and $N \geq 3$. Moreover in $d\geq 2$  the addition of the constraint of the antisymmetry on $\psi$ is non-trivial (instead in $d=1$ it is sufficient to consider carefully the signs): however we believe that, as in \cite{Lewin}, this issue can be circumvented considering mixed states.

\paragraph{Outline of the paper.}
%We present here overviews of the proof of \autoref{thm:nonfixed-marginals} and  \autoref{thm:rec}, as well as a proof of \autoref{thm:main_intro}.

The main idea in \autoref{thm:nonfixed-marginals}, once we show that every finite energy sequence needs to concentrate on the set $\{V=0\}$,  is to analyze the energy locally around the zero set of $V$: being at a minimum, the first non-zero term in the Taylor expansion of $V$ is the quadratic form corresponding to the second derivative. But then the energy is analogous to that of an anisotropic harmonic oscillator, of which we know the ground state (a Gaussian with proper covariance) and ground energy (see \autoref{lem:main} and \autoref{lem:delta}).
On one hand, an approximation of $V$ with its Taylor expansion will provide the lower bound, while in order to get an upper bound, we use Gaussians elongated along the zero set of $V$, and then we superpose them:
\begin{equation}
\label{eqn:ansa}
f(x) = \int f_{\eps,x'}(x-x') \d\gamma(x')
\end{equation}
where $f_{\eps,x'}(x)$ is a suitable truncation of the Gaussian
$C_{\eps,x'} \exp\oleft(-\eps^{-1/2}x^t {\sqrt{\nabla^2V(x')}}{} x\right)$.
The convexity of $E_{\ep}$ when written in terms of $\gamma= | \psi|^2$ lets us conclude the upper bound with no marginal constraint.
Along our recovery sequence the contribution of the two competing functionals, namely the  rescaled kinetic and potential energies, is asymptotically equal. This equipartition of energy is expected since, in the case of convergence to a single $\delta_0$ and $V$ quadratic, the kinetic energy rescales as $\KE(\gamma_\lambda)=\lambda^2\KE(\gamma)$ and the potential energy as $\int V\d\gamma_\lambda=\lambda^{-2}\int V\d\gamma$, where $\gamma(x)=\lambda\gamma(\lambda x)$.

%The proof of \autoref{thm:nonfixed-marginals} reveals how to construct a recovery sequence without fixing the marginals. 
To tackle \autoref{thm:rec}, namely to modify our ansatz \eqref{eqn:ansa} to keep the marginals constraint, 
we % decide to go for a divide and conquer strategy. We
 split the total mass in two components: a main portion of the mass and a remaining mass.
The main part consists of most of the mass and must have the correct asymptotic energy.
In \autoref{sec:rec} we prove that a suitable superposition of truncated Gaussians (introduced in \autoref{sec:gaussians}), similar to the one adopted for \autoref{thm:nonfixed-marginals}, produces a sequence with asymptotically optimal energy. The superposition is carefully chosen so as to leave a remaining mass with good properties for the subsequent construction.

The remaining mass is used to correct the marginals (see \autoref{sec:deconvolution}) %without disrupting the energy estimates of the main part 
with asymptotically vanishing cost (see \autoref{sec:remaining}). In the literature general methods for re-instating the marginals have been proposed, but they usually miss the correct order of magnitude for the kinetic energy or the potential energy: the only construction that controls in a fine way both of them has been proposed by \cite{Bindini}, and then used then by \cite{Lewin} to retrieve in fact an upper bound with the correct order of magnitude. However the problem of this approximation is that it does not catch the correct local behavior because it is a superposition of gaussians with a fixed variance throughout the space. The latter issue motivates our choice of applying the deconvolution procedure in an unconventional way, namely only to the remaining mass, to connect the remaining marginals.
%In order to treat it, we want to build a plan which has the correct remaining marginals and has negligible energy.
%The idea is to apply a deconvolution to an optimal transport plan connecting the remaining marginals.
%There are several deconvolution procedures available in the literature, but some are not good enough to obtain the right potential bounds, that we need to be of order $\eps^{-1/2}$.
Here we use the construction by Bindini in \cite{Bindini} since we are not interested in the sharp local energy but only the correct order of magnitude: in this way we do not have to worry about the kinetic energy, which is bounded by the kinetic energy of the remaining marginals.
This construction is presented in \autoref{sec:deconvolution}, where we show that this deconvolution procedure does not degrade too much the energy estimates. In particular, we estimate the potential energy energy of the deconvolved plan in terms the potential energy of the original plan connecting the remaining marginals up to the correct order in $\eps$. This section unveals the most delicate issue of the paper, needed to meet the assumption of the deconvolution: in  \autoref{sec:remaining} we prove that there exists $\tilde{\gamma}_H$ such that $\int V \, d \tilde{\gamma} = o( \sqrt{\ep})$.

This last estimate of the potential energy means that the remaining marginals $\rho-\rho^1_\eps$ and $\rho-\rho^2_\eps$ can be connected by a plan which is concentrated near $\{V=0\}=\graph T$. To prove it, we show that $T_\#(\rho-\rho^1_\eps)$ is close to $\rho-\rho^2_\eps$ in the Wasserstein sense. The marginals are obtained by convolution with very anisotropic kernels with size $\eps^{1/4}$ by $\beta^{1/2}$, which is much larger. A naive argument about balance of mass can show that the Wasserstein distance is at most $\beta^{1/2}$, but this is not sufficient. To show the desired estimate, we therefore need to employ a delicate analysis based on the Benamou-Brenier estimate of the transport cost, the linearization of the map $T$ and a technical \autoref{prop:conv-prod} to take advantage of the special structure of the projections of our truncated Gaussian kernels, introduced in \autoref{sec:gaussians}.

\section{$\Gamma$-limit of the unconstrained problem}

In this section we prove \autoref{thm:nonfixed-marginals}. The strategy is to prove that for a finite energy sequence, the measure $| \psi_{\ep}|^2 \,d x$ will concentrate around $\{ V=0\}$, and here the potential term can be well approximated by its second order Taylor expansion. More precisely, at each minimum point $V$ will have some degenerate directions, and some other nondegenerate directions along which it is bounded from below by a suitable quadratic form.
 %is a semidefinite quadratic form on the minimal points of $V$.
%
%TO DO
%The following lemma analyzes the behavior of a $C^2$ nonnegative function locally around the zero set: in particular it shows that up to a $C^2$ change of variable we have that $V$ resembles a positive semi-definite quadratic form, even one that has not a full rank. The crucial aspect is an exact bound from below: after an orthogonal change of coordinate $(x,y)$ ($x$ is the \emph{active} coordinate and $y$ is the \emph{passive} coordinate) we can estimate from below $V$ by a fixed quadratic form in the $x$ direction, which is centered for each $y$ around a point which is "close" to the spine $x=0$.
%
%We call it straightening lemma since the idea is to understand the change of variable that straighten this locus of minima, which should be the set $\{V=0\}$, to the spine.

We can thus look at the linearized problem: using the explicit solution  (a gaussian with known variance) when $V$ is a positive definite quadratic form, we estimate locally from below the energy around minimal points (\autoref{lem:main}), proving in particular  \autoref{thm:nonfixed-marginals}(a) for $V$ a nondegenerate quadratic potential. In \autoref{lem:delta} for any semidefinite quadratic form we build a sequence that attains the minimum energy in the limit, while concentrating the mass at a point. This proves  \autoref{thm:nonfixed-marginals}(b) for $\gamma=\delta_0$ and the construction is inspired by the equality cases in \autoref{lem:main}, but needs to deal with the possible degeneracy of $V$ at $0$ and with a localization, since we have no information on the potential outside the origin.
%We start with the estimate from below in the case of a positive definite quadratic form: We will consider the localized minimal energy
%
%$$ E^{\varepsilon}_{\Omega} = \min \left\{ \int_{\Omega} \ep | \nabla  \psi|^2  +  | A x|^2 | \psi|^2 \, dx \; : \; \int_{\Omega} | \psi|^2 =1, \, \psi \in H^1(\Omega) \right\}.$$
%
%Our main estimate will be the following.
\begin{lemma}\label{lem:main}  Let $A$ be a positive definite symmetric matrix such that $\lambda Id \leq A \leq \Lambda Id$ for some $\Lambda > \lambda > 0$.  Whenever $B_{2r}(0) \subset \Omega \subseteq \R^n$, there exist a constant $ C$ depending only on $\lambda$ and $r$,  such that  if $\sqrt{\ep} \leq \lambda^2 r^2/n\Lambda$, we have% $ E^{\varepsilon}_{\Omega} \geq \frac{\tr (A)\sqrt{\ep}}{1+C \sqrt{\ep} }$, that is 
\begin{equation}\label{eqn:Omega}
 \int_{\Omega} \Bigl( \ep |\nabla \psi|^2  + | Ax|^2 | \psi|^2 \Bigr) \, dx \geq  \frac { \tr (A)\sqrt{\ep}}{1+C \sqrt{\ep}} \int_{\Omega} | \psi|^2\, dx \qquad \forall \psi \in H^1(\Omega).
\end{equation}
Moreover, if $\Omega= \R^n$ we can take $C=0$ in the previous formula and equality holds if and only if $\psi$ is a multiple of $e^{-x^tAx/2\sqrt{\ep}}$.
\end{lemma}

\begin{proof} First let us prove that we can take $C=0$ if $\Omega=\R^n$. In fact, for every $\psi \in C_c^{\infty}$ we can define $g(x)=\psi(x)e^{x^tAx/2}$. Then we have
$$ |\nabla \psi|^2= \left| \nabla g e^{-x^tAx/2} - g Ax e^{-x^tAx/2} \right|^2 = \left( |\nabla g|^2 - \nabla ( g^2) \cdot Ax  + g^2 |Ax|^2 \right)e^{-x^tAx}.$$
Integrating this identity $\R^n$ and then using the integration by parts formula we get
\begin{align*} \int_{\R^n} | \nabla \psi|^2 \, dx &= \int_{\R^n}  | \nabla g|^2 e^{-x^tAx} \, dx+ \int_{\R^n} |Ax|^2 g^2 e^{-x^tAx}\, dx  + \int_{\R^n} g^2 \nabla \cdot ( Ax e^{-x^tAx}) \, dx\\
&\geq  \int_{\R^n} |Ax|^2 g^2 e^{-x^tAx}\, dx + \tr (A) \int_{\R^n} g^2 e^{-x^tAx}\, dx - 2 \int_{\R^n} |Ax|^2 g^2 e^{-x^tAx}\, dx \\
& =   - \int_{\R^n} | Ax|^2 | \psi|^2 \, dx + \tr (A),
\end{align*}

where we used that $|\nabla g|^2$ is nonnegative and that $g^2e^{-x^tAx} = | \psi|^2$. Using this inequality with $A/\sqrt{\ep}$ and using the density of $C_c^{\infty}$ in $H^1_0(\R^n)$ we have that
\begin{equation}\label{eqn:Rn}
\int_{\setR^n} \Bigl( \ep |\nabla \psi|^2  + | Ax|^2 | \psi|^2 \Bigr) \, dx \geq  \tr (A)\sqrt{\ep} \int_{\setR^n} | \psi|^2\, dx \qquad \forall \psi \in H^1(\setR^n);
\end{equation}
moreover we can have equality in \eqref{eqn:Rn} using $g \equiv C_{\ep}$. In particular  the previous computation shows that $\psi_{\varepsilon} = C_{\varepsilon} e^{-x^tAx/2\sqrt{\ep}}$ is the minimizer for $E^{\ep}_{\R^n}$ (here $C_{\varepsilon}$ is the normalization constant such that $\int | \psi_{\varepsilon}|^2 =1$).

In the general case let us consider the cut off function $f_{r} (x)= \min\{1, \bigl( \frac {|x|}{r} - 2\bigr)_+\}$. We have that $|\nabla f_r | \leq \frac 1r \chi_{B_{2r} \setminus B_r }$ and $f_r=1$ on $B_r$ while $f_r=0$ on $B_{2r}^c$. In particular for every function $\psi \in H^1(\Omega)$ we have $f_r \psi \in H^1_0(\R^n)$. But then we can apply \eqref{eqn:Rn} in order to obtain
\begin{equation}\label{eqn:cutoff} \int_{\Omega} \Bigl( \ep| \nabla ( \psi f_r )|^2 +  |Ax|^2 | \psi|^2f_r^2 \Bigr)\, dx  \geq \sqrt{\ep}\tr (A) \int_{\Omega} | \psi|^2 f_r^2.
\end{equation}
Now, by Young inequality we have that for every $\delta <1$
\begin{equation}\label{eqn:grad}| \nabla ( \psi f_r )|^2 = | \nabla \psi|^2 f_r^2 + | \nabla f_r|^2 | \psi|^2 + 2 \nabla f_r \cdot \nabla \psi f_r \psi \leq ( 1+ \delta) | \nabla \psi|^2 f_r^2 + (1+\delta^{-1}) | \nabla f_r|^2 | \psi|^2.
\end{equation}
Notice also that, since $|Ax|^2 \geq \lambda |x|^2 \geq \lambda^2 r^2 $ if $|x|\geq r$, we have
\begin{equation}\label{eqn:gradient}
 (1+\delta^{-1})| \nabla f_r|^2 = \frac{\delta+ 1}{\delta r^2 } \chi_{B_{2r} \setminus B_r}  \leq \frac 2{\delta r^2 } \cdot \frac{|Ax|^2}{\lambda^2 r^2};
 \end{equation}
moreover  we have also that 
\begin{equation}\label{eqn:easy} \tr (A) (1-f_r^2) \leq  n \Lambda \chi_{B_{2r} \setminus B_r}(1-f_r^2)  \leq  \frac{n \Lambda |Ax|^2}{\lambda^2 r^2}(1-f_r^2).
\end{equation}
Now, starting from Equation \eqref{eqn:cutoff} and then adding the term $\sqrt{\ep} \tr (A) \int_{\Omega} | \psi|^2 (1-f_r^2)$ on both sides, we can use \eqref{eqn:easy}, \eqref{eqn:grad}   and \eqref{eqn:gradient} to obtain
$$  \int_{\Omega} \ep (1+\delta) |\nabla \psi|^2 f_r^2  + \left( f_r^2+ \frac {n\Lambda \sqrt{\ep}}{\lambda^2 r^2 }(1-f_r^2) + \frac{ 2 \ep }{\delta r^4 \lambda^2} \right) | Ax|^2 | \psi|^2 \, dx \geq  \sqrt{\ep} \tr (A) \int_{\Omega} | \psi|^2\, dx;$$
In particular, as long as $r^2\lambda^2 \geq n \Lambda \sqrt{\ep}$ we have $f_r^2+ \frac {n\Lambda \sqrt{\ep}}{\lambda^2 r^2 }(1-f_r^2)  \leq 1$ and then we optimize in $\delta$, choosing $\delta= \frac{\sqrt{2 \ep}}{r^2 \lambda}$, getting
$$ \left(1+ \frac{\sqrt{2 \ep} }{r^2 \lambda}\right) \cdot \int_{\Omega} \Bigl( \ep |\nabla \psi|^2  +  | Ax|^2 | \psi|^2  \Bigr)\, dx \geq \sqrt{\ep} \tr (A) \int_{\Omega} | \psi|^2\, dx,$$
which   is precisely \eqref{eqn:Omega}, with $C=\frac{ \sqrt{2}}{r^2 \lambda}$.
\end{proof}
%{ \color{blue}
%$\lambda_1 = \tr A /\ep$ while $\lambda_2 = \tr A / \ep + 2\lambda_{min} /\ep$. In particular if $ \psi =  \sum_{i} h_i \psi_i$ with $\psi_i$ eigenfunctions, $h_1 >0$, we have
%$$ E(\psi) = \sum \lambda_i h_i^2, \qquad \| \psi - \psi_1\|^2_2  = (h_1-1)^2 + \sum_{i>1} h_i^2 $$
%$$ 1 = \sum h_i^2, \qquad \sum_{i>1} h_i^2 \leq \frac \ep{2 \lambda_{min}} ( E(\psi)- \lambda_1 ) $$
%$$ \| \psi - \psi_1\|^2_2 \leq (h_1^2-1)^2 + \sum_{i>1} h_i^2 \leq \frac \ep{\lambda_{min}} ( E(\psi) - \lambda_1 ) $$}
%

\begin{lemma}\label{lem:delta}
Let us consider a continuous potential $V \geq 0$ such that  $V(0)=0$ and $V$ is twice differentiable in $0$. Then there exist a sequence $\psi_{\ep}$ such that $|\psi_{\ep}|^2 \weakto \delta_0$ and 
\begin{equation}\label{eqn:limsup} \limsup_{\ep \to 0}  \int_{\R^d} \Bigl( \ep^{1/2} \frac{ | \nabla \psi_{\ep}|^2 }2 +  \ep^{-1/2} V(x) | \psi_{\ep}|^2  \Bigr)\, dx \leq \frac 12 \tr \left( \sqrt{ D^2V(0)} \right).
\end{equation}

\end{lemma}
\begin{proof} % and let us denote $x=(x_1, x_2, \ldots, x_n)$ and $y=(x_{n+1} , \ldots, x_d)$ and let $\Lambda$ be the diagonal matrix $n \times n$ with eigenvalues $\lambda_1 , \ldots, \lambda_n$.
First, for a positive definite symmetric matrix $A$, let us consider consider the function
$$ f_{\ep} = \left(e^{ -  \frac{x^t A x }{ 2\sqrt{\ep}} } - e^{-N} \right)_+.$$
Then let $\Omega_{\ep}= \{ f_{\ep} \neq 0 \} =\{ x^t  A x   < 2N\sqrt{\ep} \}$: we claim that
\begin{equation}\label{eqn:fepN}\frac{ \int_{\R^d} \ep | \nabla f_{\ep}|^2 +   |Ax|^2 |f_{\ep}|^2 \, dx   }{ \int_{\R^d} |f_{\ep}|^2 \, dx  } =   \sqrt{\ep} \tr (A)   \cdot h(N),\end{equation}
for some universal function $h(N)$ such that $h(N) \to 1$ as $N \to \infty$. We postpone to the end the proof of the claim.

Now, since $V$ has an absolute minimum at $x=0$ we have that $D^2V(0)$ is a positive semidefinite symmetric matrix%; in particular $D^2V (0 )=\Lambda^2$ where $\Lambda$ is a  positive semidefinite symmetric matrix as well
.

Since $V$ is twice differentiable at $0$, for $\eta>0$ there exists an increasing continuous function $\delta(\eta)>0$ with $\delta(\eta)\to0$ for $\eta\to0$ such that
\[
2V(x) \leq x^t D^2V(0) x + \eta |x|^2
\leq \abs{(\Lambda + \sqrt{\eta}\Id)x}^2  \quad \forall x \in B_{\delta_n}(0).
\]
Since the function $\frac{\sqrt\eta}{-\log\eta}$ is continuous and strictly increasing, for every $\eps$ sufficiently small there exists a unique $\eta=\eta(\eps)$ such that $\sqrt\eps = \frac{\sqrt\eta\delta(\eta)^2}{-2\log\eta}$.
Now, let $A=\sqrt{D^2V(0)} + \sqrt\eta\Id$, $N=-\log\eta$ and $\psi_{\ep}= f_{\ep}/ \sqrt{ \int |f_{\ep}|^2 }$. With this choice of parameters, we get $\Omega_{\ep} \subseteq  B_{\delta}(0)$ (using that $A \geq \sqrt\eta\Id$). Since $\psi_{\ep}$ is also normalized, we get that $|\psi_{\ep}|^2 \weakto \delta_0$.
Moreover $2V(x) \leq |Ax|^2$ in $\Omega_{\ep}$ and hence,  we can say
\begin{align*} 
 \int_{\R^d} \Bigl( \ep^{1/2} \frac{ | \nabla \psi_{\ep}|^2}2 + \ep^{-1/2} V(x)   |\psi_{\ep}|^2  \Bigr) \, dx  & \leq \frac 1{2\sqrt{\ep }} \int_{\R^d}\Bigl( \ep | \nabla \psi_{\ep}|^2 + |Ax|^2|\psi_{\ep}|^2\Bigr) \, dx  \\ & = \frac 12 \tr (A) \cdot h(-\log{\eta}) \\ 
& = \frac 12 \left( \tr \Bigl(\sqrt{D^2V(0)}\Bigr)+d \cdot \sqrt{\eta}  \right) \cdot h(-\log{\eta}).
\end{align*}
Letting $\eps \to 0$, hence $\eta(\eps) \to 0$, we deduce \eqref{eqn:limsup}. 

We now prove the claim \eqref{eqn:fepN}. First of all
$$ |\nabla f_{\ep}|^2 =\frac 1{\ep}  |Ax|^2  e^{ -  \frac{ x^t A x }{ \sqrt{\ep}} } \chi_{\Omega_{\ep}},$$
$$ |f_{\ep}|^2 = \chi_{\Omega_{\ep}}(e^{ -  \frac{x^t A x}{ \sqrt{\ep}} } - 2 e^{ -  \frac{ x^t A x } {2\sqrt{\ep}} }e^{-N} + e^{-2N}).$$
We  make the change of variable $z=\frac{\sqrt{A}}{\sqrt[4]{\ep}} x$ which maps $\Omega_{\ep})$ in $B_{\sqrt{2N}}(0)=B$. We obtain
\begin{align*}
\frac{\det \sqrt A}{\eps^{1/2}} \int_{\Omega_{\ep}} \ep | \nabla f_{\ep}|^2 +   |Ax|^2 |f_{\ep}|^2 \, dx &=\int_{B} \sqrt{\ep} z^t A z \cdot ( 2e^{-|z|^2} - 2e^{-\frac{|z|^2}2 -N} + e^{-2N} ) \, dz\\
%&=  \sqrt{\ep} \int_{B} \Bigl( \sum_i \lambda_i |z_i|^2  \Bigr) \cdot ( 2e^{-z^2} - 2e^{-\frac{z^2}2 -N} + e^{-2N} ) \, dz \\
& = \sqrt{\ep} \cdot \frac {\tr (A)}{d} \int_B |z|^2 (2e^{-|z|^2} -2e^{-\frac{|z|^2}2-N}  + e^{-2N})\, dz,
\end{align*}
where in the second step we used that in $z^tAz=\sum_{i,j}A_{ij}z_jz_i$ every term with $i\neq j$ integrates to $0$ (when multiplied by the function in brackets) because the function is antisymmetric and every term $z_i^2$ gives the same integral as $\abs{z}^2/d$ by symmetry.
Performing a similar calculation for $\int | f_{\ep}|^2 $ we find that

$$\frac{ \int_{\Omega_{\ep}} \ep | \nabla f_{\ep}|^2 + \frac 1{\ep}  |Ax|^2 |f_{\ep}|^2 \, dx   }{ \int_{\Omega_{\ep}} |f_{\ep}|^2 \, dx  } =  \frac {\tr (A)}{d} \frac {\int_B z^2 (2e^{-z^2} -2e^{-\frac{z^2}2-N}  + e^{-2N})\, dz}{\int_B  (e^{-z^2} -2e^{-\frac{z^2}2-N}  + e^{-2N})\, dz}.$$
Since $\int_{\R^d} 2z^2e^{-z^2}\, dz = d \int_{\R^d} e^{-z^2}$ it is easy to see that $h(N):=  \frac {\int_B z^2 (2e^{-z^2} -2e^{-\frac{z^2}2-N}  + e^{-2N})\, dz}{d\int_B  (e^{-z^2} -2e^{-\frac{z^2}2-N}  + e^{-2N})\, dz}$ is such that $h(N) \to 1$ as $N\to \infty$.

%Now we want to compute the integrals: denoting with $\omega_d$ the measure of the $d$-dimensional ball of radius $1$ we have
%
%$$A_{2N}=\int_{B} e^{-z^2} \,dz= \int_0^{\sqrt{2N}} d\omega r^{d-1}e^{-r^2} \, dr$$
%$$ \int_{B} e^{-z^2/2}\, dz=  \int_0^{\sqrt{2N}} d\omega r^{d-1}e^{-r^2/2} \, dr = 2^{d/2} A_{4N}$$
%$$B_{2N}= \int_{B} z^2 e^{-z^2} \, dz= -\frac 12 \int_B z \cdot \nabla ( e^{-z^2})  = \frac d2 A_{2N} + \frac{d \omega_d (2N)^{d/2} }2e^{-2N}$$
%$$\int_{B} z^2 e^{-z^2/2} \, dz = - \int_B z \cdot \nabla ( e^{-z^2/2}) = d 2^{d/2} A_{4N}+ d \omega_d (2N)^{d/2} e^{-N}$$

\end{proof}

\begin{proof}[Proof of \autoref{thm:nonfixed-marginals}]
Starting from point (a) first of all we have that $| \psi_{\ep}|^2 \weakto \gamma$ for some $\gamma \in \mathcal{M}_+(\R^n)$ with $\gamma ( \R^d) \leq 1$. Then, assuming $\liminf_{\ep \to 0} E_{\ep}(\psi_{\ep})=C < \infty$, we have that for every $\delta>0$, as $\{V > \delta\}$ is an open set
$$ \gamma ( \{ V > \delta\}) \leq \liminf_{\ep \to 0}  \int_{ \{V > \delta\}} | \psi_{\ep}|^2 \leq \liminf_{\ep \to 0}\frac {\sqrt{\ep}} \delta  \cdot \frac 1{\sqrt{\ep}}\int_{ \R^d} V(x)| \psi_{\ep}^2| \leq \lim_{\ep \to 0} \frac {\ep} \delta C =0.$$
By the arbitrariness of $\delta$ we conclude that $\gamma ( \{V>0\} ) =0$. %Moreover, since $\{ V < \delta\} \subset B_R$ for some $R$ and $\delta$ small enough we get also that $\int_{B_R} | \psi_{\ep}|^2 \geq 1- \frac {2\ep}\delta C$; in particular the measures $|\psi_{\ep}|^2$ are tight and this proves that $\gamma$ is a probability measure.

Now we define the energy density $\mu_{\ep} = \left( \ep^{1/2} \frac{| \nabla \psi_{\ep}|^2}2 +\ep^{-1/2} V(x) |\psi_{\ep}|^2\right) \mathcal{L}^d$; by hypotesis we know that the sequence $\mu_{\ep}$ is uniformly bounded and in particular, up to subsequences, there exists a weak limit $\mu$. By the semicontinuity of the mass, it will be sufficient to prove that the Radon-Nikodym derivative of $\mu$ with respect to $\gamma$ is greater or equal than $\frac 12 \tr ( \sqrt {D^2V} )$. In particular, using a result by \cite{morse} (see Theorem 5.3 and Theorem 5.7 in \cite{bliedtnerloeb} for a clearer explanation) we have that
$$ \frac { d\mu}{d \gamma} (v_0)= \lim_{ r \to 0} \frac { \mu (C_r)}{\gamma(C_r)},$$
where $C_r$ are cylinders centered in $v_0$, that is a set of the form $B^X_r(x_0) \times B^Y_r(y_0)$ for an orthogonal splitting $\R^d= X \times Y$ and $v_0=(x_0,y_0)$. It is crucial to notice that any such cylinder is convex and $B_r(v_0) \subseteq C_r \subseteq B_{\sqrt{2}r}(v_0)$: in particular $C_r \in \mathscr{S}_{\sqrt{2}}(v_0)$ in the notation of \cite{bliedtnerloeb}.

Thus we are done if for every $v_0 \in \{V=0\}$ and every $0<\delta<1$, for $r>0$ small enough there is an open cylinder $C_r$ such that
 \begin{align}\label{eqn:ts-rad-nic}
 \mu(\overline{C_r}) & \geq \liminf_{\ep \to 0} \mu_{\ep}( C_r)  = \liminf_{\ep \to 0} \left\{ \int_{C_r} \ep^{1/2} \frac{| \nabla \psi_{\ep}|^2}2 +\ep^{-1/2} V(x) |\psi_{\ep}|^2 \, dx   \right\}  \\
 & \geq \frac{\sqrt{1-\delta}}2  \tr \left( \sqrt{D^2V(v_0)} \right) \gamma(C_r).
 \end{align}
 
In order to do this let us suppose that $D^2V \neq 0$ otherwise the statement is trivial. In that case, denote $A= \sqrt{ D^2 V (v_0)}$ and make an orthonormal change of coordinates such that $D^2V(v_0)$ becomes diagonal with the eigenvalues in decreasing order. Then, letting $n={\rm rk} (D^2V (v_0))$ we consider $X$ the span of the first $n$ eigenvectors and $Y$ the remaining ones. In particular we have
%\begin{equation*}%\label{eqn:coerciveX}
$c \Id \leq D^2V (v_0)|_X \leq C \Id $ { for some } $C>c>0$ (depending on $v_0$)
%\end{equation*}
and $D^2V|_Y=0$. Moreover, we have that $\nabla_XV(x_0,y_0)=0$ and $\nabla_X\nabla_XV(x_0,y_0)=D^2V(v_0)|_X$ is invertible, therefore by the implicit function theorem there exists $r>0$ and a $C^1$ function $m:B^y_r(y_0)\to B^x_r(x_0)$ such that $\nabla_XV(m(y),y)=0$ for every $y\in B^y_r(y_0)$. The gradient of the implicit function is
\[
\nabla_Ym(y)=-(\nabla_{XX}V(m(y),y))^{-1} \nabla_{XY}V(m(y),y),
\]
which vanishes for $y=y_0$ since $\nabla_{XY}V(v_0)=0$. Therefore, up to possibly reducing $r$, we have that $|m(y)-x_0|\leq\delta|y-y_0|$ for every $y\in B^y_r(y_0)$. By the continuity of $D^2V$, up to reducing $r$ a second time, we have also $D^2V(x,y)|_X\geq(1-\delta)D^2V(x_0,y_0)|_X$ for every $(x,y)\in B^x_r(x_0)\times B^y_r(y_0)$. Therefore, for every $(x,y)\in B^x_r(x_0)\times B^y_r(y_0)$, we consider the second order Taylor expansion of $V$ at $(y,m(y))$ (recall that $V \geq 0$ and that $\nabla_X V$ vanishes at this point) to find% there is $\xi\in B^x_r(x_0)$ such that
\[
\begin{split}
V(x,y)% &= V(m(y),y) + \nabla_XV(m(y),y) + \frac12 (x-m(y))^T D^2V(\xi,y)|_X(x-m(y)) \\
&\geq \frac{1-\delta}{2}(x-m(y))^T D^2V(v_0)|_X(x-m(y))
= \frac{1-\delta}{2} |A(x-m(y))|^2
\end{split}
\]
where we denoted $Ax=A(x,0)$; since $A$ is non degenerate only in the $x$ directions in particular we have $|Ax| \geq c|x|$ for some $c>0$ and $\tr_x (A)=\tr (A)$. 

Since $B^x_{r(1-\delta)}(m(y)) \subseteq B^x_r(x_0)$, we can use \autoref{lem:main} in order to say that there exists a constant $C$ such that, for every $y \in B_r^y(y_0)$:
\begin{align*} \int_{B^x_r(x_0)} \ep \frac{| \nabla_x \psi_{\ep} (x,y)|^2}2 +   V(x,y) | \psi_{\ep}(x,y)|^2 \, dx & \geq  \frac 12 \int_{B^x_r(x_0)} \ep | \nabla_x \psi_{\ep}|^2 + (1-\delta)    |A(x-m(y))|^2 | \psi_{\ep}|^2 \, dx  \\&  \geq \sqrt{\ep}  \frac{ \sqrt{1-\delta} }{1+C \sqrt{\ep} } \cdot \frac{\tr (A)}2 \int_{B^x_r(x_0)} | \psi_{\ep}(x,y)|^2 \, dx . 
\end{align*}
Integrating with respect to $y$ and then passing to the limit as $\ep \to 0$ we obtain \eqref{eqn:ts-rad-nic}.

In order to prove part (b) we use \autoref{lem:delta} in case $\gamma= \delta_{v_0}$ and in the rest of the cases we argue by convexity: let us consider $\psi_{\ep}^v(x):=\psi_\eps(x-v)$ a recovery sequence for $\delta_v$ and let us define $f_{\ep}^v = | \psi_{\ep}^v|^2$.
We notice that 
$$ \frac 1{\sqrt{\ep}}\int_{\R^d} \Bigl( \ep \frac{| \nabla \psi_{\ep}^v|^2}2 +  V(x) | \psi_{\ep}^v|^2\Bigr) \, dx = \frac 1{\sqrt{\ep}}\int_{\R^d} \Bigl(\ep \frac { | \nabla f_{\ep}^v|^2}{8f_{\ep}^v}  +  V(x) f_{\ep}^v \Bigr)\, dx$$
denoting by $H_{\ep} (f_{\ep}^v)$ the functional in the right-hand side, we observe that $H_{\ep}$ is convex.  Let us consider moreover the set $V_{\ep} \subseteq  \{ V=0\}$ defined as 
$$ V_{\ep} = \left\{ v \in \{ V=0 \} \; : \; H_{\ep}(f_{\ep}^v) \leq \frac 12 \tr \Bigl(\sqrt{ D^2V}(v) \Bigr) + 1 \right\}.$$
By construction we have that $	\limsup_{\ep} H_{\ep}(f_{\ep}^v) \leq \frac 12 \tr \Bigl(\sqrt{ D^2V}(v) \Bigr)$ for every $v$; in particular this implies that $\chi_{V_{\ep}} \to \chi_{ \{V=0\} }$ pointwise.

Now given a measure $\gamma \in \mathcal{P}(\R^d)$ with $\gamma(\R^d)$ supported in $\{V=0\}$, by dominated convergence we get $\gamma(V_{\ep}) \to 1$. In particular for $\ep$ sufficiently small $\gamma(V_{\ep}) \geq \frac 12$  and for such $\eps$ we define the functions
$$ f^{\gamma}_{\ep} (x) = \frac 1{ \gamma(V_{\ep})}  \int_{ V_{\ep}} f^{v}_{\ep} (x) \, d \gamma.$$
Consider now $\psi_{\ep} (x)= \sqrt{ f^{\gamma}_{\ep} (x) }$. By construction we have $\int | \psi_{\ep}|^2 =1$; moreover, by testing with a continuous function we see that $|\psi_{\ep}|^2 \weakto \gamma$ and thanks to the convexity of $H_{\ep}$ we get that
\begin{align*}
\limsup_{\ep \to 0} E_{\ep} (\psi_{\ep}) &= \limsup_{\ep \to 0} H_{\ep} (f^{\gamma}_{\ep}) \leq\limsup_{\ep \to 0} \frac 1{\gamma(V_{\ep})} \int_{V_{\ep}} H_{\ep}( f^v_{\ep}) \, d \gamma \\
&= \limsup_{\ep \to 0} \int_{\R^d} \frac{\chi_{V_{\ep}}}{\gamma(V_{\ep})}  E_{\ep}( \psi^v_{\ep}) \, d \gamma \end{align*}
Now let us consider the functions $g_{\ep} (v) = \frac{\chi_{V_{\ep}}(v)}{\gamma(V_{\ep})} E_{\ep}( \psi^v_{\ep}) $ and $g(v) = \tr\bigl(\sqrt{ D^2V(v) } \bigr)$; if $g$ is not $\gamma$-integrable there is nothing to prove; so we can assume $g \in L^1(\gamma)$. By definition of $V_{\ep}$, for $\eps$ sufficiently small we have
$$  \frac{\chi_{V_{\ep}}(v)}{\gamma(V_{\ep})} E_{\ep}( \psi^v_{\ep}) \leq  \tr\Bigl(\sqrt{ D^2V(v) } \Bigr) + 1 $$
and we notice that the right-hand side is  $\gamma$-integrable, otherwise there is nothing to prove. 
By Fatou lemma, using also $\gamma(V_{\ep}) \to 1$, $\chi_{V_{\ep}} \to 1$, we have 
\[
\limsup_{\ep \to 0} E_{\ep}(\psi_{\ep}) \leq \limsup_{\ep \to 0} \int_{\R^d} \frac{\chi_{V_{\ep}}}{\gamma(V_{\ep})}  E_{\ep}( \psi^v_{\ep}) \, d \gamma \leq  \int_{\R^d} \limsup_{\ep \to 0} \frac{\chi_{V_{\ep}}}{\gamma(V_{\ep})} E_{\ep}( \psi^v_{\ep}) \, d \gamma \leq \frac 12 \int_{\R^d} \tr ( \sqrt{ D^2V} ) \, d \gamma. \qedhere
\]
\end{proof}

\section{Rectangular truncation of Gaussians}\label{sec:gaussians}

As described in the introduction, Gaussian densities are asymptotically optimal for the energy of a single delta, therefore, to approach the general case, it is natural to construct a recovery sequence which is a superposition of Gaussian kernels. However, since we do not have global estimates on the potential $V$, only short-range interactions can be allowed, hence we need the kernels to have compact support. To achieve this, we truncate the densities with a slowly growing parameter $N$. The construction we present in \autoref{sec:construction} of the recovery sequence for \autoref{thm:rec} is then given by a convolution of $x$-dependent suitably rescaled and truncated Gaussian-like kernels, whose major axis is parallel at each point to $\graph T$.

The requirements discussed above leads us to introduce some kernels which resemble Gaussian densities but have compact support, built in a suitable way that guarantees good properties for both the kinetic and potential energy.

Let $M$ be a positive definite symmetric $2\times 2$ matrix with eigenvalues $a\geq b>0$ and let $w$ and $z$ be the orthonormal coordinates in the direction of the corresponding eigenvectors. 
For $N\in\setR$ we define the unnormalized Gaussian and unnormalized truncated Gaussian
\begin{align*}\label{eq:tilde-Gamma}
\tilde\Gamma_{M,\infty}(\x)
&= e^{-\x^TM\x}
= e^{-aw^2} e^{-bz^2}, \\
\tilde\Gamma_{M,N}(\x) &=
	\left(e^{-aw^2/2} - e^{-N/2}\right)_+^2 \left(e^{-bz^2/2} - e^{-N/2}\right)_+^2 .
\end{align*}

 In the definition of these truncated kernels, we adopt a rectangular truncation, which leads to kernels of product form; we call it rectangular because the support of the truncated Gaussian is a tilted rectangle, aligned with the eigenvectors of the matrix. This structure is crucial to do the computations in \autoref{sec:marg-prod}. Moreover, instead of the obvious truncation, we have to truncate the square root of the Gaussian, and then put a square outside, to guarantee the finiteness of the kinetic energy. 

In the proof of the main theorem we will apply the results of this section to the matrix
\[
M = M_{\eps,\beta}(x) = \frac{A(x)}{\eps^{1/2}} + \frac{I}{\beta}
\]
where $A(x)=\sqrt{\nabla^2V\bigl(x,T(x)\bigr)}$, hence $a=q(x)/\eps^{1/2}+1/\beta$ and $b=1/\beta$, with $q(x)$ being the non-zero eigenvalue of $A(x)$.
The parameter $\beta \gg \eps^{1/2}$ is required to ensure that the resulting Gaussian has compact support, since $\nabla^2V\bigl(x,T(x)\bigr)$ is a singular matrix, hence the natural choice of $M=A(x)/\eps^{1/2}$ wouldn't work properly. We will discuss more about $\beta$ when we present all the parameters in \autoref{sec:parameters}.

%\old{Notice that a natural choice of $M$ may seem the one that does not involve $\beta$, but since $D^2V$ is a degenerate matrix in order to ensure that the Gaussian has compact support, we need to add a multiple of the identity to it, weighted with the parameter $\beta \gg \eps^{1/2}$.}

%As seen in the introduction, the Gaussians are asymptotically optimal for the energy of a single delta. Therefore, it would be natural to approach the general case by constructing a recovery sequence which is a superposition of Gaussians. However, since we do not have global estimates on $V$, we need to localize these Gaussians. We adopt a rectangular truncation, which leads to kernels of product form. This structure is crucial to do the computations in \autoref{sec:marg-prod}. Moreover, instead of the obvious truncation, we necessarily have to use a more involved one to guarantee that the kinetic energy is finite. In order to ensure that the Gaussian has compact support, we need to add a small multiple of the identity to the matrix defining it.

Define now the one dimensional integrals
\begin{align*}
G_{\alpha,\infty} &= \int_\setR e^{-\alpha t^2} \d t = \frac{\sqrt\pi}{\sqrt\alpha}, \\
G_{\alpha,N} &= \int_\setR \left(e^{-\alpha t^2/2}-e^{-N/2}\right)_+^2 \d t
= \int_{-\sqrt{N/\alpha}}^{\sqrt{N/\alpha}}
	\left(e^{-\alpha t^2/2}-e^{-N/2}\right)^2 \d t,
\end{align*}
and the integrals of the two dimensional densities
\begin{align*}
G_{M,\infty} &= \int_{\setR^2} \tilde\Gamma_{M,\infty}(\x) \d \x
	= G_{a,\infty} G_{b,\infty}
	= \frac{\pi}{\sqrt{\det M}}, \\
G_{M,N} &= \int_{\setR^2} \tilde\Gamma_{M,N}(\x) \d \x = G_{a,N} G_{b,N} .
\end{align*}
We can now introduce the normalized Gaussian and truncated Gaussian, which are probability densities, given by 
\begin{equation}\label{eq:gamma-truncated}
\Gamma_{M,N}(\x) = \frac{\tilde\Gamma_{M,N}(\x)}{G_{M,N}},
\qquad N\in\setR\cup\{\infty\}.
%\Gamma_{M,N}(\x)
%&= \frac{\left(e^{-aw^2/2} - e^{-N/2}\right)_+^2} {G_{a,N}}
%	\cdot \frac{\left(e^{-bz^2/2} - e^{-N/2}\right)_+^2} {G_{b,N}} .
\end{equation}
Finally, let $\eta_{M,\infty}$ and $\eta_{M,N}$ be the first marginal of $\Gamma_{M,\infty}$ and $\Gamma_{M,N}$ respectively.

We present here two lemmas, but we postpone their proof to \autoref{sec:proof-lemmas}.

\begin{lemma}\label{lem:norm-diff-gauss}
With the definition above, there exists a constant $C>0$ such that for every $N \geq 3$ we have
\begin{align}
G_{M,N} < G_{M,\infty} &< G_{M,N} + Ce^{-N/2}G_{M,\infty} ,
	\label{eq:GN-infty} \\
\norm{\Gamma_{M,N}-\Gamma_{M,\infty}}_\infty &\leq C\sqrt{\det M} e^{-N/2} ,
	\label{eq:Linfty-norm-trunc} \\
\norm{\Gamma_{M,N}-\Gamma_{M,\infty}}_1 &\leq C N e^{-N/2} ,
	\label{eq:L1-norm-trunc} \\
\norm{\eta_{M,N}-\eta_{M,\infty}}_\infty &\leq C\sqrt a \sqrt{N} e^{-N/2} .
	\label{eq:Linfty-norm-trunc-marg}
\end{align}
\end{lemma}

Let us now compute the kinetic and potential energies of a Gaussian.
A direct computation shows that, if $Be_w=fe_w$ and $Be_z=ge_z$, then
\[
\int_{\setR^2} \abs{B\x}^2 \Gamma_{M,\infty}(\x) \d\x
= \int_{\setR^2} (f^2w^2+g^2z^2) \frac{\sqrt{ab}}{\pi} e^{-aw^2-bz^2} \d w \d w
= \frac12 \left(\frac{f^2}{a} + \frac{g^2}{b}\right) .
\]
Therefore, using $\abs*{\nabla\sqrt{\Gamma_{M,\infty}}}^2
= \frac12\abs*{\nabla[\log\Gamma_{M,\infty}]}^2 \Gamma_{M,\infty}
= \abs{M\x}^2 \Gamma_{M,\infty}(\x)$ and the previous identity with $B=M$ we get
\begin{equation}
\label{eqn:daqualcheparte}
\begin{split}
\KE(\Gamma_{M,\infty})
%&= \frac12 \int_{\setR^2} \abs*{\nabla\sqrt{\Gamma_{M,\infty}(\x)}}^2 \d\x
%= \frac12 \int_{\setR^2} \abs*{\frac{\nabla\Gamma_{M,\infty}(\x)}
%	{2\Gamma_{M,\infty}(\x)} \sqrt{\Gamma_{M,\infty}(\x)}}^2 \d\x \\
%&= \frac12 \int_{\setR^2} \abs*{\frac12\nabla[\log\Gamma_{M,\infty}(\x)]}^2
%	\Gamma_{M,\infty}(\x) \d x \\
&= \frac12 \int_{\setR^2} \abs{M\x}^2 \Gamma_{M,\infty}(\x) \d x
= \frac{\tr M}{4} .
\end{split}
\end{equation}

The following lemma compares the potential energy associated to the quadratic potential induced by the matrix $M$ and the kinetic energy of the Gaussian $\Gamma_{M,\infty}$ and the truncated Gaussian $\Gamma_{M,N}$.

\begin{lemma}\label{lem:diff-energies}
There is a universal constant $C>0$ such that
%\begin{equation}\label{eq:truncated-quadratic-energy_}\color{gray}
%\abs*{\int_{\setR^2} \abs{M\x}^2 \Gamma_{M,N}(\x) \d\x
%	- \int_{\setR^2} \abs{M\x}^2 \Gamma_{M,\infty}(\x) \d\x}
%\leq \norm{M}^2 N e^{-N/2} ,
%\end{equation}
\begin{equation}\label{eq:truncated-quadratic-energy}
\int_{\setR^2} \abs{M\x}^2 \abs{\Gamma_{M,N}(\x)-\Gamma_{M,\infty}(\x)} \d\x
\leq C\tr(M) N e^{-N/2} ,
\end{equation}
\begin{equation}\label{eq:truncated-kinetic-energy}
\abs{\KE(\Gamma_{M,N}) - \KE(\Gamma_{M,\infty})}
\leq C \KE(\Gamma_{M,\infty}) e^{-N/2}
= C \tr(M) e^{-N/2} .
\end{equation}
%{\color{gray}
%\begin{equation}
%\abs*{ \abs*{\nabla\sqrt{\Gamma_{M,N}(\x)}}^2
%	- \abs*{\nabla\sqrt{\Gamma_{M,\infty}(\x)}}^2 }
%\leq \dots
%\end{equation}
%\begin{equation}%\label{eq:norm-gradient-sqrt-gamma}
%\abs*{\nabla\sqrt{\Gamma_{M,N}(\x)}}^2
%\leq \frac{G_{M,\infty}}{G_{M,N}} \abs{M_{\eps,\beta}\x}^2
%	\Gamma_{M,\infty}(\x) \bm{1}_{\supp\Gamma_{M,N}}(\x)
%\end{equation}
%} % end color
%and
%\begin{equation*}
%\norm{(h^2_\eps)'(z)}_\infty \leq \dots ,\qquad
%\norm{(h^2_\eps)''(z)}_\infty \leq \dots .
%\end{equation*}
\end{lemma}

\section{Construction of the recovery sequence} \label{sec:rec}

\subsection{Structure of optimal plans, maps and potentials in one dimension}
%{Assumptions on $\rho$ and boundedness of the second derivative of the Kantorovich potential $u$}
\label{sec:rho}

In this section we comment on the assumption on the boundedness of the second derivative of the Kantorovich potential $u$ in \autoref{thm:rec}
and we provide a class of $\rho$ for which this assumption is satisfied and the main \autoref{thm:rec} is applicable: we will use the structural results for the $1$D Coulomb multimarginal optimal transport problem contained in \cite{CDD}, which allows us to transfer the regularity of $\rho$ to information on $V$ and the optimal maps.

Let $\rho\in C^1_{\rm loc}(\setR)$ be a strictly positive probability density. Without loss of generality, by translating we may assume that the origin $0$ is the median of the probability $\rho$, i.e.\ $\rho\bigl((-\infty,0]\bigr)=\rho\bigl([0,\infty)\bigr)=1/2$. By the result in \cite{CDD}, the unique optimal plan is induced by a map $T$ from $\rho$ to itself which can be written explicitly in terms of the repartition function of $\rho$, which is $C^2$ with strictly positive first derivative. As a consequence, we have that $T\in C^2(\setR\setminus\{0\})$, $T$ is increasing in $\setR_-$ and $\setR_+$, $T(\setR_-)=\setR_+$, $T(\setR_+)=\setR_-$ and
\begin{align*}
\lim_{x\to-\infty} T(x) &= 0, & \lim_{x\to0^-} T(x) &= \infty, &
\lim_{x\to0^+} T(x) &= -\infty, & \lim_{x\to+\infty} T(x) &= 0.
\end{align*}

The Kantorovich potential $u:\setR\to\setR$ satisfies by definition $u(x) + u(y) \leq |x-y|^{-1}$ with equality if and only if $y= T(x)$; hence it is determined up to a constant by
\begin{equation}\label{eqn:uprimo}
u'(x) = - \frac{\sign(x)}{\bigl(T(x)-x\bigr)^2} .
\end{equation}
Let $V$ be as in \eqref{eq:V}, which in $d=1$ and $N=2$ reads as
\begin{equation}\label{eqn:V}
V(x,y) = |x-y|^{-1} - u(x) - u(y) \qquad \mbox{for }x,y\in \R.
\end{equation}
Assuming $x<0$ and $y>0$ (the other case is analogous), then we can compute the gradient of the potential
\[
\nabla V(x,y) = \begin{pmatrix}
\frac{1}{(y-x)^2} - u'(x) \\
-\frac{1}{(y-x)^2} - u'(y)
\end{pmatrix} ;
\]
notice that $\nabla V\bigl(x,T(x)\bigr)=0$. Differentiating again we get the Hessian
\begin{equation}\label{eqn:hess-v}
\nabla^2V(x,y)
= \frac2{(y-x)^3}\begin{pmatrix}1 & -1 \\ -1 & 1\end{pmatrix}
	- \begin{pmatrix} u''(x) & 0 \\ 0 & u''(y) \end{pmatrix}.
\end{equation}
In the proof of \autoref{thm:rec} we are interested in its behavior in a neighborhood of the graph of $T$, i.e.\ when $y=T(x)$.
%As $x\to0^-$, we have $y\to\infty$, hence
%\[
%\lim_{x\to0^-} \nabla^2V\bigl(x,T(x)\bigr)
%= - \begin{pmatrix} u''(0^-) & 0 \\ 0 & u''(\infty) \end{pmatrix}.
%\]
%Let's compute its limits.
%\[
%\begin{split}
%\lim_{x\to0^-} u''(x)
%&= \lim_{x\to0^-} \frac{\d}{\d x}\left[\frac1{\bigl(T(x)-x\bigr)^2}\right] \\
%&= 2\lim_{x\to0^-} \frac1{\bigl(T(x)-x\bigr)^3}
%	-2\lim_{x\to0^-} \frac{T'(x)}{\bigl(T(x)-x\bigr)^3} \\
%&= 0 - 2 \lim_{x\to0^-} \frac{T'(x)}{T(x)^3} .
%\end{split}
%\]

%We say that our probability $\rho$ has \emph{heavy tails} if $\rho(x)\geq C/x^3$ for $x\to\pm\infty$; we say that our probability $\rho$ has \emph{controlled tails} if $\rho(x)\leq Cx^3$ for $x\to\pm\infty$.
%In such case, we have both
%\[
%\lim_{x\to0^-} \frac{T'(x)}{T(x)^3}
%= \lim_{x\to0^-} \frac{\rho(x)}{T(x)^3\rho\bigl(T(x)\bigr)} < \infty
%\]
%and
%\[
%\lim_{x\to-\infty} \frac{T'(x)}{T(x)^3}
%= \lim_{x\to-\infty} \frac{\rho(x)}{T(x)^3\rho\bigl(T(x)\bigr)} < \infty.
%\]
%This implies that $u''$ is bounded on $(-\infty,0)$, therefore $\nabla^2V$ stays bounded outside of a neighborhood of the diagonal, and in particular near the graph of $T$.
The aim of the following lemma is twofold: first, we provide an assumption on the tails of $\rho$ which is sufficient to obtain boundedness of the second derivatives of $u$. Secondly, we show how the assumption on the boundedness of $u''$ will be used in the rest of the paper: through the computation of the Hessian of $V$, it will provide a sufficient condition to control the growth of $V$ with a (uniform) parabola around the graph of $T$. We expect \autoref{thm:rec} to hold even without the uniform growth condition \eqref{eqn:Vgrowth} below, but at the price of several technical complications that we don't address here.

\begin{lemma}\label{lem:good-rho}
Let $\rho\in C^1_{\rm loc}(\setR)$ be a strictly positive probability density, let $u$ be its Kantorovich potential and let $V$ as in \eqref{eqn:V}.
\begin{enumerate}[(i)]
\item If $u'' \in L^\infty(\setR)$, then there exists $\eps_0,C>0$ such that
\begin{equation}\label{eqn:Vgrowth}
V(x,y) \leq C \dist\bigl((x,y),\graph(T)\bigr)^2 \qquad
\forall (x,y) \text{ s.t.\ } \dist\bigl((x,y),\graph(T)\bigr)< \eps_0;
\end{equation}
\item if $\rho$ satisfies also \eqref{hp:tails}, then $u'' \in L^\infty(\setR)$ and $\nabla^2V$ is locally Lipschitz in a neighborhood of $\graph(T)$.
\end{enumerate}
\end{lemma}
%We are interested in its behavior in a neighborhood of the graph of $T$, i.e.\ when $y=T(x)$. The first entry is locally bounded in $(-\infty,0)$. 

\begin{proof}
Since $\rho \in L^1(\R)$, there exists $\delta>0$ such that every interval $I$ of length $<\delta$ has $\int_I \rho(x) \, dx <1/2$. Since $ \int_x^{T(x)} \rho(x')\, dx'=1/2$ for every $x <0$ and  $ \int_{T(x)}^x \rho(x')\, dx' =1/2$ for every $x >0$, we obtain that $|T(x) - x|>\delta$, namely that $\graph(T)$ has positive distance from the diagonal. The hessian of $V$ is
given by \eqref{eqn:hess-v} and, %\[
%\nabla^2V(x,y)
%= \frac2{(y-x)^3}\begin{pmatrix}1 & -1 \\ -1 & 1\end{pmatrix}
%	- \begin{pmatrix} u''(x) & 0 \\ 0 & u''(y) \end{pmatrix}.
%\]
in a neighborhood of $\graph(T)$, the first term is bounded because we are far from the diagonal $y=x$, whereas the second term is bounded everywhere since we have proved that $u''\in L^\infty(\setR)$. This implies that $V$ has controlled quadratic growth in a neighborhood of $\graph(T)$, namely that $V(x,y) \leq C \dist\bigl((x,y),\graph(T)\bigr)^2$ in $\R^2 \setminus \{|x-y|< \eps $ for any $\eps>0$.% Hence in particular it holds in a uniform neighborhood of $\graph T $, where t
The size $\eps$ of the neighborhood only depends on $\rho$ (and in particular on the positive distance between the graph of the optimal map $T$ and the diagonal $\{x=y\}$). This establishes (i).

%To prove (ii), we first notice that, by integrating, \eqref{hp:tails} implies
%\[
%\liminf_{\abs{x}\to\infty} \abs{x}^\alpha \rho(x) > 0 \qquad \text{and} \qquad
%\limsup_{\abs{x}\to\infty} \abs{x}^\beta \rho(x) < \infty .
%\]
Let us now prove (ii).
From the formula \eqref{eqn:uprimo}
%\[
%u'(x)=-\frac{\sign(x)}{\bigl(T(x)-x\bigr)^2}
%\]
and the fact that $T\in C^1_{\rm loc}(\setR\setminus\{0\})$ we deduce that $u\in C^2_{\rm loc}(\setR\setminus\{0\})$. We need to check that $u''(x)$ stays bounded as $x\to0^\pm$ and $x\to\pm\infty$.
For $x<0$, we have that
\[
u''(x) = \frac{\d}{\d x}\left[\frac1{\bigl(T(x)-x\bigr)^2}\right]
= 2\frac1{\bigl(T(x)-x\bigr)^3} - 2 \frac{T'(x)}{\bigl(T(x)-x\bigr)^3}.
\]
In both cases when $x\to-\infty$ or $x\to0^-$, the first fraction goes to $0$, therefore $u''(x)$ has the same asymptotic behavior as $T'(x)/\bigl(T(x)-x\bigr)^3$, let it be having a limit, being bounded or diverging. By the Monge-Amp\`ere equation, we have moreover that
\[
\frac{T'(x)}{\bigl(T(x)-x\bigr)^3}
= \frac{\rho(x)}{\bigl(T(x)-x\bigr)^3\rho\bigl(T(x)\bigr)},
\]
so our goal is to show that this fraction stays bounded as $x\to-\infty$ or $x\to0^-$.
The quantity is clearly non-negative, so  we need to prove that the $\limsup$ is finite. It is immediate to see that
\[
\limsup_{x\to0^-} \frac{T'(x)}{\bigl(T(x)-x\bigr)^3}
= \limsup_{x\to0^-} \frac{\rho(x)}{T(x)^3\rho\bigl(T(x)\bigr)}
= \limsup_{y\to\infty} \frac{\rho(0)}{y^3\rho(y)}
= \frac{\rho(0)}{\liminf\limits_{y\to\infty} y^3\rho(y)}
< \infty
\]
because of \eqref{hp:tails}. Let's now turn to studying the limit for $x\to-\infty$. Since $T(x)\to0^+$ for $x\to-\infty$, we have
\[
\limsup_{x\to-\infty} \frac{T'(x)}{\bigl(T(x)-x\bigr)^3}
= \limsup_{x\to-\infty} \frac{\rho(x)}{-x^3\rho\bigl(T(x)\bigr)}
= \limsup_{x\to-\infty} \frac{\rho(x)}{-x^3\rho(0)}
%= \frac{\limsup\limits_{x\to-\infty}{} \abs{x}^\beta \rho(x)}
%	{\rho(0)\liminf\limits_{x\to-\infty}{} \frac{\bigl(\abs{x}^{\beta/3}R(x)\bigr)^3}{\rho(0)^3}}
< \infty
\]
because $\rho\in L^\infty(\setR)$.
%Define
%\begin{align*}
%R(x) &= \int_{-\infty}^x \rho(x)\d x , &
%S(x) &= \int_0^x \rho(x)\d x.
%\end{align*}
%We have that $S(y)=\rho(0)y+o(y)$ for $y>0$ small, hence $S^{-1}(y)=y/\rho(0)+o(y)$ The optimal map $T$ is characterized by $S\bigl(T(x)\bigr)=R(x)$ for $x<0$, therefore $T(x)=S^{-1}\bigl(R(x)\bigr)=R(x)/\rho(0)+o\bigl(R(x)\bigr)$ for $x\to-\infty$.
%By integration, we have that
%\[
%\liminf_{x\to-\infty}{} \abs{x}^{\alpha-1}R(x) > 0,
%\]
%therefore
%\[
%\limsup_{x\to-\infty} \frac{\rho(x)}{T(x)^3\rho\bigl(T(x)\bigr)}
%= \limsup_{x\to-\infty} \frac{\rho(x)}{T(x)^3\rho(0)}
%= \frac{\limsup\limits_{x\to-\infty}{} \abs{x}^\beta \rho(x)}
%	{\rho(0)\liminf\limits_{x\to-\infty}{} \frac{\bigl(\abs{x}^{\beta/3}R(x)\bigr)^3}{\rho(0)^3}}
%< \infty
%\]
%because $\alpha-1\leq\beta/3$.
Exactly the same can be said for $x\to0^+$ and $x\to\infty$, therefore we have that $u\in C^2(\setR)$ and $u''\in L^\infty(\setR)$.

Notice moreover that by \eqref{eqn:uprimo} and the fact that $T\in C^2(\setR\setminus\{0\})$ it follows that $u\in C^3(\setR\setminus\{0\})$. Hence, denoting $D\subset\setR^2$ the diagonal, by \eqref{eqn:hess-v} we have that $\nabla^2V$ is locally Lipschitz in $(\setR\setminus\{0\})^2\setminus D$, which is a neighborhood of $\graph(T)$.
%Let's now estimate $\KE(\rho)$. First of all, from \eqref{hp:ab} we deduce that $2+2\beta-\alpha>1$. Let $M,C>0$ be such that $\rho'(x) < C \abs{x}^{-\beta-1}$ and $\rho(x) > C\abs{x}^{-\alpha}$ for $\abs{x}>M$. Then
%\[
%\begin{split}
%\KE(\rho)
%&= \int_M^M \frac{\rho'(x)^2}{\rho(x)} \d x
%	+ \int_{\abs{x} > M} \frac{\rho'(x)^2}{\rho(x)} \d x \\
%&\leq 
%%\int_M^M \frac{\rho'(x)^2}{\rho(x)} \d x
%%	+ C\int_{\abs{x} > M} \frac{\abs{x}^{-2\beta-2}}{\abs{x}^{-\alpha}} \d x \\
%%&= 
%\int_M^M \frac{\rho'(x)^2}{\rho(x)} \d x
%	+ C\int_{\abs{x} > M} \frac{1}{\abs{x}^{2\beta+2-\alpha}} \d x < \infty. \qedhere
%\end{split}
%\]
%because $2\beta+2-\alpha>1$.
\end{proof}
%If $T$ has heavy tails, then as we saw $\lim_{x\to0^-}u''(x)=0$, and moreover also
%\[
%\begin{split}
%\lim_{x\to-\infty} u''(x)
%&= 2\lim_{x\to-\infty} \frac1{\bigl(T(x)-x\bigr)^3}
%	-2\lim_{x\to-\infty} \frac{T'(x)}{\bigl(T(x)-x\bigr)^3} \\
%&= 0 - 2\lim_{x\to-\infty} \frac{\rho(x)}{\rho\bigl(T(x)\bigr)\bigl(T(x)-x\bigr)^3}
%= 0
%\end{split}
%\]

Observe that \eqref{hp:tails} and $\KE(\rho)<\infty$ are simultaneously satisfied if
\begin{equation*}
\liminf_{\abs{x}\to\infty}{} \abs{x}^{\alpha+1} \rho'(x) > 0
\qquad \text{and} \qquad
\limsup_{\abs{x}\to\infty}{} \abs{x}^{\beta+1} \rho'(x) < \infty
\end{equation*}
for $\alpha, \beta \in (1,3)$ with $1 < \beta \leq \alpha < 2\beta+1$.

%{\color{gray}
%Moreover, we need to assume the following property:
%a neighborhood of the graph of $T$ is contained in the region
%\[
%\{(x,y)\in\setR^2:\abs{T'}<M \text{ near $x$}\} \cup
%\{(x,y)\in\setR^2:\abs{T'}<M \text{ near $y$}\}
%\]
%for some constant $M>0$.
%}

\subsection{Choice of the parameters}\label{sec:parameters}

The construction of $\bar\gamma_\eps$ depends on the choice of many parameters. Here we introduce them and specify which inequalities between them are necessary for the construction. First of all, we fix a parameter $H>1$.
All the constructions in the following sections depend on $H$, which is considered as a fixed parameter, but to ease the notation we don't always make this dependence explicit.
This parameter $H$ will play a role in the proof of \autoref{thm:rec} where we then use a diagonal argument to extract a recovery sequence.

Our construction relies on several structural assumptions. First of all, we need to split the mass in two components: the main part, where most of the energy comes from, and a remaining mass, used to fix the marginals of the recovery sequence. In order to do so, we need to lower the mass $\rho$ by some constant, and for this reason we have to restrict ourselves to the bulk of the mass and avoid the tails at $\pm\infty$. Moreover, in order to be able to estimate the potential energy, we require a linearization of the map $T$ (for the application of \autoref{prop:conv-prod}), whose error depends on $T''$. Again for the application of \autoref{prop:conv-prod}, we want the slope of the linearized map to be far from $0$ and $\infty$. Since $T$ has asymptotes at $\pm\infty$ and $0$, we want then to work in a domain that stays away from those regions. Finally, we need a control on the eigenvalue of $\nabla^2V$ evaluated on $\graph T$, as this will determine the width of the Gaussians used in the construction.

All these reasons motivate the introduction of a set where we have good estimates of the objects intervening in the construction.
%
%\old{The invariant domain that we will use is the following. Given $H>1$, consider the set
%\[
%\Omega_H = [T(1/H),T(H)] \cup [1/H,H]
%\]
%and the enlarged set
%\[
%\Omega_H' = [T(1/H)-1,T(H)/2] \cup [1/(2H),H+1].
%\]}
The invariant domain that we will use is the following. Given $H>0$ such that $\rho\bigl([H,\infty)\bigr)<1/4$, define $0<r_H<H$ such that $\rho\bigl([0,r_H]\bigr)=\rho\bigl([H,\infty)\bigr)$ and consider the set
\begin{equation}\label{eqn:def-omegaH}
\Omega_H = [T(r_H),T(H)] \cup [r_H,H]
\end{equation}
and the enlarged set
\begin{equation}\label{eqn:def-omegaH'}
\Omega_H' = [T(r_{H+1}),T(H+1)] \cup [r_{H+1},H+1] \supset \Omega_H.
\end{equation}
We have that $T(\Omega_H)=T^{-1}(\Omega_H)=\Omega_H$ and $\Omega_H\nearrow\setR\setminus\{0\}$, hence $\int_{\Omega_H^c}\rho(x)\d x\to0$. Moreover $\rho\bigl(\Omega_H\cap(-\infty,0]\bigr)=\rho\bigl(\Omega_H\cap([0,\infty)\bigr)$ and similarly for $\Omega_H'$.

Let $q(x) = \lambda_1 \bigl(\nabla^2V(x,T(x))\bigr)$ the positive eigenvalue. Notice that $q(x)>0$ because, as already mentioned, from \eqref{eqn:hess-v} follows that $\nabla^2 V$ cannot be the zero matrix. Moreover, $q$ is a continuous function because the roots of a polynomial are continuous with respect to the coefficients.

As mentioned above, we want a control on some crucial quantities that influence the estimates of our construction. We therefore introduce the constant
%Let $L=L(H)>1$ be the smallest constant such that
%\begin{align*}
%L^{-1} &\leq \abs{T'(x)} \leq L, &
%\rho(x) &\geq \norm{\rho}_\infty / L, \\
%\norm{T}_{C^2} &\leq L \text{ at } x, &
%\norm{V}_{C^2} &\leq L \text{ at } (x,T(x)), \\
%L^{-1} &\leq q=\lambda_1(D^2V(x,T(x))) \leq L
%\end{align*}
%for every $x\in\Omega_H'$ and
%\[
%\Lip\oleft(\rho\rvert_{\Omega_H'}\right) \leq L, \qquad
%\Lip\oleft(u''\rvert_{\Omega_H'}\right) \leq L.
%\]
\begin{equation}\label{eq:L}
\begin{split}
L = L(H) &= \max\Bigl\{
	\norm{\rho}_{L^\infty(\Omega_H')},\ 
	\norm{1/\rho}_{L^\infty(\Omega_H')},\ 
	\norm{1/T'}_{L^\infty(\Omega_H')},\ 
	\norm{T}_{C^2(\Omega_H')},
\spliteq \phantom{\max\Bigl\{}
	\norm{u}_{C^2(\Omega_H')},\ 
	\norm{q}_{L^\infty(\Omega_H')},\ 
	\norm{1/q}_{L^\infty(\Omega_H')},
\spliteq \phantom{\max\Bigl\{}
	\Lip\oleft(\rho\rvert_{\Omega_H'}\right),\ 
	\Lip\oleft(u''\rvert_{\Omega_H'}\right)
\Bigr\}.
\end{split}
\end{equation}

%{\color{gray}
%Let $L>1$ be fixed and consider the set
%\[
%\Omega_L' = \set*{x\in\setR}{
%\begin{gathered}
%	L^{-1} \leq \abs{T'(x)} \leq L, \\
%	\rho(x) \geq \norm{\rho}_\infty / L, \\
%	\norm{T}_{C^2} \leq L \text{ at } (x,T(x)), \\
%	\norm{V}_{C^2}\leq L \text{ at } (x,T(x)),
%\end{gathered}
%}.
%\]
%Let then $\Omega_L''=\set{x\in\setR}{\dist\bigl(x,(\Omega_L')^c\bigr)>1/L}$ and $\Omega_L = \Omega_L'' \cap T^{-1}(\Omega_L'')$. Since $T\circ T=\Id$, we have $T(\Omega_L)=T^{-1}(\Omega_L)=\Omega_L$.
%Also, we immediately have that $\int_{\Omega_L^c}\rho(x)\d x\to0$ because $\Omega_L\nearrow\setR$.
%}

We now introduce the parameters on which our construction depends. They all depend on $\eps>0$, the variable that indexes the sequence of energies and the recovery sequence.
First of all, we fix a parameter $N = N(\eps) \to \infty$ faster than $\log \eps$, which will be used to truncate the Gaussian-like kernels so that they have compact support, such that
\begin{equation}\label{eqn:choice-N}
N(\eps ) = \abs{\log\eps}^{5/4} ;
\end{equation}
then a parameter $\beta = \beta(\eps) \to 0$ such that
\begin{equation}\label{eqn:choice-beta}
\eps^{1/2}N \ll \beta \ll \eps^{2/5} ,
\end{equation}
used to desingularize the matrix $\sqrt{\nabla^2V(x,T(x))}$;
then a parameter $\delta = \delta(\eps) \to 0$ such that
\begin{equation}\label{eqn:choice-delta}
(\beta N )^{1/2}\ll \delta \ll \eps^{1/8} N^{-3/5} ,
\end{equation}
which controls the resolution of the linearization of the map $T$;
and finally a parameter $\tau = \tau(\eps) \to 0$ such that
\begin{align}\label{eqn:choice-tau}
\tau &\gg \frac{\delta^2}{\eps^{1/4}} , &
\tau &\gg \frac{\eps^{1/2}N^2}{\beta} , &
\tau &\gg \frac{\beta\delta N}{\eps^{1/2}} \gg \frac{\beta\delta^2}{\eps^{1/2}} , &
\tau &\gg \beta^{1/2}N^{1/2} \gg \eps^{1/4} \gg \beta ,
\end{align}
which controls the amount of mass that we subtract from $\rho$ in order to subdivide it in main mass and remaining mass.

In the following sections, we will write many of the estimates relying only on these inequalities, so that their use will be easier in the future. Finally, we provide a choice of the parameters that fulfills the previous inequalities:
\begin{align}
\label{eqn:choice}
N(\eps) &= \abs{\log\eps}^{5/4}, &
\beta(\eps) &= \eps^{1/2} \abs{\log(\eps)}^3, &
\delta(\eps) &= \eps^{1/8} \abs{\log\eps}^{-1}, &
\tau(\eps) &= \abs{\log\eps}^{-1/3}.
\end{align}
%{\color{blue} Expected choice of parameters: $\delta = \ep^{1/8} |\log\ep|^{-1}$, $\beta(\ep) = \ep^{1/2} |\log(\ep)|$, $N(\ep) = |\log\eps|^2$. }

\subsection{Construction of \texorpdfstring{$\bar\gamma_\eps$}{\gamma\eps}}\label{sec:construction}
Let $T$ be the optimal map between $\rho$ and itself with respect to the Coulomb cost, which induces the plan $\gamma$. For every $\bigl(x,T(x)\bigr) \in \setR^2$ on the graph of $T$, let us consider the truncated Gaussian $\Gamma_{M_{\eps,\beta(x)},N}$ given by \eqref{eq:gamma-truncated} with covariance matrix $M_{\eps,\beta}(x)=\sqrt{D^2V(x,T(x))}/\eps^{1/2}+ \Id/\beta$ and truncated through $N$.

%Let $\Omega_L=\omega_1\cup\dots\cup\omega_p$ be the connected components. We subdivide each interval $\omega_j=[A_j,B_j]$ into intervals $I^1_i=[a_i,b_i]$ of length $\delta/2<b_i-a_i<\delta$.
We subdivide $\Omega_H$ into intervals $I^1_i=[a_i,b_i]$ of length $\delta/2<b_i-a_i<\delta$. Let $T_\delta$ be the piecewise linear interpolation of $T$ on the intervals $I^1_i=[a_i,b_i]$, namely a map that is affine in each $I^1_i$ and which coincides with $T$ at the boundary of each interval. We define also $T_\delta=T$ outside of $\Omega_H$. With this definition, writing $T(a_i)$ and $T(b_i)$ in terms of the Taylor expansion of $T$ at $x$, we observe that
%\begin{equation}
%\label{eqn:Tdelta-Tclose}
%|T_\delta(x)- T(x)| \leq \Big| \frac{x-x_i}{x_{i+1}-x_i} T(x_i)+ \frac{x_{i+1}-x}{x_{i+1}-x_i} T(x_{i+1}) -T(x) \Big|
%\leq \|T'' \|_{L^\infty} \delta^2 
%\end{equation}
\begin{equation}\label{eqn:Tdelta-Tclose}
\abs{T_\delta(x)- T(x)}
\leq \abs*{ \frac{x-a_i}{b_i-a_i} T(a_i) + \frac{b_i-x}{b_i-a_i} T(b_i) - T(x) }
\leq \norm{T''}_{L^\infty(\Omega_H)} \delta^2
\leq L\delta^2
\end{equation}
for every $x\in I^1_i $ and for every $ i$. Let $I^2_i=T(I^1_i)=T_\delta(I^1_i)$. Fix also points $x_i\in I_i^1$ such that $T'(x_i)=T_\delta'(x_i)$ that will be used to freeze the coefficients of the truncated Gaussians. The existence of such a point is ensured by Lagrange Theorem.
%and let $A_i=\sqrt{D^2V(x_i, T(x_i))}$.

%$$
%\Gamma[\ep, N, A,\x](\x'):= \frac{1}{C(\ep, N, A)}\left( \exp \left( -(\x'-\x)^T \frac{A}{\ep^{1/2}} (\x'-\x) - \frac{|\x'-\x|^2}{\beta(\ep)} \right) - \exp({-N})\right)_+ .
%$$
%$$
%\Gamma[M,N, \x](\x'):= \frac{1}{}\left( \exp \left( -(\x'-\x)^T \frac{A}{\ep^{1/2}} (\x'-\x) - \frac{|\x'-\x|^2}{\beta(\ep)} \right) - \exp({-N})\right)_+ .
%$$
%$$
%\Gamma[M,N, \x](\x'):= \frac{1}{}\left( \exp \left( -(\x'-\x)^T \frac{A}{\ep^{1/2}} (\x'-\x) - \frac{|\x'-\x|^2}{\beta(\ep)} \right) - \exp({-N})\right)_+ .
%$$
%Here $N:= N(\ep)$ is such that the expression in brackets has $1+[\log(\ep)]^{-1}$ mass with respect to $\exp \left( -(\x')^T \frac{A}{\ep^{1/2}} \x' - \frac{|\x'|^2}{\beta(\ep)}\right)$, and $C(\ep, N, A)$ is a renormalization which makes $\Gamma[\ep, N, A,\x]$ a probability density.

For convenience, in the following we use the notation $\Gamma_{M,N}^\x(\x')=\Gamma_{M,N}(\x'-\x)$ to denote the Gaussian with covariance matrix $M$, truncated at level $N$, centered at $\x$.

Ideally, we would like to build a recovery sequence which is a superposition of Gaussians on $\supp\gamma=\graph T$ as
\[
\gamma_\eps(\x')
= \int_{\setR^2} \Gamma_{M_{\eps,\beta}(x),N}^{\x}(\x') \d\gamma(\x)
= \int_\setR \Gamma_{M_{\eps,\beta}(x),N}^{(x,T(x))}(\x') \d\rho(x) .
%= \int_\setR \Gamma_{M_{\eps,\beta}(x),N}\bigl(\x'-(x,T(x))\bigr) \d\rho(x).
\]
Unfortunately, this plan has the wrong marginals, and there isn't a simple way to fix them without changing too much the energy. For instance one can notice that the construction in \cite{Bindini} pays an excessive price in the deconvolution step (this is analogous to applying \autoref{prop:cost-deconvolved-plan} to the whole transport plan). We follow therefore a different route.

The actual construction of the recovery sequence is more involved. We split the total mass $\rho$ in two components: a main portion of the mass and a remaining mass. The main part is dealt with in this section by superposing some truncated Gaussian kernels, whereas the remaining mass is used to bring back the marginals to what they need to be (see \autoref{sec:deconvolution}) without disrupting the energy estimates (see \autoref{sec:remaining}).% \old{The advantage of introducing  $\bar \gamma_\eps$ with respect to  $ \gamma_\eps$ relies then on the gain presented in \autoref{sec:remaining}, where we perform estimates on the remaining mass which allow to prescribe the marginals without destroying the potential and kinetic energy bounds.}

The main part of the recovery sequence is a suitable approximation of $\gamma_\eps$, which is given by the plan
%$$
%\gamma_\ep = \int \Gamma[\ep, \sqrt{D^2V(\x)}, (x, T(x))](\x') \, d\rho(x').
%$$ 
\begin{equation}\label{eq:bar-gamma-eps}
\bar\gamma_\eps(\x') = \sum_{i=1} \bar\gamma_{\eps,i}(\x'),
\end{equation}
where each piece is built by superposing along the graph of $T_\delta$ Gaussians with frozen covariance matrix $M_{\eps,\beta}(x_i)$ with $x_i\in I^1_i$, wheighing them with the density $\rho-\tau$ (which is positive in $\Omega_H'$ for $\tau$, and hence $\eps$, sufficiently small):
\begin{equation}\label{eq:bar-gamma-eps-i}
\begin{split}
\bar\gamma_{\eps,i}(\x')
&= \int_{I^1_i} \Gamma_{M_{\eps,\beta}(x_i),N}^{(x,T_\delta(x))}(\x')
	(\rho(x)-\tau) \d x \\
&= \int_{I^1_i} \Gamma_{M_{\eps,\beta}(x_i),N}\bigl(\x'-(x,T_\delta(x))\bigr)
	(\rho(x)-\tau) \d x \\
&= (\Id,T_\delta)_\#\bigl((\rho-\tau)\bm{1}_{I^1_i}\bigr)
	* \Gamma_{M_{\eps,\beta}(x_i),N} (\x').
\end{split}
\end{equation}
We denote its two marginals by
\begin{equation}\label{eq:rho-eps-i}
\rho^1_\eps = \proj^1_\# \bar\gamma_\eps, \qquad
\rho^2_\eps = \proj^2_\# \bar\gamma_\eps.
\end{equation}

%\begin{lemma}\label{lem:remaining-mass}
%Let $\gamma_\eps$ and $\bar\gamma_\eps$ be as defined above, with $\bar\gamma_\eps$ constructed with parameter $H$.
%Let $\rho^1_\eps = \proj^1_\# \bar\gamma_\eps$ and $\rho^2_\eps = \proj^2_\# \bar\gamma_\eps$ be the two marginals.
%Then, \end{lemma}
The following proposition collects the properties of $\bar \gamma_\eps$. Firstly, it shows that its marginals are quantitatively below $\rho$, a property which crucially relies on the presence of $\tau$ in the definition of  $\bar \gamma_\eps$ and on the regularity of $\rho$ and $V$. Secondly, it proves that this is still a recovery sequence (without prescribed marginals) as the one build in the proof of \autoref{thm:nonfixed-marginals} was.

\begin{proposition}\label{lem:approx-plan}
Under the same assumptions of \autoref{thm:rec}, let $\bar\gamma_\eps$ be as defined in \eqref{eq:bar-gamma-eps}, with $\bar\gamma_\eps$ constructed with parameter $H$.
Then,
%\begin{equation}\label{eqn:order-marg}
%\proj^i_\#(\bar\gamma_\eps) \leq \rho \qquad \text{for } i=1,2
%\end{equation}
%and
%\begin{equation}
%%\label{eqn:en-approx-plan}
%\lim_{\eps \to 0}
%	\bigl( E_\eps(\sqrt{\bar\gamma_\eps}) - E_\eps(\sqrt{\gamma_\eps}) \bigr) \leq 0 .
%\end{equation}
\begin{equation}\label{eqn:en-approx-plan}
\lim_{\eps\to0} E_\eps(\sqrt{\bar\gamma_\eps})
%= \lim_{\eps\to0} E_\eps(\sqrt{\gamma_\eps})
%= \lim_{\eps\to0} \int_{\Omega_H}
%	E_\eps\oleft(\sqrt{\Gamma_{M_{\eps,\beta}(x),N}^{(x,T(x))})}\right) \rho(x) \d x
%=
\leq \frac12 \int_{\Omega_H} \tr\oleft(\sqrt{\nabla^2V(x, T(x))}\right) \rho(x) \d x .
\end{equation}
Moreover, there exists $0<c_H<1$ such that for $\eps$ sufficiently small
\begin{equation}\label{eq:lower-marginal}
\rho^i_\eps(x) \leq \rho(x) - c_H \tau
\qquad \text{for } i=1,2 \text{ and } x\in\Omega_H',
\end{equation}
and
\begin{equation}\label{eq:remaining-mass}
(\rho - \rho^i_\eps)(\setR)
= 1 - \bar\gamma_\eps(\setR^2)
= \rho(\Omega_H^c) + \tau\abs{\Omega_H} .
\end{equation}

%{\color{gray}
%For every $\ep>0$, let $\delta:= \delta(\ep), N, \beta$. (in particolare dovremmo usare che $T$ e' $C^{1,1}$; moralmente questo corrisponde a $V\in C^{2,1}$ che corrisponde alla richiesta nella \eqref{eqn:ass-A-x}).
%% Assume that
%%\begin{equation}
%%\label{ratios}\lim_{\ep \to 0} \frac{\delta(\ep)^4+ \beta(\ep)N(\ep)}{ \sqrt \ep}=0. 
%%\end{equation}
%In particular, with the choice of parameters in \eqref{eqn:choice}, the right-hand side of \eqref{eqn:en-approx-plan} converges to $0$ as $\eps \to 0$.
%}
\end{proposition}

\begin{proof} We use here the notation introduced in \autoref{sec:gaussians} and \autoref{sec:rec}.

\textit{Step 1: proof of \eqref{eq:lower-marginal}}.
Let us deal first with the case of $\rho^1_\eps$. First of all, notice that $\supp \rho^1_\eps\subset \Omega_H'$ and $\rho>\tau$ on $\Omega_H'$.
%, hence
%\begin{equation}\label{eq:rho-tau-plus}
%\bigl(\rho(x)-\tau\bigr)_+ \bm{1}_{\Omega_H}(x)
%= \bigl(\rho(x)-\tau\bigr) \bm{1}_{\Omega_H}(x).
%\end{equation}
We have that
\begin{equation}\label{eq:rho1eps-sum}
\rho^1_\eps(x) = \sum_i \Bigl(\bigl[(\rho-\tau) \bm1_{I^1_i}\bigr]*\eta^1_i\Bigr)(x)
%= \sum_i \Bigl(\bigl[(\rho-\tau) \bm1_{I^1_i}\bigr]*\eta^1_i\Bigr)(x)
\end{equation}
where $\eta^1_i = \proj^1_\# \Gamma_{M_{\eps,\beta}(x_i),N}$. For every $x$ there can be at most two non zero terms contributing to the sum, because the convolutions kernels $\eta^1_i$ have support with diameter $\diam(\supp \eta^1_i) \leq C\beta^{1/2}N^{1/2}$, whereas two non-consecutive intervals are further away than $\dist(I^1_{i-1},I^1_{i+1}) \geq C \delta \gg \beta^{1/2}N^{1/2}$, thanks to \eqref{eqn:choice-delta}.

If $\dist\bigl(x,(I^1_i)^c\bigr)\geq C\beta^{1/2}N^{1/2}$, there is actually only the $i$-th term in the sum. Since the functions $(\rho-\tau) \bm1_{I^1_i}$ and $(\rho-\tau)$ coincide in $B(x,C\beta^{1/2}N^{1/2})$, when we convolve them with $\eta^1_i$ they give the same value at $x$, therefore we can compute
\[
\rho^1_\eps(x)
= \Bigl(\bigl[(\rho-\tau) \bm1_{I^1_i}\bigr] * \eta^1_i\Bigr)(x)
= \Bigl((\rho-\tau) * \eta^1_i\Bigr)(x)
\leq \rho(x) - \tau + C \Lip\oleft(\rho\rvert_{\Omega_H'}\right) \beta^{1/2} N^{1/2}.
\]

Otherwise, let $x$ be a point where the two non zero terms in the sum \eqref{eq:rho1eps-sum} are $i-1$ and $i$. Then, since the functions $\bm1_{I^1_{i-1}\cup I^1_i}$ and $1$ coincide in $B(x,C\beta^{1/2}N^{1/2})$, we have
\[
\begin{split}
\rho^1_\eps(x)
&= \Bigl(\bigl[(\rho-\tau) \bm1_{I^1_{i-1}}\bigr] * \eta^1_{i-1}\Bigr)(x)
	+ \Bigl(\bigl[(\rho-\tau) \bm1_{I^1_i}\bigr] * \eta^1_i\Bigr)(x) \\
&= \Bigl(\bigl[(\rho-\tau) \bm1_{I^1_{i-1}\cup I^1_i}\bigr] * \eta^1_{i-1}\Bigr)(x)
	+ \Bigl(\bigl[(\rho-\tau) \bm1_{I^1_i}\bigr] * (\eta^1_i-\eta^1_{i-1})\Bigr)(x) \\
&= \Bigl((\rho-\tau) * \eta^1_{i-1}\Bigr)(x)
	+ \Bigl(\bigl[(\rho-\tau) \bm1_{I^1_i}\bigr] * (\eta^1_i-\eta^1_{i-1})\Bigr)(x) \\
&\leq \rho(x) - \tau + C \Lip\oleft(\rho\rvert_{\Omega_H'}\right) \beta^{1/2} N^{1/2}
	+ \norm{\rho}_\infty \norm{\eta^1_i-\eta^1_{i-1}}_1 \\
&\leq \rho(x) - \tau + C L \beta^{1/2} N^{1/2}
	+ C_H \eps^{-1/4}\delta^2 + C_H \eps^{-1/2}\beta\delta
\end{split}
\]
because $\norm{\eta^1_i-\eta^1_{i-1}}_1 \leq \norm{\Gamma_{M_{\eps,\beta}(x_i),N} - \Gamma_{M_{\eps,\beta}(x_{i-1}),N}}_1 \leq C_H \eps^{-1/4}\delta^2 + C_H \eps^{-1/2}\beta\delta$ by \autoref{lemma:distance-gaussian-wrt-matrix} (with $\x_1=\x_2=0$).
Therefore \eqref{eq:lower-marginal} follows (with any $c_H<1$) because, from \eqref{eqn:choice-N}-\eqref{eqn:choice-tau}, we have $\tau \gg \eps^{-1/2}\beta\delta N \gg \eps^{-1/2}\beta\delta$, $\tau \gg \beta^{1/2}N^{1/2}$ and $\tau \gg \eps^{-1/4}\delta^2$.

Let us now deal with the case of $\rho^2_\eps$. From the definitions \eqref{eq:bar-gamma-eps}-\eqref{eq:rho-eps-i}, we have
\[
\begin{split}
\rho^2_\eps(y)
&= \sum_i \proj^2_\# \bar\gamma_{\eps,i}(y)
= \sum_i \proj^2_\# \left( (\Id,T_\delta)_\#\bigl((\rho-\tau)\bm{1}_{I^1_i}\bigr)
	* \Gamma_{M_{\eps,\beta}(x_i),N} \right) (y) \\
&= \sum_i \left(\proj^2_\# (\Id,T_\delta)_\# \bigl((\rho-\tau)\bm{1}_{I^1_i}\bigr)\right)
	* \left(\proj^2_\# \Gamma_{M_{\eps,\beta}(x_i),N}\right) (y) \\
&= \sum_i \left(\bigl[{T_\delta}_\# \bigl((\rho-\tau)\bm{1}_{I^1_i}\bigr)\bigr]
	* \eta^2_i \right)(y)
= \sum_i \Bigl(\bigl[{T_\delta}_\#(\rho-\tau) \bm1_{I^2_i}\bigr]*\eta^2_i\Bigr)(y)
\end{split}
\]
where $\eta^2_i = \proj^2_\# \Gamma_{M_{\eps,\beta},N}$. For $y\in\Omega_H$ we have
%\[
%{T_\delta}_\#\rho(y) = \rho\bigl(T_\delta^{-1}(y)\bigr)(T_\delta^{-1})'(y)
%= \rho\bigl(T_\delta^{-1}(y)\bigr) / T'(x_i)
%\]
\[
\begin{split}
\abs{{T_\delta}_\#\rho(y) - \rho(y)}
&= \abs{{T_\delta}_\#\rho(y) - T_\#\rho(y)}
= \abs*{\frac{\rho\bigl(T_\delta^{-1}(y)\bigr)}{T_\delta'\bigl(T_\delta^{-1}(y)\bigr)}
	- \frac{\rho\bigl(T^{-1}(y)\bigr)}{T'\bigl(T^{-1}(y)\bigr)}} \\
&\leq \frac{\abs{\rho\bigl(T_\delta^{-1}(y)\bigr)-\rho\bigl(T^{-1}(y)\bigr)}}
		{T'\bigl(T^{-1}(y)\bigr)}
	+ \rho\bigl(T_\delta^{-1}(y)\bigr)
		\frac{\abs{T_\delta'\bigl(T_\delta^{-1}(y)\bigr)-T'\bigl(T^{-1}(y)\bigr)}}
		{T_\delta'\bigl(T_\delta^{-1}(y)\bigr)T'\bigl(T^{-1}(y)\bigr)} \\
&\leq \Lip(\rho\rvert_{\Omega_H'})C\delta L + L \Lip(T'\rvert_{\Omega_H'})C\delta L^2
\leq C_H \delta,
\end{split}
\]
where we used \eqref{eq:L} to bound the norms of $\rho$ and $T'$ and the Lipschitz constants. In particular, if $y\in I^2_i$, then ${T_\delta}_\#(\rho-\tau)(y) = {T_\delta}_\#\rho(y)-\tau/T'(x_i) \leq \rho(y) + C_H\delta - \tau/L$.
%\[
%{T_\delta}_\#(\rho-\tau)(y) = \bigl(\rho(T_\delta^{-1}(y))-\tau\bigr)(T_\delta^{-1})'(y)
%= \bigl(\rho(T_\delta^{-1}(y))-\tau\bigr) / T'(x_i)
%\]
We have therefore
\[
\rho^2_\eps(y)
\leq \sum_i \Bigl(\bigl[\bigl(\rho-\tau/L+C_H\delta\bigr) \bm1_{I^2_i}\bigr]
	* \eta^2_i\Bigr)(y).
\]
We observe the similarity with \eqref{eq:rho1eps-sum}, where $\rho-\tau$, $I^1_i$ and $\eta^1_i$ have been replaced by $\rho-\tau/L+C_H\delta$, $I^2_i$ and $\eta^2_i$ respectively. We have $\diam(I^2_i) \geq C\delta/L \gg \beta^{1/2}N^{1/2}$ and $\diam(\supp\eta^2_i)\leq C\beta^{1/2}N^{1/2}$, so with the same argument as before we obtain
\[
\rho^2_\eps(y)
\leq \rho(y) - \tau/L + C_H\bigl( \delta + \beta^{1/2}N^{1/2}
	+ \eps^{-1/4}\delta^2 + \eps^{-1/2}\beta\delta \bigr),
\]
from which the thesis follows with any $c_H<1/L(H)$ because $\tau$ is asymptotically larger than all the error terms.

\textit{Step 2: proof of \eqref{eq:remaining-mass}}.
%Let's now prove \eqref{eq:remaining-mass}.
Since $\Gamma_{M_{\eps,\beta}(x_i),N}^{(x,T_\delta(x))}$ is a probability measure, we have
\[
\begin{split}
\bar\gamma_\eps(\setR^2)
&= \sum_i \int_\setR \int_{I^1_i} \Gamma_{M_{\eps,\beta}(x_i),N}^{(x,T_\delta(x))}(\x')
	(\rho(x)-\tau)\d\x' \\
&= \sum_i \int_{I^1_i} (\rho(x)-\tau) \d x
 = \int_{\Omega_H} (\rho(x)-\tau) \d x
 = \rho(\Omega_H) - \tau\abs{\Omega_H} ,
\end{split}
\]
from which
\[
1 - \bar\gamma_\eps(\setR^2)
= 1 - [\rho(\Omega_H) - \tau\abs{\Omega_H}]
= \rho(\Omega_H^c) + \tau\abs{\Omega_H} .
\]

\textit{Step 3}.
%For every $x\in\setR$, let $\x_1 = \bigl(x, T_\delta(x)\bigr)$ and $\x_2 = \bigl(x, T(x)\bigr)$.
For every $i$ and $x\in I^1_i$,
%consider the Gaussians $\Gamma_{M_{\eps,\beta}(x),N}\bigl(\x'-(x,T(x))\bigr)$ and $\Gamma_{M_{\eps,\beta}(x_i),N}\bigl(\x'-(x,T_\delta(x))\bigr)$.
we claim the estimate
\begin{equation}\label{eqn:step1-ts}
\begin{split}
E_\eps&\oleft(\sqrt{\Gamma^{(x,T_\delta(x))}_{M_{\eps,\beta}(x_i),N}}\right)
	- E_\eps\oleft(\sqrt{\Gamma_{M_{\eps,\beta}(x),N}^{(x,T(x))}}\right) \\
%&\leq C_H \delta^2\eps^{-3/4}
%	\bigl( \delta^4 + \eps^{1/2}N + \beta N \delta^2 \bigr)
%	+ C_H \eps^{1/2} e^{-N/2} + C_H\delta . \\
& \leq C_H (\eps^{-3/4}\delta^2 + \eps^{-1}\beta\delta)
	\bigl( \delta^4 + \eps^{1/2}N + \beta N \delta^2 \bigr)
	+ C_H \eps^{1/2} e^{-N/2} + C_H\delta .
%\frac{C_H \delta} {\sqrt \eps}\bigl( \delta^3+ \beta N + \delta \beta^{1/2} N^{1/2} \bigr).
\end{split}
\end{equation}

Notice that the right hand side of \eqref{eqn:step1-ts} goes to $0$. In fact, thanks to \eqref{eqn:choice-delta} and the fact that $N\to\infty$ we have
\[
\beta N\delta^2 \ll \delta^4 \ll \eps^{1/2} \ll \eps^{1/2}N,
\]
therefore, thanks to \eqref{eqn:choice-N}-\eqref{eqn:choice-delta}, the right hand side of \eqref{eqn:step1-ts} is less than
%\[
%C_H\frac{\delta^2}{\eps^{1/4}} N + C_He^{-N/2} + C_H\delta \ll \frac{C_H}{N}.
%\]
\[
\begin{split}
C_H(\eps^{-3/4}\delta^2 &+ \eps^{-1}\beta\delta) \eps^{1/2}N
	+ C_H\eps^{1/2}e^{-N/2} + C_H\delta \\
&\leq C_H(\eps^{-1/4}\delta^2 + \eps^{-1/2}\beta\delta)N
	+ C_H\eps^{1/2}e^{-N/2} + C_H\delta \\
&\ll N \abs{\log\eps}^{-2} + \eps^{1/40}N^{2/5} + e^{-N/2} + \delta \ll 1 .
\end{split}
\]

%{\color{gray}Notice that, with the particular choice of the parameters given by \eqref{eqn:choice}, asymptotically the right hand side of \eqref{eqn:step1-ts} is less than $C_H\abs{\log\eps}^{-1/2}$, which converges to $0$ as $\eps\to0$.}

%Notice that, applying only the upper bounds on $\beta$ and $\delta$ appearing in \eqref{eqn:choice-beta} and \eqref{eqn:choice-delta}, the right-hand side is bounded by $\eps^{1/4}N^2$, which converges to $0$ as $\eps \to 0$.
%We prove the following claim.
%Let $v = (v_1, v_2)$ be such that $|v| \leq \delta^2$ and let 
% $$\eta_1^\ep = \left(e^{ -  ( x^t A_1 x + \delta |y|^2 ) / 2\ep } - e^{-N} \right)_+, \qquad \eta_1^\ep = \left(e^{ -  ( x^t \sqrt{D^2V(x,T(x))} x + \delta |y|^2 ) / 2\ep } - e^{-N} \right)_+.$$
We observe that by \eqref{eqn:Tdelta-Tclose} the points $\x_1=\bigl(x,T_\delta(x)\bigr),\x_2=\bigl(x,T(x)\bigr) \in \R^2$ and the symmetric matrices $A_1=\sqrt{\nabla^2V\bigl(x_i,T(x_i)\bigr)}$, $A_2=\sqrt{\nabla^2V(\x_2)}$ which generate $M_1=M_{\eps,\beta}(x_i)$ and $M_2=M_{\eps,\beta}(x)$ satisfy
\begin{equation}
\label{eqn:ass-A-x}
\abs{\x_1-\x_2} \leq L\delta^2 \qquad \text{and} \qquad
\abs{A_1-A_2} \leq C_H\delta.
\end{equation}
In fact, $\bigl(x_i,T(x_i)\bigr)$ and $\bigl(x,T(x)\bigr)$ belong to a region of $\graph(T)$ where $\nabla^2 V$ is Lipschitz.
As a consequence,
\begin{equation}\label{eq:norm-difference-M}
\abs{M_1-M_2}
= \abs*{\left(\frac{A_1}{\eps^{1/2}}-\frac{I}{\beta}\right)
	-\left(\frac{A_2}{\eps^{1/2}}-\frac{I}{\beta}\right)}
= \eps^{-1/2} \abs{A_1-A_2} \leq C_H \eps^{-1/2}\delta .
\end{equation}
%Then
%\begin{equation}
%\label{eqn:step1-ts}
%%\lim_{\ep \to 0} 
%E_\ep((\Gamma[\ep, N, A_1,\x_1])^{1/2}) - E_\ep((\Gamma[\ep, N, A_2,(x,T(x))])^{1/2}) \leq ...
%\end{equation}
%uniformly in $x, v$ and $\bar x$.

%In order to prove this claim, we observe that the optimal transport map between $\Gamma[\ep, N, A_1,\x_1]$ and $\Gamma[\ep, N, A_2,(x,T(x))]$ is an affinity. TO DO

As regards the potential energy, denoting by $E$ the union of the supports of the probability measures appearing in \eqref{eqn:step1-ts}
\[
E = \supp \Gamma^{(x,T_\delta(x))}_{M_{\eps,\beta}(x_i),N}
	\cup \supp \Gamma_{M_{\eps,\beta}(x),N}^{(x,T(x))}
\]
we have
\begin{equation}\label{eqn:approx-pot}
\begin{split}
\int_{\setR^2} & V(\x') \Gamma^{(x,T_\delta(x))}_{M_{\eps,\beta}(x_i),N} (\x') \d\x'
	- \int_{\setR^2} V(\x') \Gamma_{M_{\eps,\beta}(x),N}^{(x,T(x))}(\x') \d\x' \\
&\leq \left(\sup_E V\right) \int_{\setR^2}
	\abs*{\Gamma^{(x,T_\delta(x))}_{M_{\eps,\beta}(x_i),N} (\x') \d\x'
	- \Gamma_{M_{\eps,\beta}(x),N}^{(x,T(x))}(\x')} \d\x'
\end{split}
\end{equation}
%we notice that the Lipschitz constant of $V$ in a $\delta_0$-neigborhood of $(x, T(x))$ (which belongs to the $0$-level set of $V$) is estimated by $C\delta_0$.
%We notice that both the probability measures appearing in \eqref{eqn:step1-ts} are supported in a set of diameter less than $C \beta(\ep)^{1/2}N(\eps)^{1/2}$ and that their centers (and hence their supports) differ by at most $\delta^2$ in view of \eqref{eqn:ass-A-x}. Hence we can estimate 
%\begin{equation}
%\label{eqn:claim1}
%\sup_E V \leq C (\delta^2 + \beta(\ep)^{1/2}N(\eps)^{1/2})^2,
%\end{equation}
%
%OSSERVAZIONE: questa stima si puo' migliorare un po' perche' questo ellisse e' molto schiacciato nella direzione in cui V fa 0. Questo migliora la stima nel caso in cui $\beta(\ep)$ e $\sqrt \ep $ sono abbastanza diversi. Dovremmo ottenere:
Observe that, by \autoref{lemma:distance-gaussian-wrt-matrix}, we have
\begin{equation}
\label{eqn:claim2}
\int \Big|\Gamma^{(x,T_\delta(x))}_{M_{\eps,\beta}(x_i),N}(\x') d\x' - \Gamma_{M_{\eps,\beta}(x),N}^{(x,T(x))}(\x')\Big| d\x' 
\leq C_H \eps^{-1/4}\delta^2 + C_H \eps^{-1/2}\beta\delta .
\end{equation}

We claim that $E$ is contained in a strip of size $C (\delta^2 + \eps^{1/4}N^{1/2} + \beta^{1/2} N^{1/2} \delta \Lip(T'))$ around the graph of the affine map which equals $T_\delta$ in $I^1_i$
\begin{equation}
\label{eqn:Econtained}
E \subseteq S= \{ \x \in \R^2 : \dist\bigl(\x, \graph(T_\delta)\bigr) \leq C (\delta^2 + \eps^{1/4}N^{1/2} + \beta^{1/2} N^{1/2} \delta \Lip(T'))\}.
\end{equation}
Indeed, by construction $\supp \Gamma^{(x,T_\delta(x))}_{M_{\eps,\beta}(x_i),N} $ is a rectangle centered on the graph of $T_\delta$ and with the longer side aligned with  the graph of $T_\delta$. Moreover, its width is $\eps^{1/4}N^{1/2}$; hence the claim \eqref{eqn:Econtained} is proved for the first rectangle in the definition of $E$. 

On the other side, we observe that the rectangle $\supp \Gamma_{M_{\eps,\beta}(x),N}^{(x,T(x))}$ is centered in the point $(x,T(x))$, whose distance from the graph of $T_\delta$ is estimated by $ \|T'' \|_{L^\infty} \delta^2 $ by \eqref{eqn:Tdelta-Tclose}; moreover, this rectangle is tilted with respect to the previous one proportionally to $|T'(x)-T'_\delta(x)| \leq \delta \Lip(T')$, which in turn gives a contribution of $ \beta^{1/2} N^{1/2} \delta \Lip(T')$. This proves \eqref{eqn:Econtained}.

Since for every $x$ the point $(x,T(x))$ and the point are at distance at most $(x,T_\delta(x))$ $\|T'' \|_{L^\infty} \delta^2$ by \eqref{eqn:Tdelta-Tclose}, from \eqref{eqn:Econtained} we deduce that $S$ is contained in the strip of width $C (1+\|T'' \|_{L^\infty} ) (\delta^2 + \eps^{1/4}N^{1/2} + \beta^{1/2} N^{1/2} \delta )$ around the graph of $T$. Hence, recalling that both $V$ and $\nabla V$ vanish along the graph of $T$, and more precisely recalling that $V$ has uniform quadratic growth around the graph of $T$ has stated in  we obtain that
\begin{equation}
\label{eqn:claim1}
\sup_E V \leq C (1+\|T'' \|_{L^\infty} )^2 (\delta^2 + \eps^{1/4}N^{1/2} + \beta^{1/2} N^{1/2} \delta )^2.
\end{equation}
Hence we deduce by \eqref{eqn:approx-pot}, \eqref{eqn:claim2}, and \eqref{eqn:claim1} that
\begin{equation}\label{eqn:view}
\begin{split}
\int V(\x')& \Gamma^{(x,T_\delta(x))}_{M_{\eps,\beta}(x_i),N} (\x') \d\x'
	- \int V(\x') \Gamma_{M_{\eps,\beta}(x),N}^{(x,T(x))}(\x') \d\x' \\
&\leq C_H (\eps^{-1/4}\delta^2 + \eps^{-1/2}\beta\delta)
	(\delta^2 + \eps^{1/4}N^{1/2} + \beta^{1/2} N^{1/2} \delta)^2 \\
&\leq C_H (\eps^{-1/4}\delta^2 + \eps^{-1/2}\beta\delta)
	(\delta^4 + \eps^{1/2}N + \beta N \delta^2).
\end{split}
\end{equation}
In view of \eqref{eqn:view}, we deduce the corresponding statement to \eqref{eqn:step1-ts} for the potential energy.

As regards the kinetic energy, we observe that it equals
\[
\begin{split}
\KE(\Gamma_{M_1,N}^{\x_1}) - \KE(\Gamma_{M_2,N}^{\x_2})
&= \KE(\Gamma_{M_1,\infty}^{\x_1}) - \KE(\Gamma_{M_2,\infty}^{\x_2})
		\spliteq
	+ [\KE(\Gamma_{M_1,N}^{\x_1}) - \KE(\Gamma_{M_1,\infty}^{\x_1})]
		\spliteq
	- [\KE(\Gamma_{M_2,N}^{\x_2}) - \KE(\Gamma_{M_2,\infty}^{\x_2})] .
\end{split}
\]
The two terms in square brackets are estimated by \eqref{eq:truncated-kinetic-energy} and they are negligible with respect to the other contribution, whereas the first difference is
\[
\KE(\Gamma_{M_1,\infty}^{\x_1}) - \KE(\Gamma_{M_2,\infty}^{\x_2})
= \frac{\tr M_1}{4} - \frac{\tr M_2}{4}
= \frac{q_1-q_2}{4\eps^{1/2}}
\leq \frac{\delta}{4\eps^{1/2}}
\]
because, by the variational characterization of the largest eigenvalue, we have
\[
q_1 = \max_{\abs{v}\leq1}v^TA_1v
\leq \max_{\abs{v}\leq1}v^TA_2v + \max_{\abs{v}\leq1}v^T(A_1-A_2)v
\leq q_2 + C_H \delta .
\]
Therefore
\[
\KE(\Gamma_{M_1,N}^{\x_1}) - \KE(\Gamma_{M_2,N}^{\x_2})
\leq C(\tr M_1 + \tr M_2) e^{-N/2}
	+ \frac{C_H\delta}{\eps^{1/2}}\leq C_H e^{-N/2}
	+ \frac{C_H\delta}{\eps^{1/2}} .
\]

\textit{Step 4}. We conclude the proof. By convexity of the kinetic energy, since the potential energy is linear and thanks to \eqref{eqn:step1-ts}, with the remark following it, we have that
\[
\begin{split}
\limsup_{\eps\to0} E_\eps(\sqrt{\bar\gamma_\eps})
&\leq \limsup_{\eps\to0} \sum_{i=1} \int_{I^1_i}
	E_\eps\oleft(\sqrt{\Gamma_{M_{\eps,\beta}(x_i),N}^{(x,T_\delta(x))}}\right)
	\rho(x) \d x \\
&\leq \limsup_{\eps\to0} \int_{\Omega_H}
	E_\eps\oleft(\sqrt{\Gamma_{M_{\eps,\beta}(x),N}^{(x,T(x))})}\right) \rho(x) \d x .
\end{split}
\]
We now estimate the integrand in the right-hand side with $\tr\oleft(\sqrt{\nabla^2V(x,T(x))}\right) \rho(x) $ up to small errors. To this end, we bound the kinetic energy by passing to the non-truncated Gaussians with the estimate \eqref{eq:truncated-kinetic-energy} and by the exact computation \eqref{eqn:daqualcheparte}.
We obtain
\[
\begin{split}
%\biggl\lvert E_\eps\oleft(\sqrt{\Gamma_{M_{\eps,\beta}(x),N}^{(x,T(x))})}\right)
%	&- \tr M_{\eps,\beta}(x)
%	\biggr\rvert \\
\eps^{1/2} KE&\oleft(\sqrt{\Gamma_{M_{\eps,\beta}(x),N}^{(x,T(x))})}\right)
=  \eps^{1/2}\KE\oleft(\sqrt{\Gamma_{M_{\eps,\beta}(x),N}}\right) \\
&\leq \eps^{1/2}\KE\oleft(\sqrt{\Gamma_{M_{\eps,\beta}(x),\infty}}\right)
	+ C \eps^{1/2} \tr\bigl(M_{\eps,\beta}(x)\bigr) e^{-N/2} \\
&= \eps^{1/2}\frac{\tr M_{\eps,\beta}(x)}{4} + Ce^{-N/2}
\leq \frac 14\tr\oleft(\sqrt{\nabla^2V(x,T(x))}\right)
	+ \frac{\eps^{1/2}}{2\beta} + Ce^{-N/2}
% \\
%&\leq \eps^{1/2} \abs*{\KE\oleft(\sqrt{\Gamma_{M_{\eps,\beta}(x),N}}\right)
%	- \KE\oleft(\sqrt{\Gamma_{M_{\eps,\beta}(x),\infty}}\right)}
%	\spliteq
%	+ \eps^{-1/2} \int_{\setR^2} \abs*{\sqrt{\nabla^2V(x)}\y}^2
%		\abs{\Gamma_{M_{\eps,\beta}(x),N}(\y)
%		-\Gamma_{M_{\eps,\beta}(x),\infty}(\y)} \d\y
%	\spliteq
%	+ \eps^{-1/2} \int_{\setR^2} \abs*{V(\x+\y) - \abs*{\sqrt{\nabla^2V(x)}\y}^2}
%		\Gamma_{M_{\eps,\beta}(x),N}(\y) \d\y \\
%&\leq C \eps^{1/2} \tr\bigl(M_{\eps,\beta}(x)\bigr) e^{-N/2}
%	+ \norm*{\sqrt{\nabla^2V(x)}}_\infty^2 N e^{-N/2}
%	\spliteq
%	+ C_H \eps^{-1/2} \int_{\setR^2} \abs{\y}^3 \Gamma_{M_{\eps,\beta}(x),N}(\y) \d\y \\
%&\leq C e^{-N/2} + C N e^{-N/2}
%	+ C_H \eps^{-1/2} \beta^{3/2} N^{3/2} \\
%%	+ C_H \eps^{-1/2} \frac{\sqrt{ab}}{\pi} \int_{\setR^2}
%%		(w^2+z^2)^{3/2} e^{-aw^2-bz^2} \d w \d z \\
%%&\leq C N e^{-N/2} + C_H \eps^{-1/2} \frac{\sqrt{ab}}{\pi} \int_{\setR^2}
%%		2\bigl(\abs{w}^3+\abs{z}^3\bigr) e^{-aw^2-bz^2} \d w \d z \\
%%&\leq C N e^{-N/2} + C_H \frac{2}{\sqrt\pi} \eps^{-1/2} \left(a^{-3/2}+b^{-3/2}\right) \\
%&\leq C_H \eps^{-1/2} \beta^{-3/2} N^{3/2}
.
\end{split}
\]
%where in the last inequality we just point out that all errors appearing from the previous line go to zero better than $ \abs{\log\eps}^{-1/2} $.
In order to estimate the potential energy, we need to compare the $V$ with its second order Taylor expansion.
Since $\nabla^2V$ is locally Lipschitz in a neighborhood of $\graph T$, there is $r>0$ such that for every $\x\in\supp\gamma\cap\Omega_H^2$ and $\y\in B_r(0)$ we have
\[
\abs*{V(\x+\y) - \frac12\abs*{\sqrt{\nabla^2V(\x)}\y}^2}
\leq \norm{\nabla^3V}_{L^\infty} \cdot \abs{\y}^3
\leq C_H \abs{\y}^3,
\]
where the $L^\infty$ norm is taken in a neighborhood of size $(\beta N)^{1/2}$ of $\supp\gamma\cap\Omega_H^2$.
Then, thanks to \eqref{eq:truncated-quadratic-energy}\footnote{When applying \eqref{eq:truncated-quadratic-energy} to the integral of $\abs*{\sqrt{\nabla^2V(\x)}\y}^2$ times the difference of the Gaussians, we use the fact that the matrix $\sqrt{\nabla^2V(\x)}$ is bounded from above by $\eps^{1/2}M_{\eps,\beta}(x)$.} and \eqref{eqn:daqualcheparte}, and since $|\y| \leq (\beta N)^{1/2}$ for $\y \in \supp \Gamma_{M_{\eps,\beta}(x),N}$, we have
\[
\begin{split}
%\biggl\lvert E_\eps\oleft(\sqrt{\Gamma_{M_{\eps,\beta}(x),N}^{(x,T(x))})}\right)
%	&- \tr M_{\eps,\beta}(x)
%	\biggr\rvert \\
&\eps^{-1/2}\int V(\y) \Gamma_{M_{\eps,\beta}(x),N}^{(x,T(x))}(\y)\, d\y \\
&\leq \frac{\eps^{-1/2}}2 \int_{\setR^2} \abs*{\sqrt{\nabla^2V(x,T(x))}\y}^2
	\Gamma_{M_{\eps,\beta}(x),N}(\y) \d\y
%	\spliteq
	+ C_H \eps^{-1/2} \int_{\setR^2} |\y|^3 \Gamma_{M_{\eps,\beta}(x),N}(\y) \d\y \\
&\leq \frac{\eps^{-1/2}}2 \int_{\setR^2} \abs*{\sqrt{\nabla^2V(x,T(x))}\y}^2
	\Gamma_{M_{\eps,\beta}(x),\infty}(\y) \d\y
%	\spliteq
	+ C \eps^{-1/2} N e^{-N/2} + C_H \eps^{-1/2} \beta^{3/2} N^{3/2} \\
&\leq \frac14 \tr\oleft(\sqrt{\nabla^2V(x,T(x))}\right)
	+ C \eps^{-1/2} N e^{-N/2} + C_H \eps^{-1/2} \beta^{3/2} N^{3/2} .
\end{split}
\]
%Therefore
%\[
%\begin{split}
%E_\eps\oleft(\sqrt{\Gamma_{M_{\eps,\beta}(x),N}^{(x,T(x))})}\right)
%&\leq \eps^{1/2}\KE\oleft(\sqrt{\Gamma_{M_{\eps,\beta}(x),\infty}}\right)
%	+ \eps^{-1/2} \int_{\setR^2} \abs*{\sqrt{\nabla^2V(x)}\y}^2
%		\Gamma_{M_{\eps,\beta}(x),\infty}(\y) \d\y
%	\spliteq
%	+ C_H \eps^{-1/2}\beta^{3/2}N^{3/2} \\
%&\leq \eps^{1/2}\KE\oleft(\sqrt{\Gamma_{M_{\eps,\beta}(x),\infty}}\right)
%	+ \eps^{1/2} \int_{\setR^2} \abs*{M_{\eps,\beta}(x)\y}^2
%		\Gamma_{M_{\eps,\beta}(x),\infty}(\y) \d\y
%	\spliteq
%	+ C_H \eps^{-1/2}\beta^{3/2}N^{3/2} \\
%&= \eps^{1/2}\frac{\tr M_{\eps,\beta}(x)}{2}
%	+ \eps^{1/2}\frac{\tr M_{\eps,\beta}(x)}{2}
%	+ C_H \eps^{-1/2}\beta^{3/2}N^{3/2} \\
%&= \tr\oleft(\sqrt{\nabla^2V(x)}\right) + 2\eps^{1/2}\beta^{-1}
%	+ C_H \eps^{-1/2}\beta^{3/2}N^{3/2}  \\
%&= \tr\oleft(\sqrt{\nabla^2V(x)}\right) + O_L\oleft(\abs{\log\eps}^{-3}\right) .
%\end{split}
%\]
%
%%We now put together all the estimates.
%\[
%\begin{split}
%E_\eps(\sqrt{\bar\gamma_\eps})
%&\leq \sum_{i=1} \int_{I^1_i}
%	E_\eps\oleft(\sqrt{\Gamma_{M_{\eps,\beta}(x_i),N}^{(x,T_\delta(x))}}\right)
%	\rho(x) \d x \\
%&\leq \int_{\Omega_H}
%	E_\eps\oleft(\sqrt{\Gamma_{M_{\eps,\beta}(x),N}^{(x,T(x))})}\right) \rho(x) \d x
%	+ C_H \abs{\log\eps}^{-1/2}
%	 \\
%&\leq \int_{\Omega_H} \tr\oleft(\sqrt{\nabla^2V(x)}\right) \rho(x) \d x
%	+ O_L\oleft(\abs{\log\eps}^{-3}\right)
%	+ C_H \abs{\log\eps}^{-1/2} ,
%\end{split}
%\]
Overall, from the estimates in Step 4 we obtain that 
\[
\begin{split}
%\limsup_{\eps\to0}
E_\eps(\sqrt{\bar\gamma_\eps})
&\leq \frac12 \int_{\Omega_H} \tr\oleft(\sqrt{\nabla^2V(x,T(x))}\right) \rho(x) \d x
	\splitleq
	+ C_H \left( \frac{\eps^{1/2}}\beta + Ce^{-N/2}
		+ C \eps^{-1/2} N e^{-N/2} + C_H \eps^{-1/2} \beta^{3/2} N^{3/2} \right) ,
\end{split}
\]
which concludes the proof of \eqref{eqn:en-approx-plan} since by \eqref{eqn:choice-N} and \eqref{eqn:choice-beta} we have
\[
\begin{split}
\frac{\eps^{1/2}}\beta + Ce^{-N/2}
	&+ C \eps^{-1/2}N\eps^{-1/2} N e^{-N/2} + C_H \eps^{-1/2} \beta^{3/2} N^{3/2} \\
&\ll \frac1N + C e^{-N/2} + \eps^{1/10}N^{3/2} \ll 1 .
\end{split}
\]
As a side note, it could be easily seen by means of the kind of computations performed in this proof that \eqref{eqn:en-approx-plan} is actually an equality, but since we will not need this in the sequel we don't pursue this matter here.
\end{proof}

\begin{lemma}\label{lemma:distance-gaussian-wrt-matrix}
Let $\x_1,\x_2\in\setR^2$ with $\abs{\x_1-\x_2}\leq C_0\delta^2$ and let $A_1,A_2$ be two degenerate positive-semidefinite symmetric $2\times2$ matrices with $\abs{A_i}\leq L$ and $\abs{A_1-A_2}\leq C_0\delta$; finally let $M_i=M_{\eps,\beta}(A_i)=\frac{A_i}{\eps^{1/2}}+\frac{I}{\beta}$. Then there exists a positive constant $C_H$ only depending on $H$ and $C_0$ such that
\[
\norm{\Gamma_{M_1,N}^{\x_1}-\Gamma_{M_2,N}^{\x_2}}_1
\leq C_H \eps^{-1/4}\delta^2 + C_H \eps^{-1/2}\beta\delta.
\]
%and
%\begin{equation}\label{eqn:dist-gauss2}
%\norm*{
%	\frac{\tilde\Gamma_{M_1,\infty}^{\x_1}}{G_{M_1,N}}
%		\bm{1}_{\supp \Gamma_{M_1,N}}
%	-\frac{\tilde\Gamma_{M_2,\infty}^{\x_2}}{G_{M_2,N}}
%		\bm{1}_{\supp \Gamma_{M_2,N}}
%	}_1 \leq C \frac{\beta\delta N}{\eps^{1/2}} .
%\end{equation}
\end{lemma}

Observe that, under the assumptions \eqref{eqn:choice-N}-\eqref{eqn:choice-delta}, the right hand side goes to $0$ when $\eps\to0$.

\begin{proof}
Thanks to \eqref{eq:L1-norm-trunc} we have
\[
\norm{\Gamma_{M_1,N}^{\x_1}-\Gamma_{M_2,N}^{\x_2}}_1
\leq \norm{\Gamma_{M_1,\infty}^{\x_1}-\Gamma_{M_2,\infty}^{\x_2}}_1 + CNe^{-N/2}.
\]
Notice that, by \eqref{eqn:choice-N}-\eqref{eqn:choice-delta}, $Ne^{-N/2} \ll N^2\eps^{1/4}=\eps^{1/2}N^2\eps^{-1/4}\ll\beta N\eps^{-1/4}\ll\delta^2\eps^{-1/4}$, so we can forget about this error term because it is dominated by the right hand side of the thesis.

We split the estimate into first changing only the center point, and then changing the matrix. By the triangle inequality we have
\begin{equation}\label{eq:L1-norm-point-then-matrix}
\norm{\Gamma_{M_1,\infty}^{\x_1}-\Gamma_{M_2,\infty}^{\x_2}}_1
\leq \norm{\Gamma_{M_1,\infty}^{\x_1}-\Gamma_{M_1,\infty}^{\x_2}}_1
	+ \norm{\Gamma_{M_1,\infty}^{\x_2}-\Gamma_{M_2,\infty}^{\x_2}}_1 .
\end{equation}

The first norm can be estimated in the following way.
If $a$ and $b$ are the eigenvalues of $M_1$, we have $\nabla\Gamma_{M_1,\infty}(w,z) = 2 (aw,bz)^T \Gamma_{M_1,\infty}(w,z)$,
therefore
\[
\begin{split}
\norm{\nabla\Gamma_{M_1,\infty}}_1
%= \int_{\setR^2} \abs{\nabla\Gamma_{M_1,\infty}(\y)} \d\y
&\leq G_{M_1,\infty}^{-1} \int_{\setR^2} (2a\abs{w}+2b\abs{z})e^{-aw^2-bz^2} \d w \d z \\
&= \frac{\sqrt{ab}}\pi \left(\frac{2\sqrt\pi}{\sqrt b}+\frac{2\sqrt\pi}{\sqrt a}\right)
= \frac{2}{\sqrt\pi}(\sqrt a+\sqrt b)
\leq C q(x)^{1/2} \eps^{-1/4}
\leq C_H \eps^{-1/4} .
\end{split}
\]
From this we get
\[
\begin{split}
\norm{\Gamma_{M_1,\infty}^{\x_1}-\Gamma_{M_1,\infty}^{\x_2}}_1
&= \int_{\setR^2} \abs{\Gamma_{M_1,\infty}(\y+\x_2-\x_1)-\Gamma_{M_1,\infty}(\y)} \d\y \\
&= \int_{\setR^2} \abs*{\int_0^1 \nabla\Gamma_{M_1,\infty}\bigl(\y+s(\x_2-\x_1)\bigr)
	\cdot (\x_2-\x_1) \d s} \d\y \\
&\leq \int_{\setR^2} \int_0^1 \abs*{\nabla\Gamma_{M_1,\infty}\bigl(\y+s(\x_2-\x_1)\bigr)}
	\cdot \abs{\x_2-\x_1} \d s \d\y \\
&\leq C\delta^2 \norm{\nabla\Gamma_{M_1,\infty}}_1
\leq C_H \delta^2 \eps^{-1/4}.
\end{split}
\]

Let's now turn to the second norm in \eqref{eq:L1-norm-point-then-matrix}. For $t\in[0,1]$, define the matrix $B_t=M_1+t(M_2-M_1)$. Notice that
\[
B_t = \frac1{\eps^{1/2}}\oleft(A_1+\frac{\eps^{1/2}}{\beta}I+t(A_2-A_1)\right)
\]
and $\abs{A_2-A_1}<\delta\ll\eps^{1/2}/\beta$, because $\beta\delta\ll \eps^{2/5+1/8}N^{-3/5}\ll\eps^{21/40}\ll\eps^{1/2}$ by \eqref{eqn:choice-beta}-\eqref{eqn:choice-delta}. Therefore we have that for $\eps$ small $1/2\det M_1 \leq \det B_t \leq 2\det M_1$ and $1/2\tr M_1 \leq \tr B_t \leq 2\tr M_1$.
Moreover,
\[
\begin{split}
\abs*{\frac{\d}{\d t}\log \det B_t}
&= \frac1{\det B_t}\abs*{\tr\oleft(\adj(B_t) \frac{\d B_t}{\d t}\right)}
= \frac{\abs*{\tr\bigl(\adj(B_t)(M_2-M_1)\bigr)}}{\det B_t} \\
&\leq C_H \frac{\abs{B_t} \cdot \abs{M_2-M_1}}{\det B_t}
\leq C_H \frac{\eps^{-1/2}\cdot\delta\eps^{-1/2}}{\eps^{-1/2}\beta^{-1}}
= C_H \frac{\beta\delta}{\eps^{1/2}}
\end{split}
\]
and
\[
\begin{split}
\frac{\d}{\d t}\Gamma_{B_t,\infty}(\x)
&= \frac{\d}{\d t} \oleft(\frac{\sqrt{\det B_t}}{\pi} e^{-\x^TB_t\x}\right) \\
&= \frac{\frac{\d}{\d t}\det B_t}{2\pi\sqrt{\det B_t}} e^{-\x^TB_t\x}
	- \x^T\frac{\d B_t}{\d t}\x \frac{\sqrt{\det B_t}}{\pi} e^{-\x^TB_t\x} \\
%&= \left(\frac{\d}{\d t}\det B_t\right) \frac{\Gamma_{B_t,\infty}}{2\det B_t}
%	- \x^T(M_2-M_1)\x \Gamma_{B_t,\infty} \\
&= \frac12 \left(\frac{\d}{\d t}\log\det B_t\right) \Gamma_{B_t,\infty}
	- \x^T(M_2-M_1)\x \Gamma_{B_t,\infty} .
\end{split}
\]
Therefore
\[
\begin{split}
\norm{\Gamma_{M_1,\infty}^{\x_2}-\Gamma_{M_2,\infty}^{\x_2}}_1
&= \norm{\Gamma_{M_1,\infty}-\Gamma_{M_2,\infty}}_1
%&= \int_{\setR^2} \abs*{\int_0^1 \frac{\d}{\d t}\Gamma_{B_t,\infty}(\x) \d t} \d\x
\leq \int_0^1 \int_{\setR^2} \abs*{\frac{\d}{\d t}\Gamma_{B_t,\infty}(\x)} \d\x \d t \\
&\leq C_H \int_0^1 \abs*{\frac{\d}{\d t}\log\det B_t}
	\int_{\setR^2} \Gamma_{B_t,\infty} \d\x \d t
		\spliteq
	+ \int_0^1 \abs{M_2-M_1} \frac{\sqrt{\det B_t}}{\pi}
		\int_{\setR^2} \abs{\x}^2 e^{-\x^TB_t\x}\d\x\d t \\
&= C_H \frac{\beta\delta}{\eps^{1/2}}
	+ \frac{\delta}{\eps^{1/2}} \int_0^1 \frac{\tr B_t}{2\det B_t} \d t \\
&\leq C_H \frac{\beta\delta}{\eps^{1/2}} + C_H \frac{\delta}{\eps^{1/2}}\beta
= C_H \frac{\beta\delta}{\eps^{1/2}} . \qedhere
\end{split}
\]

\end{proof}

\section{Deconvolution of plans}
\label{sec:deconvolution}

Let $\sigma^1,\sigma^2\in\Meas_+(\setR)$ with the same mass and $\Pi_0\in\Pi(\sigma^1,\sigma^2)$ be a transport plan. Fix a radial convolution kernel $\Theta\in C^\infty(\setR^2;[0,1])$ with $\supp\Theta\subset B(0,1)$. Define the rescaling $\Theta_\eps(x)=\eps^{-1/2}\Theta(\eps^{-1/4}x)$ and the marginals $\theta(x)=(\proj^1_\#\Theta)(x)$, $\theta_\eps(x)=(\proj^1_\#\Theta_\eps)(x)=\eps^{-1/4}\theta(\eps^{-1/4}x)$.
Define the convolved plan
\[
\Pi_\eps=\Pi_0*\Theta_\eps.
\] 
Let $\sigma^i_\eps = \proj^i_\#\Pi_\eps = \sigma^i*\theta_\eps$, for $i=1,2$, be the two marginals of the convolved plan.

Define now the deconvolved plan as introduced in \cite[Theorem 6.3]{Bindini}
\[
\begin{split}
\tilde\Pi_\eps(x, y)
&= \int_{\setR^2}
	\frac{\sigma^1(x)\theta_\eps(x'-x)}{\sigma^1_\eps(x')}
	\frac{\sigma^2(y)\theta_\eps(y'-y)}{\sigma^2_\eps(y')}
	\Pi_\eps(x',y') \d x' \d y' \\
&= \int_{\setR^2} P(x, y, x', y') \d x' \d y' .
\end{split}
\]
We  verify that $\tilde\Pi_\eps \d x \d y\in\Pi(\sigma^1,\sigma^2)$. Indeed, for any integrable function $F: \R^2 \to \R$ and for every $x\in \R$ we perform this useful computation, integrating first the variable $y$ and then $y'$
\begin{equation}\label{eqn:useful0int}
\begin{split}
 \int_{\setR^3} \! F(x,x') P(x, y, x', y') 
 \d x' \d y' \d y
% &= \int_{\setR^3} \! F(x,x') \frac{\sigma^1(x)\theta_\eps(x'-x)}{\sigma^1_\eps(x')}
%	\frac{\sigma^2(y)\theta_\eps(y'-y)}{\sigma^2_\eps(y')}
%	\Pi_\eps(x',y') \d x' \d y' \d y \\
&= \int_{\setR^2} \! F(x,x') \frac{\sigma^1(x)\theta_\eps(x'-x)}{\sigma^1_\eps(x')}
	\cancel{\frac{\sigma^2_\eps(y')}{\sigma^2_\eps(y')}}
	\Pi_\eps(x',y') \d x' \! \d y' \\
&= \int_\setR F(x,x') \frac{\sigma^1(x)\theta_\eps(x'-x)}{\cancel{\sigma^1_\eps(x')}}
	\cancel{\sigma^1_\eps(x')} \d x'.
%= \sigma^1(x)
\end{split}
\end{equation}
By taking the function $F\equiv 1$  this reads%we find in particular 
\[
\begin{split}
\int_{\setR} \tilde\Pi_\eps(x,y) \d y
&= {\sigma^1(x)}.
%= \sigma^1(x)
\end{split}
\]
%{\color{gray}\[
%\begin{split}
%\int_\setR \tilde\pi_\eps(x,y) \d y
%&= \int_{\setR^3} \frac{\sigma^1(x)\theta_\eps(x'-x)}{\sigma^1_\eps(x')}
%	\frac{\sigma^2(y)\theta_\eps(y'-y)}{\sigma^2_\eps(y')}
%	\pi_\eps(x',y') \d x' \d y' \d y \\
%&= \int_{\setR^2} \frac{\sigma^1(x)\theta_\eps(x'-x)}{\sigma^1_\eps(x')}
%	\cancel{\frac{\sigma^2_\eps(y')}{\sigma^2_\eps(y')}}
%	\pi_\eps(x',y') \d x' \d y' \\
%&= \int_\setR \frac{\sigma^1(x)\theta_\eps(x'-x)}{\cancel{\sigma^1_\eps(x')}}
%	\cancel{\sigma^1_\eps(x')} \d x'
%= \sigma^1(x)
%\end{split}
%\]
%}
Analogously one can show that the second marginal of $\tilde\Pi_\eps$ is $\sigma^2$. Moreover, $\tilde\Pi_\eps$ is supported in a neighborhood of size $ \eps^{1/4}$ of the support of $\Pi_\eps$ because, in the definition of $\tilde\Pi_\eps$, $(x',y')$ is in the support of $\Pi_\eps$ and $x' -x$ is less than $\eps^{1/4}$. Therefore 
\begin{equation}
\label{eqn:supp-gran-pi}
\mbox{$\tilde\Pi_\eps$ is supported in a neighborhood of size $ 2\eps^{1/4}$ of the support of $\Pi_0$.}
\end{equation}

We claim in the next proposition that the operation of deconvolving a transport plan $\Pi_0$ worsens its potential energy in a controlled way, namely proportionally to the total mass of $\Pi_0$, while on the other hand the kinetic energy of the deconvolution  is controlled only in terms of the kinetic energy of the marginals of $\Pi_0$. %(which may be infinite for $\Pi_0$, for instance if it consists of a singular measure) is . More precisely, as $\eps \to 0$ the kinetic energy is the same, whereas the potential energy is worsened proportionally to the total mass of $\Pi_0$ (assuming that it lies in a precise neighborhood of the zero level set of $V$).
The loss in potential energy proportional to the mass of $\Pi_0$ is the reason why in the proof of \autoref{thm:rec} we use this deconvolution procedure only on the remaining mass, rather than on the full plan $\bar \gamma_\eps$. Let us also point out that, in the proof of \autoref{thm:rec}, $\Pi_0$ will in turn be an $\eps$-dependent plan, built in Section~\ref{sec:remaining} below. For this reason we keep explicit all dependences on $\Pi_0$ in the next proposition.

\begin{proposition}\label{prop:cost-deconvolved-plan}
Let $V$ be as in \autoref{thm:rec}, $H>0$ and $\Omega_H, \Omega_H'$ be the sets introduced in \eqref{eqn:def-omegaH}, \eqref{eqn:def-omegaH'}. Let  $\sigma^1,\sigma^2\in\Meas_+(\setR)$ with the same mass, let $\Pi_0\in\Pi(\sigma^1,\sigma^2)$, and let $\tilde \Pi_\eps$ be the plan introduced at the beginning of this section. Then there exists a universal constant $C>0$
\begin{equation}
\label{eqn:kin-est-deconv}
\KE(\tilde\Pi_\eps)
%= \frac12 \int_{\setR^2} \abs*{\nabla\sqrt{\tilde\Pi_\eps}}^2
\leq \KE(\sigma^1) + \KE(\sigma^2) + \frac{C}{\sqrt\eps}\norm{\Pi_0}_1
\qquad \text{for every $\eps>0$}.
\end{equation}

Moreover, there exists a constant $c_H>0$ depending only on $H$ such that if
\begin{equation}
\label{eqn:deconv-hp}
\dist\bigl((x,y),\graph(T)\bigr)
%\leq C_H\eps^{1/4}N^{1/2}\bm{1}_{\Omega_H'}(x)
\leq c_H\bm{1}_{\Omega_H'}(x)
\qquad \text{for every } (x,y)\in\supp\Pi_0,
\end{equation}
%then for $\eps$ sufficiently small and for a universal constant $C>0$ (independent on $L$), we have
then there exists a universal constant $C>0$ and a constant $C_H>0$ depending only on $H$ such that for $\eps$ sufficiently small
\begin{equation}
\label{eqn:ts-deconv-pote}\frac{1}{ \sqrt{\eps}} \int_{\setR^2} V\d\tilde\Pi_\eps
%\leq C_H\int_{\setR^2} V\d\Pi_0 + C\sqrt{\eps}\norm{\Pi_0}_1 %+ CL^2\delta^4\norm{\Pi_0}_1
%	+ CL^2\sqrt\eps\norm{\Pi_0}_1
\leq \frac{C_H}{ \sqrt{\eps}} \int_{\setR^2} V\d\Pi_0 + C\norm{\Pi_0}_1.
\end{equation}
%If moreover we assume that %{\color{red} CAMBIAMO?}
%\[
%\supp\Pi_0 \subset \set*{(x,y)\in\setR^2}
%	{\dist\bigl((x,y),\graph(T)\bigr) < 2L\eps^{1/4}%(\delta^2+\eps^{1/4})
%	}
%	,
%\]
%then for $\eps$ sufficiently small, we have
%\[
%\int_{\setR^2} V\d\tilde\pi_\eps
%\leq %\int_{\setR^2} V\d\pi + C_H\sqrt{\eps}\norm{\pi}_1 + %CL^2\delta^4\norm{\pi}_1+ CL^2\sqrt\eps\norm{\pi}_1
%%=
% \int_{\setR^2} V\d\pi + C_H\sqrt{\eps}\norm{\pi}_1+ o_L(\eps^{1/2}) .
%\]
\end{proposition}
We notice that the dependence of the constants in \eqref{eqn:ts-deconv-pote} is of fundamental importance: in particular
%, when we renormalize the potential energy by the factor $\sqrt \eps$ 
the first term will go to $0$ (at $H$ fixed) by the particular choice of $\Pi_0$, while the second term will be estimated with \eqref{eq:remaining-mass} by the mass of $\rho$ outside $\Omega_H$ (times the constant $C$ independent of $H$); hence this will be small sending $H \to \infty$ uniformly in $\eps$.

\begin{remark} %Several deconvolution procedures have been proposed in the literature. 
Another estimate of the potential energy of the deconvolved plan in terms of the original one was proposed by Lewin \cite{Lewin}. 
%In our notations, %introducing his deconvolved plan (see ... )   $\bar{\Pi}_{\ep}$ and 
%tracing in his proof the dependence on  $\| \Pi_0\|_1$, we would get an estimate of the form
%$$  \left| \int_{\R^2} V\, d \bar{\Pi}_{\ep} - \int_{\R^2} V \, d \Pi_0 \right|
%= \left| \int_{\R^2} V_{ee}\, d \bar{\Pi}_{\ep} - \int_{\R^2} V_{ee} \, d \Pi_0 \right|
% \leq C\sqrt{\ep} \|\rho\|_{\infty}^3 \cdot \left( \int |\nabla \sigma_1| + { \| \Pi_0\|_1} \right),$$
%where the first equality simply follows from the fact that   $\bar{\Pi}_{\ep}$ and $\Pi_0$ have the same marginals.

%Another approach that could work in our case is the one proposed by Lewin ... . 
In our notations, %introducing his deconvolved plan (see ... )   $\bar{\Pi}_{\ep}$ and 
tracing in his proof the dependence on  $\| \Pi_0\|_1$, we would get an estimate of the form
$$
\left| \int_{\R^2} \Phi \, d \bar{\Pi}_{\ep} - \int_{\R^2} \Phi \, d \Pi_0 \right| \leq \sqrt{\ep} \left( \| \nabla \Phi \|_{\infty} \int |\nabla \sigma_1| + \| D^2 \Phi \|_{\infty} { \| \Pi_0\|_1} \right) 
$$
where the norms of $\nabla \Phi$ and $D^2\Phi$ are calculated in the neighbourhood of the graph of $T$ where $\bar{\Pi}_{\ep}$ is supported. Since   $\bar{\Pi}_{\ep}$ and $\Pi_0$ have the same marginals, we notice that
$$   \int_{\R^2} V \, d \bar{\Pi}_{\ep} - \int_{\R^2} V \, d \Pi_0   =     \int_{\R^2} V_{ee} \, d \bar{\Pi}_{\ep} - \int_{\R^2} V_{ee} \, d \Pi_0,$$

 Applying the estimate above for $\Phi=V_{ee}$ and using that  $|x-y| \geq 1/ (2\| \rho\|_{\infty})$ on the graph of $T$, we would get that for $\ep$ small enough% we obtain
$$  \left| \int_{\R^2} V\, d \bar{\Pi}_{\ep} - \int_{\R^2} V \, d \Pi_0 \right| \leq C\sqrt{\ep} \|\rho\|_{\infty}^3 \sqrt{ \| \Pi_0\|_1} \cdot \left( \int |\nabla \sigma_1|+ \sqrt{ \| \Pi_0\|_1} \right).$$
However, this estimate is not good enough for our application since the $BV$ norm of $\sigma_1$ does not go to $0$ as $\eps \to 0$ (recall that $\sigma_1$ coincides with $\rho$ in a neighborhood of $0$).
%We choose instead to use the deconvolution of Bindini and make our result self-contained.
%Another approach that could work in our case is the one proposed by Lewin ... . In our notations, %introducing his deconvolved plan (see ... )   $\bar{\Pi}_{\ep}$ and 
%tracing in his proof the dependence on  $\| \Pi_0\|_1$, we would get an estimate of the form
%$$
%\left| \int_{\R^2} \Phi \, d \bar{\Pi}_{\ep} - \int_{\R^2} \Phi \, d \Pi_0 \right| \leq \sqrt{\ep} \sqrt{ \| \Pi_0\|_1} \cdot \left( \| \nabla \Phi \|_{\infty} \int |\nabla \sigma_1|} + \| D^2 \Phi \|_{\infty} \sqrt{ \| \Pi_0\|_1} \right) 
%$$
%where the norms of $\nabla \Phi$ and $D^2\Phi$ are calculated in the neighbourhood of the graph of $T$ where $\bar{\Pi}_{\ep}$ is supported. Next, since   $\bar{\Pi}_{\ep}$ and $\Pi_0$ have the same marginals, we notice that
%$$   \int_{\R^2} V \, d \bar{\Pi}_{\ep} - \int_{\R^2} V \, d \Pi_0   =     \int_{\R^2} V_{ee} \, d \bar{\Pi}_{\ep} - \int_{\R^2} V_{ee} \, d \Pi_0,$$
%
% Applying the estimate above for $\Phi=V_{ee}$ and using that  $|x-y| \geq 1/ (2\| \rho\|_{\infty})$ on the graph of $T$, we would get that for $\ep$ small enough% we obtain
%$$  \left| \int_{\R^2} V\, d \bar{\Pi}_{\ep} - \int_{\R^2} V \, d \Pi_0 \right| \leq C\sqrt{\ep} \|\rho\|_{\infty}^3 \sqrt{ \| \Pi_0\|_1} \cdot \left( \sqrt{ KE(\sigma_1)} + \sqrt{ \| \Pi_0\|_1} \right).$$
%We choose instead to use the deconvolution of Bindini and make our result self-contained.
\end{remark}

\begin{proof}
We have that
\[
\begin{split}
\KE(\tilde\Pi_\eps)
&= \frac12 \int_{\setR^2} \abs*{\nabla\sqrt{\tilde\Pi_\eps(\x)}}^2 \d\x
= \int_{\setR^2} \frac{\abs{\partial_x\tilde\Pi_\eps}^2
	+ \abs{\partial_y\tilde\Pi_\eps}^2}{8\tilde\Pi_\eps} \d x.
\end{split}
\]
Let us deal with the first term; the second one is treated analogously. % and we just have to remember to add it at the end.
 By H\"older inequality and since $	\int_{\setR^2} P \d x' \d y'= \tilde \Pi_\eps(x,y)$, we have
\[
\begin{split}
\abs{\partial_x\tilde\Pi_\eps(x,y)}^2
&= \abs*{\int_{\setR^2}
	\frac{\partial}{\partial x}\left(
		\frac{\sigma^1(x)\theta_\eps(x'-x)}{\sigma^1_\eps(x')}\right)
	\frac{\sigma^2(y)\theta_\eps(y'-y)}{\sigma^2_\eps(y')}
	\Pi_\eps(x',y') \d x' \d y' }^2 \\
&= \abs*{\int_{\setR^2} \left(
	\frac{\partial_x\sigma^1(x)}{\sigma^1(x)}
	+ \frac{\partial_x\theta_\eps(x'-x)}{\theta_\eps(x'-x)}\right)
	P(x,y,x',y') \d x' \d y' }^2 \\
%&\leq \left[ \int_{\setR^2} \left(
%		\frac{\partial_x\sigma^1(x)}{\sigma^1(x)}
%		+ \frac{\partial_x\theta_\eps(x'-x)}{\theta_\eps(x'-x)} \right)^2
%		P \d x' \d y'\right]
%	\int_{\setR^2} P \d x' \d y' \\
&\leq \tilde\Pi_\eps(x,y) \int_{\setR^2} \left[
		\left(\frac{\partial_x\sigma^1(x)}{\sigma^1(x)}\right)^2
		+ \left(\frac{\partial_x\theta_\eps(x'-x)}{\theta_\eps(x'-x)}\right)^2
	\right] P(x,y,x',y') \d x' \d y' .
\end{split}
\]
Then
\[
\begin{split}
\int_{\setR^2}
	& \frac{\abs{\partial_x\tilde\Pi_\eps(x,y)}^2}{8\tilde\Pi_\eps(x,y)} \d x \d y \\
&= \int_{\setR^4} \frac18 \oleft[
		\left(\frac{\partial_x\sigma^1(x)}{\sigma^1(x)}\right)^2
		+ \left(\frac{\partial_x\theta_\eps(x'-x)}{\theta_\eps(x'-x)}\right)^2
	\right] P(x,y,x',y') \d x \d y\d x' \d y' .
\end{split}
\]
Applying \eqref{eqn:useful0int} with $F(x,x')$ given by the expression in square brackets above and we find
\[
\begin{split}
\int_{\setR^2} \frac{\abs{\partial_x\tilde\Pi_\eps(x,y)}^2}{8\tilde\Pi_\eps(x,y)} \d x \d y
&= \int_{\setR^4} \frac18 \oleft[
		\left(\frac{\partial_x\sigma^1(x)}{\sigma^1(x)}\right)^2
		+ \left(\frac{\partial_x\theta_\eps(x'-x)}{\theta_\eps(x'-x)}\right)^2
	\right]  {\sigma^1(x)\theta_\eps(x'-x)}\d x' \d x \\
&= \int_\setR \frac{\abs{\partial_x\sigma^1(x)}^2}{8\sigma^1(x)} \d x
	+ \norm{\sigma^1}_1 \int_\setR
		\frac{\abs{\partial_x\theta_\eps(z)}^2}{8\theta_\eps(z)} \d z \\
&= \KE(\sigma_1) + \KE(\theta_\eps) \norm{\Pi_0}_1
%	\int_{\setR^3}  \left[ 
%	\frac{\abs{\partial_x\sigma^1(x)}^2}{\sigma^1(x)}
%	\frac{\theta_\eps(x'-x)}{\sigma^1_\eps(x')}
%	\cancel{\frac{\sigma^2_\eps(y')}{\sigma^2_\eps(y')}}
%	\pi_\eps(x',y')
%	+
%	\frac{\abs{\partial_x\theta_\eps(x'-x)}^2}{\theta_\eps(x'-x)}
%	\frac{\sigma^1(x)}{\sigma^1_\eps(x')}
%	\cancel{\frac{\sigma^2_\eps(y')}{\sigma^2_\eps(y')}}
%	\pi_\eps(x',y')
%	\right] \d x \d x' \d y' \\
 .
\end{split}
\]

The kinetic energy of the marginal of the rescaled convolution kernel scales like
\begin{equation*}\label{eqn:kinet-kernel}
\begin{split}
\KE(\theta_\eps)
= \frac12 \int_\setR \abs*{\partial_x\sqrt{\theta_\eps(x)}}^2 \d x
= \int_\setR \frac{\abs{\theta_\eps'(x)}^2}{8\theta_\eps(x)} \d x
%&= \int_\setR \frac{\abs*{\eps^{-1/2}\theta(\eps^{-1/4}x)}^2}
%	{\eps^{-1/4}\theta(\eps^{-1/4}x)} \d x
= \int_\setR \eps^{-3/4} \frac{\abs{\theta'(y)}^2}{8\theta(y)} \eps^{1/4}\d y
= \eps^{-1/2} \KE(\theta) 
\end{split}
\end{equation*}
and the kinetic energy of $\theta$ is a universal constant since $\theta$ is fixed.
In conclusion, putting everything together, we obtain \eqref{eqn:kin-est-deconv}.
%\[
%\KE(\tilde\pi_\eps)
%\leq \KE(\sigma^1) + \KE(\sigma^2) + \frac{C}{\sqrt\eps}\norm{\pi}_1 . \qedhere
%\]

Let $c_{0, H}$ be any constant strictly less than the distance of $\graph T$ from the diagonal $\{ x=y\}$. For $\eps$ sufficiently small, $\supp \Pi_\eps$ will be supported at positive distance from the diagonal thanks to \eqref{eqn:supp-gran-pi}, as soon as the constant $c_0$ appearing there is less or equal than $c_{0,H}$.
By using the Taylor expansion of $V$ centered at $(x',y')$ and the fact that thanks to \autoref{lem:good-rho} $\nabla^2V$ is bounded by a universal constant $C$ independent of $H$ in the region of interest because we are far from the diagonal, we have that
\[
\begin{split}
%\int_{\setR^2} V(x,y)\d\pi_\eps(x,y)
\int_{\setR^2} V\d\Pi_\eps
%&= \int_{\setR^4} V(x,y)\Theta_\eps(x-x',y-y')\d\pi(x',y')\d x\d y \\
%&\leq \int_{\setR^4} \Bigl[V(x',y') + \nabla V(x',y')(x-x',y-y')
%	\\ &\mathrel{\phantom{\leq}}{\phantom{\int_{\setR^4} \Bigl[}}
%	+C\bigl((x-x')^2+(y-y')^2\bigr)\Bigr] \Theta_\eps(x-x',y-y')\d\pi(x',y')\d x\d y \\
&\leq \int_{\setR^4} V(x',y') \Theta_\eps(x-x',y-y')\d\Pi_0(x',y')\d x\d y
	\\ &\mathrel{\phantom{\leq}}{}
	+ \int_{\setR^4} \nabla V(x',y')\cdot(x-x',y-y')
		\Theta_\eps(x-x',y-y')\d\Pi_0(x',y')\d x\d y
	\\ &\mathrel{\phantom{\leq}}{}
	+ \int_{\setR^4} C\bigl(\abs{x-x'}^2+\abs{y-y'}^2\bigr)
		\Theta_\eps(x-x',y-y')\d\Pi_0(x',y')\d x\d y \\
&\leq \int_{\setR^2} V\d\Pi_0
	+ C\int_{\setR^4} \bigl(\abs{x-x'}^2+\abs{y-y'}^2\bigr)
		\Theta_\eps(x-x',y-y')\d\Pi_0(x',y')\d x\d y
\end{split}
\]
because the  integral in the second line vanishes (it's the integral of a linear function in $x,y$ with respect to a symmetric kernel).
The terms involving $\abs{x-x'}^2$ and $\abs{y-y'}^2$ in the last integral can be both computed as
\[
\begin{split}
%\mathrel{\phantom{=}}
\int_{\setR^4} \abs{y-y'}^2 \Theta_\eps(x-x',y-y') \d\Pi_0(x',y') \d x \d y& = \int_{\setR^3} \abs{y-y'}^2 \theta_\eps(y-y') \d\Pi_0(x',y') \d y \\
&= \int_{\setR^2} \abs{y-y'}^2 \theta_\eps(y-y') \d y \d\sigma^2(y') \\
&= C\sqrt\eps \norm{\sigma^2}_1
= C\sqrt\eps \norm{\Pi_0}_1 .
\end{split}
\]
Therefore
\[
\int V \d\Pi_\eps
\leq \int V \d\Pi_0 + C\sqrt\eps\norm{\Pi_0}_1 .
\]

Using again the Taylor expansion of $V$ centered at $(x',y')$ we can proceed to estimate
\begin{equation}
\label{eqn:maledetto}
\begin{split}
\int V(x,y)\d\tilde\pi_\eps(x,y)
&\leq \int_{\setR^4} V(x',y') P(x,y,x',y')\d x\d y\d x'\d y'
	\\ &\mathrel{\phantom{\leq}}{}
	+ \int_{\setR^4} \nabla V(x',y')\cdot(x-x',y-y')
		P(x,y,x',y')\d x\d y\d x'\d y'
	\\ &\mathrel{\phantom{\leq}}{}
	+ \int_{\setR^4} C\bigl(\abs{x-x'}^2+\abs{y-y'}^2\bigr)
		P(x,y,x',y')\d x\d y\d x'\d y' .
%&\leq \int_{\setR^4} V(x',y') P(x,y,x',y')\d x\d y\d x'\d y'
%	\\ &\mathrel{\phantom{\leq}}{}
%	+ \int_{\setR^4} \nabla V(x',y')\cdot(x-x',y-y')
%		P(x,y,x',y')\d x\d y\d x'\d y'
%	\\ &\mathrel{\phantom{\leq}}{}
%	+ \int_{\setR^4} C\bigl(\abs{x-x'}^2+\abs{y-y'}^2\bigr)
%		P(x,y,x',y')\d x\d y\d x'\d y' \\
\end{split}
\end{equation}
By the definition of $P$ and integrating first the variable $x$ and secondly the variable $y$, the first integral in the right-hand side is
\[
\begin{split}
\int_{\setR^4} V(x',y') P%(x,y,x',y')
\d x\d y\d x'\d y' 
%&= \int_{\setR^4} V(x',y')
%	\frac{\sigma^1(x)\theta_\eps(x'-x)}{\sigma^1_\eps(x')}
%	\frac{\sigma^2(y)\theta_\eps(y'-y)}{\sigma^2_\eps(y')}
%	\Pi_\eps(x',y') \d x\d y\d x'\d y' \\
&= \int_{\setR^4} V(x',y')
	\cancel{\frac{\sigma^1_\eps(x')}{\sigma^1_\eps(x')}}
	\frac{\sigma^2(y)\theta_\eps(y'-y)}{\sigma^2_\eps(y')}
	\Pi_\eps(x',y') \d y\d x'\d y' \\
&= \int_{\setR^4} V(x',y')
	\cancel{\frac{\sigma^2_\eps(y')}{\sigma^2_\eps(y')}}
	\Pi_\eps(x',y') \d x'\d y' = \int_{\setR^2} V\d\Pi_\eps .
\end{split}
\]
%{\color{red} sicuri che non potevamo evitare di passare da $\Pi_\eps$?}
%We observe that, since $\nabla V$ vanishes on $\graph(T)$ and $D^2 V$ is bounded by Lemma ... with a constant independent of $L$, 
%\[
%\abs{\nabla V(x',y')} \leq C\dist\bigl((x',y'),\graph(T)\bigr) \leq C L(\delta^2+\eps^{1/4}).
%\]

Given %$(x,y)\in(\Omega_H')^2$ near $\graph(T)$
$x\in \Omega_H'$ and $y$ such that $d((x,y), \graph T) \leq c_{0,H}$, let $(x_0,y_0)\in\graph(T)$ be a projection, i.e.\ a point realising the minimum of the distance, and let $v=(x-x_0,y-y_0)$ be the difference vector and $d=\dist\bigl((x,y),\graph(T)\bigr)=\abs{v}$ the distance. By looking at the Taylor expansion of $V$ at $(x_0, y_0)$, recalling that $\nabla V(x_0, y_0)=0$, we have that
\[
\abs{V(x,y)-v^T \nabla^2V(x_0,y_0) v} \leq \norm{\nabla^3 V}_\infty d^3 \leq C_H d^3,
\]
where the $L^\infty$ norm is taken over the set of points $(x,y)$ such that $x\in \Omega_H'$ and $y$ such that $d((x,y), \graph T) \leq c_{0,H}$. This $L^\infty$ norm is in turn estimated by $C_H$ thanks to the explicit computation of the derivatives of $V$ (see \autoref{lem:good-rho}) and thanks to the quantities that we have under control at $H$ fixed in \eqref{eq:L}.
The difference vector $v$ is orthogonal to $\graph(T)$ at $(x_0,y_0)$, therefore it is in the direction of the eigenvector of $\sqrt{\nabla^2 V(x_0,y_0)}$ associated to the non-zero eigenvalue $q(x_0,y_0)$. Therefore
\[
V(x,y)
\geq v^T \nabla^2 V(x_0,y_0) v - \norm{\nabla^3 V}_\infty d^3
\geq q^2 d^2 - C_H d^3
\geq \frac12 \Bigl(\min_{\Omega_H'} q\Bigr)^2 d^2
\]
as soon as $d\leq \frac12 C_H^{-1} \bigl(\min_{\Omega_H'} q\bigr)^2$. Therefore, for
\begin{equation}
\label{eqn:cond-xy}
\dist\bigl((x,y),\graph(T)\bigr) < \min\big\{ c_{0,H}, \frac12 C_H^{-1} \bigl(\min_{\Omega_H'} q\bigr)^2 \big\}=:c_L,
\end{equation}
 we have
\begin{equation}\label{eq:nabla-V}
\abs{\nabla V(x,y)}^2 \leq \norm{\nabla^2 V}_{\infty}^2 d^2
\leq 2 \norm{\nabla^2 V}_{\infty}^2 \norm{1/q}_{L^\infty(\Omega_H')}^2 V(x,y)
\leq C_H V(x,y) .
\end{equation}
Using Cauchy-Schwarz, the second integral in \eqref{eqn:maledetto} can be bounded by
\[
\begin{split}
\int_{\setR^4} & \nabla V(x',y')\cdot(x-x',y-y') P(x,y,x',y')\d x\d y\d x'\d y' \\
&\leq \int_{\setR^4} \bigl(\abs{\nabla V(x',y')}^2 + (\abs{x-x'}^2+\abs{y-y'}^2)\bigr)
	P(x,y,x',y')\d x\d y\d x'\d y' .
%&\leq \bigl(CL(\delta^2+\eps^{1/4})\bigr)^2 \norm{\pi}_1
%	+ \int_{\setR^4} \bigl(\abs{x-x'}^2+\abs{y-y'}^2\bigr)
%		P(x,y,x',y')\d x\d y\d x'\d y', \\
\end{split}
\]
%{\color{gray}Using \eqref{eqn:deconv-hp}, which implies $\abs{\nabla V(x',y')}^2\leq C_H\eps^{1/2}N\bm{1}_{\Omega_H'(x)}$,}
Using \eqref{eq:nabla-V} in the point $(x' y')$, which for $\eps$ sufficiently small satisfies \eqref{eqn:cond-xy} thanks to \eqref{eqn:deconv-hp} and \eqref{eqn:supp-gran-pi}, the first term can be estimated as
\[
\begin{split}
\int_{\setR^4} & \abs{\nabla V(x',y')}^2 P(x,y,x',y') \d x \d y \d x' \d y' 
= \int_{\setR^2} \abs{\nabla V(x',y')}^2 \Pi_\eps(x',y') \d x' \d y' \\
&= \int_{\setR^4} \abs{\nabla V(x,y)}^2 \Theta_\eps(x-x',y-y')
	\d\Pi_0(x',y') \d x \d y \\
&\leq 2\int_{\setR^4} \abs{\nabla V(x',y')}^2 \Theta_\eps(x-x',y-y')
	\d\Pi_0(x',y') \d x \d y
		\spliteq
	+ 2\int_{\setR^4} \abs{\nabla V(x,y)-\nabla V(x',y')}^2 \Theta_\eps(x-x',y-y')
	\d\Pi_0(x',y') \d x \d y \\
&\leq 2\int_{\setR^2} \abs{\nabla V(x',y')}^2 \d\Pi_0(x',y')
		\spliteq
	+ 2\bigl(\Lip(\nabla V)\bigr)^2 \int_{\setR^4} \bigl(\abs{x-x'}^2+\abs{y-y'}^2\bigr)
		\Theta_\eps(x-x',y-y')\d\Pi_0(x',y')\d x\d y \\
%&\leq 2\int_{\setR^2} C_H^2\eps^{1/2}N\bm{1}_{\Omega_H'}(x') \d\pi(x',y')
%	+ C\eps^{1/2} \norm{\pi}_1 \\
%&\leq 2C_H \int_{\setR^2} V(x',y')\d\Pi_0(x',y') + C\eps^{1/2} \norm{\Pi_0}_1 \\
&\leq C_H \int_{\setR^2} V\d\Pi_0 + C\sqrt \eps\norm{\Pi_0}_1 .
\end{split}
\]
The third integral in \eqref{eqn:maledetto} can be computed exploiting \eqref{eqn:useful0int} with $F(x,x') = |x-x'|^2$ to get
\[
\begin{split}
\mathrel{\phantom{=}} \int_{\setR^4} \abs{x-x'}^2P%(x,y,x',y') 
\d x \d y \d x' \d y'
%&= \int_{\setR^3} \abs{y-y'}^2
%	\cancel{\frac{\sigma^1_\eps(x')}{\sigma^1_\eps(x')}}
%	\frac{\sigma^2(y)\theta_\eps(y'-y)}{\sigma^2_\eps(y')}
%	\pi_\eps(x',y') \d y \d x' \d y' \\
%&= \int_{\setR^2} \abs{y-y'}^2
%	\frac{\sigma^2(y)\theta_\eps(y'-y)}{\cancel{\sigma^2_\eps(y')}}
%	\cancel{\sigma^2_\eps(y')} \d y \d y' \\
&= \int_{\setR^2} \abs{x-x'}^2 \theta_\eps(x'-x)\sigma^1(x)\d x\d x'
\\& 
=\norm{\sigma^1}_1  \int_{\setR} |z|^2 \theta_\eps(z) \d z
= C\sqrt\eps \norm{\sigma^1}_1 = C\sqrt\eps \norm{\Pi_0}_1
\end{split}
\]
and analogously for the term with $|y-y'|^2$.
%{\color{gray}
%as
%\[
%\begin{split}
%&\mathrel{\phantom{=}} \int_{\setR^4} \abs{y-y'}^2
%	\frac{\sigma^1(x)\theta_\eps(x'-x)}{\sigma^1_\eps(x')}
%	\frac{\sigma^2(y)\theta_\eps(y'-y)}{\sigma^2_\eps(y')}
%	\pi_\eps(x',y') \d x \d y \d x' \d y' \\
%&= \int_{\setR^3} \abs{y-y'}^2
%	\cancel{\frac{\sigma^1_\eps(x')}{\sigma^1_\eps(x')}}
%	\frac{\sigma^2(y)\theta_\eps(y'-y)}{\sigma^2_\eps(y')}
%	\pi_\eps(x',y') \d y \d x' \d y' \\
%&= \int_{\setR^2} \abs{y-y'}^2
%	\frac{\sigma^2(y)\theta_\eps(y'-y)}{\cancel{\sigma^2_\eps(y')}}
%	\cancel{\sigma^2_\eps(y')} \d y \d y' \\
%&= \int_\setR \left(\int_\setR \abs{y'}^2\d \theta_\eps(y')\right) \d\sigma^2(y)
%= \frac{\sqrt\eps}{2} \norm{\sigma^2}_1 = \frac{\sqrt\eps}{2} \norm{\pi}_1.
%\end{split}
%\]
%}
Putting the pieces together, we obtain
\[
\begin{split}
\int_{\setR^2} V\d\tilde\Pi_\eps
&\leq C_H \int_{\setR^2} V\d\Pi_\eps + C\sqrt{\eps}\norm{\Pi_0}_1 . \qedhere
\end{split}
\]

\end{proof}

%\begin{proposition}\label{prop:kin-cost-deconvolved-plan}
%With the notation introduced at the beginning of this section, we have
%
%\end{proposition}
%
%\begin{proof}
%\end{proof}

\section{Potential cost of the remaining mass and kinetic energy of its marginals}\label{sec:remaining}

In \autoref{sec:rec} we described how to build the main portion of the recovery sequence $\bar\gamma_\eps$. As stated in \eqref{eq:lower-marginal}-\eqref{eq:remaining-mass}, this procedure leaves out a small positive mass that we still have to account for. The goal of this section is to build a transport plan for the remaining mass which satisfies suitable energy bounds and allows to apply the deconvolution procedure of \autoref{sec:deconvolution}, namely \autoref{prop:cost-deconvolved-plan}. To this end we need to build a transport plan $\pi\in\Pi(\rho-\rho^1_\eps, \rho-\rho^2_\eps)$ with potential energy smaller than $\sqrt \eps$ and to have a suitable control on how the kinetic energy of its marginals explodes in $\eps$. This is the content of the following proposition and this section is dedicated to its proof.
\begin{proposition}\label{prop:cost-rem-mass}
Let $\rho\in C^1(\setR)$, $\gamma\in\Pi(\rho,\rho)$, $u\in W^{2,\infty}$ and $V(x,y)=\abs{x-y}^{-1}-u(x)-u(y)$ be as in \autoref{thm:rec}.
Let $\bar\gamma_\eps$ be the plan introduced in \eqref{eq:bar-gamma-eps} and \eqref{eq:bar-gamma-eps-i} and let $\rho^1_\eps = \proj^1_\# \bar\gamma_\eps$ and $\rho^2_\eps = \proj^2_\# \bar\gamma_\eps$ be its two marginals introduced in \eqref{eq:rho-eps-i}.
Then for every $c>0$ and for some constant $C_H>1$ depending only on $H$
\begin{enumerate}[i)]
\item (construction of a plan with $o_L(\sqrt \eps)$ potential energy) for every $\eps $ sufficiently small
there is $\pi_\eps\in\Pi(\rho-\rho^1_\eps, \rho-\rho^2_\eps)$ such that 
\begin{equation}\label{ts:PE-bound}
\int V\d\pi_\eps
\leq C_H\left( \frac{\delta^4}{\tau} + \frac{\eps N^2}{\beta\tau}
	+ \frac{\beta\delta^2}{\tau} \right)
\end{equation}
and
\begin{equation}
\label{eqn:deconv-hp-verif}
\dist\bigl((x,y),\graph(T)\bigr)
%\leq C_H\eps^{1/4}N^{1/2}\bm{1}_{\Omega_H'}(x)
\leq c \bm{1}_{\Omega_H'}(x)
\qquad \text{for every } (x,y)\in\supp\Pi_0;
\end{equation}
\item (bound on the kinetic energy of the marginals)
we have
%\begin{equation}\color{gray}\tag{old}
%\label{ts:KE-bound}
%\KE(\rho-\rho^1_\eps) + \KE(\rho-\rho^2_\eps)
%\leq \frac{C_H}\tau \left( \KE(\rho) + \frac{\delta}{\eps^{1/4}\beta^{1/2}}
%	+ \frac{\beta^{1/2}}{\eps^{1/2}} + \frac{1}{\beta^{1/2}}\right),
%\end{equation}
\begin{equation}
\label{ts:KE-bound}
\KE(\rho-\rho^1_\eps) + \KE(\rho-\rho^2_\eps)
\leq \frac{C_H}\tau \left( \KE(\rho) + \frac1\beta \right),
\end{equation}
\end{enumerate}
Notice that the right hand side of \eqref{ts:PE-bound} is $o_H(\eps^{1/2})$ and the right hand side of \eqref{ts:KE-bound} is $o_H(\eps^{-1/2})$ as $\eps\to0$ because of \eqref{eqn:choice-beta}-\eqref{eqn:choice-tau}.
\end{proposition}

To build the plan $\pi_\eps$, we will first estimate the Wasserstein distance between two suitable measures. We reduce to this estimate since an essential tool of its proof is the Benamou-Brenier formulation of optimal transport, which is not available with a general potential $V$. As a matter of fact fact, even a slightly more precise estimate than the bound on the Wasserstein distance by $o_L(\eps^{1/2})$ is needed: namely, we need a map with  $o_L(\eps^{1/2})$ cost and with a control on the support. Since this second condition is more technical and comes naturally from the proof of \eqref{eqn:wass-est}, the reader may decide to skip it at a first reading and focus on \eqref{eqn:wass-est}.
%{\color{red} Maria aggiunge}
\begin{proposition}\label{prop:pot-cost-rem-mass}
%Let $\rho\in C^1(\setR)$, $\gamma\in\Pi(\rho,\rho)$, $u\in W^{2,\infty}$ and $V(x,y)=\abs{x-y}^{-1}-u(x)-u(y)$ be as in \autoref{thm:rec}.
%Let $\bar\gamma_\eps$ be the plan introduced in \eqref{eq:bar-gamma-eps} and \eqref{eq:bar-gamma-eps-i} and let $\rho^1_\eps = \proj^1_\# \bar\gamma_\eps$ and $\rho^2_\eps = \proj^2_\# \bar\gamma_\eps$ be its two marginals introduced in \eqref{eq:rho-eps-i}.
Let $\Omega_H''$ be the convex hull of $\Omega_H'$. Under the same assumptions as in \autoref{prop:cost-rem-mass}, there exists a constant $C_H$ depending only on $H$ such that
%\old{\begin{equation}
%\label{eqn:wass-est}
%W_2^2\oleft(\rho \bm1_{\Omega_H''}-\rho^2_\eps, {T_\delta}_\#(\rho \bm1_{\Omega_H''}-\rho^1_\eps)\right)
%\leq C_H\left( \frac{\delta^4}{\tau} + \frac{\eps N^2}{\beta\tau}
%	+ \frac{\beta\delta^2}{\tau} \right) %o_L(\sqrt\eps) .
%\end{equation}}
\begin{equation}
\label{eqn:wass-est}
W_2^2\oleft(\rho-\rho^2_\eps, {T_\delta}_\#(\rho-\rho^1_\eps)\right)
\leq C_H\left( \frac{\delta^4}{\tau} + \frac{\eps N^2}{\beta\tau}
	+ \frac{\beta\delta^2}{\tau} \right) %o_L(\sqrt\eps) .
\end{equation}
%Even more precisely, we can build a map $S_\eps$ which maps $\rho-\rho^2_\eps$ to ${T_\delta}_\#(\rho-\rho^1_\eps)$, whose cost is bounded by
%\begin{equation}
%\label{eqn:prop-rhs}
%\int_\R |S(x)-x|^2 (\rho-\rho^2_\eps)\, dx \leq C_H\Big( \frac{\delta^4}{\tau}+ \frac{\eps N^2}{\beta\tau} +  \frac{\eps^{1/2}N}{\beta}\Big) 
%\end{equation}
%\old{Moreover, the optimal map $S_\eps$ between $\rho\bm1_{\Omega_H''}-\rho^2_\eps$ and ${T_\delta}_\#(\rho\bm1_{\Omega_H''}-\rho^1_\eps)$ satisfies a pointwise bound
%\begin{equation}
%\label{eqn:deconv-hp-verif2}
%\abs{S_\eps(x)-x}
%\leq C_H\left(\frac{\delta^4}{\tau^2} + \frac{\eps N^2}{\beta\tau^2}
%	+ \frac{\beta\delta^2}{\tau^2} \right)^{1/3} \bm{1}_{\Omega_H'}(x)
%\qquad \text{for a.e. $x\in\Omega_H''$}.
%\end{equation}}
Moreover, the optimal map $S_\eps$ between $\rho-\rho^2_\eps$ and ${T_\delta}_\#(\rho-\rho^1_\eps)$ satisfies a pointwise bound
\begin{equation}
\label{eqn:deconv-hp-verif2}
\abs{S_\eps(x)-x}
\leq C_H\left(\frac{\delta^4}{\tau^2} + \frac{\eps N^2}{\beta\tau^2}
	+ \frac{\beta\delta^2}{\tau^2} \right)^{1/3} \bm{1}_{\Omega_H'}(x)
\qquad \text{for a.e. $x\in\R$}.
\end{equation}

%\alert{Credo che la strada più semplice sia giustificare nella dimpostrazione che $S_\eps$ ottima per $\rho-\rho^2_\eps$ and ${T_\delta}_\#(\rho-\rho^1_\eps)$ è anche ottima per $\rho\bm1_{\Omega_H''\setminus\Omega_H'}+\rho\bm1_{\Omega_H'}-\rho^2_\eps$ and $\rho\bm1_{\Omega_H''\setminus\Omega_H'}+{T_\delta}_\#(\rho\bm1_{\Omega_H'}-\rho^1_\eps)$.}

In particular, the right-hand sides of \eqref{eqn:wass-est} and  \eqref{eqn:deconv-hp-verif2} are $o_H(\sqrt\eps)$ and $O_H(1)$ respectively as soon as \eqref{eqn:choice-beta} %, \eqref{eqn:choice-delta}, 
and 
\eqref{eqn:choice-tau} are satisfied.
\end{proposition}

We show immediately how \autoref{prop:cost-rem-mass}(i) can be deduced from this estimate.

\begin{proof}[Proof of \autoref{prop:cost-rem-mass}(i) from \autoref{prop:pot-cost-rem-mass}]
Let $\sigma_1 = \rho-\rho^1_\eps$, $\sigma_2 = \rho-\rho^2_\eps$ and let $S_\eps$ be the map given by \autoref{prop:pot-cost-rem-mass} %extended to the identity outside $\Omega_H''$
. %We observe that, since both $\sigma_1$ and $\sigma_2$ coincide with $\rho$ outside $\Omega_H{''}$, this map  satisfies ${S_\eps}_\#{T_\delta}_\#\sigma_1=\sigma_2$.
Let then $\pi_\eps=(\Id, S_\eps\circ T_\delta)_\#\sigma_1\in\Pi(\sigma_1,\sigma_2)$. We observe that (recalling that $T_\delta$ maps $\Omega_H'$ in itself and coincides with $T$ outside, and that $S_\eps (x)= x$ outside $\Omega_H'$)
\begin{equation}
\begin{split}
\abs{S_\eps(T_\delta(x)) - T(x)}
&\leq \abs{S_\eps(T_\delta(x)) - T_\delta(x)} + \abs{T_\delta(x) - T(x)} \\
&\leq \left[C_H\left( \frac{\delta^4}{\tau^2} + \frac{\eps N^2}{\beta\tau^2}
	+ \frac{\beta\delta^2}{\tau^2} \right) + L \delta^2\right] \bm{1}_{\Omega_H'}(x)
\end{split}
\end{equation}
and the right-hand side converges to $0$ as $\eps \to 0$, therefore proving \eqref{eqn:deconv-hp-verif} for $\eps$ sufficiently small. 
Under the assumptions of \autoref{prop:cost-rem-mass} we have that \eqref{eqn:Vgrowth} of \autoref{lem:good-rho} holds, that is
$V(x,y) \leq C_H \dist\bigl((x,y),\graph(T)\bigr)^2$
in a neighborhood (independent of $\eps$) of $\graph(T)$. For $\eps$ sufficiently small, the plan $\pi_\eps$ is supported in a neighborhood of $\graph (T)$ which converges to $0$ as $\eps \to 0$ thanks to \eqref{eqn:deconv-hp-verif2}, so in particular the quadratic upper bound for $V$ applies. Therefore
\[\begin{split}
\int V(x,y) \d\pi_\eps(x,y)
&= \int_{\Omega_H'} V(x,S_\eps(T_\delta(x))) \d\sigma_1(x) \\
&\leq% C\int \abs{y-T_\delta(x)-T(x)+T_\delta(x)}^2 \d\pi_\eps(x,y) \\
%&= 
C_H\int_{\Omega_H'}  \abs{S_\eps(T_\delta(x))-T_\delta(x)-T(x)+T_\delta(x)}^2 \d\sigma_1(x) \\
&\leq 2C_H \int_{\Omega_H'} \abs{S(T_\delta(x))-T_\delta(x)}^2 \d\sigma_1(x)
	+ 2C_H \norm{T-T_\delta}_\infty^2 \\
&\leq 2C_H \int_{\Omega_H'} \abs{S(x)-x}^2 \d({T_\delta}_\#\sigma_1)(x) + 2C_H(L\delta^2)^2 \\
&\leq C_H\left( \frac{\delta^4}{\tau} + \frac{\eps N^2}{\beta\tau}
	+ \frac{\beta\delta^2}{\tau} \right) . \qedhere
%2C o_L(\sqrt\eps) + 2L^2\delta^4
%= %o_L(\sqrt\eps) + \frac{L^2\sqrt\eps}{\abs{\log\eps}^2}
\end{split}
\]

\end{proof}
The estimate \eqref{eqn:wass-est} relies on the particular structure that the marginals of the convolved plan have in the model case when $T$ is (locally) a linear map, which is the main reason why we introduce $T_\delta$ in the construction. We carry out this surprising computation in \autoref{sec:marg-prod} below. {This computation could be interpreted by seeing the convolutions with Gaussian measures as the action of the heat semigroup. However, in our context the Gaussians need to be truncated to have a finite-range interaction, since the control of the potential is only local around its $0$ level set rather than global, and hence this analogy remains only formal and we proceed in a different way.}

It can be seen that $\rho$ and $T_\delta \rho$ are very close in Wasserstein distance, as well as $\rho^2_\eps$ and  ${T_\delta}_\#\rho^1_\eps$. This sole property does not guarantee in general that they remain close when we subtract each other, namely when we consider the measures $\rho-\rho^2_\eps$ and ${T_\delta}_\#(\rho-\rho^1_\eps)$: for instance, when one considers the three measures $\nu_1 = 1_{[0,1]}\mathcal L$, $\nu_2 = 1_{[\alpha,1+\alpha]}\mathcal L$ and $\nu_0 = 1_{[\alpha,1]}\mathcal L$, then $W_2(\nu_1, \nu_2) = \alpha$ but $W_2(\nu_1- \nu_0, \nu_2-\nu_0) =\sqrt \alpha \gg \alpha$ as $\alpha \to0$. However, in our particular situation we can control the distance \eqref{eqn:wass-est} via the Benamou-Brenier formula, since we can provide an interpolating curve between $\rho-\rho^2_\eps$ and ${T_\delta}_\#(\rho-\rho^1_\eps)$ which enjoys some extra properties. For instance, we estimate how close the density of the interpolating curve is to that of the endpoints, providing then a lower bound on this density along the entire curve.

The proof of \autoref{prop:pot-cost-rem-mass} in \autoref{sec:pot-rem} deals then with this Benamou-Brenier estimate and combines it with a linearization argument, which allows to pass from $T$ to its piecewise linear approximation $T_\delta$, and with a suitable control of the errors generated at the interface between two consecutive intervals where $T_\delta$ is linear. A crucial ingredient of the proof is the particular way the two marginals of a truncated Gaussian can be mapped one onto the other at a low cost, which is explained in \autoref{sec:marg-prod}.

Finally, we dedicate \autoref{sec:kin-rem} to the kinetic energy bounds of \autoref{prop:cost-rem-mass}(ii). 

\subsection{Marginals of kernels of product type along a linear map}\label{sec:marg-prod}

In this section we study the properties of the marginals of the rectangularly truncated Gaussian \eqref{eq:gamma-truncated}, which is a convolution kernel of product type
\begin{equation*}%\label{eq:gamma-truncated_}
\Gamma_{M_{\eps,\beta},N}(\x) =
	\frac{\left(e^{-\left(\frac{q}{\eps^{1/2}}+\frac1\beta\right)w^2/2}
		- e^{-N/2}\right)_+^2}
		{G_{\frac{q}{\eps^{1/2}}+\frac1\beta,N}}
	\cdot
	\frac{\left(e^{-\frac1\beta z^2/2} - e^{-N/2}\right)_+^2}
		{G_{\frac1\beta,N}}
= h^1_\eps(w) h^2_\eps(z)
\end{equation*}
where $M_{\eps,\beta}=A/\eps^{1/2}+I/\beta$ as usual, $q$ is the positive eigenvalue of $A$ which is assumed to lie in a certain interval $[1/L, L]$, $z$ is the coordinate in the direction of $\ker(A)$ and $w$ is the transversal coordinate. The change of variables is given by the rotation
\[
(w,z) = R(x,y)
= \left( \frac{y-ax}{\sqrt{1+a^2}}, \frac{ay+x}{\sqrt{1+a^2}} \right)
= (\cos\theta\,y-\sin\theta\,x, \sin\theta\,y+\cos\theta\,x) .
\]
where $(1,a)\in\ker(A)$, hence $\ker(A)=\graph(a\Id+b)$, and $\theta=\arctan(a)$.
Observe that $c_L < \theta < \pi/2 - c_L$ for some constant $c_L>0$ depending only on $L$.

In fact, for the computations of this section the precise form of $h^i_\eps$ does not play a role: we only need $h^1_\eps$ and $h^2_\eps$ to be positive, symmetric $L^1$ functions supported in an interval of length $C_H \eps^{1/4}N^{1/2}$ and $\beta^{1/2}N^{1/2}$ respectively, and with 
\begin{equation}
\label{eqn:deriv h2}\beta^{1/2} \|\partial_z h^2_\eps(z)\|_{L^1} + \beta \|\partial^2_z h^2_\eps(z)\|_{L^1} \leq C .
\end{equation}
This property is satisfied by our kernel since 
\[
(h^2_\eps)'(z)
= -2G_{\frac1\beta,N}^{-1} \frac{z}{\beta}e^{-\frac1\beta z^2/2}
\left(e^{-\frac1\beta z^2/2} - e^{-N/2}\right)_+,
\]
\[
\begin{split}
(h^2_\eps)''(z)
&= 2G_{\frac1\beta,N}^{-1} e^{-\frac1\beta z^2/2}\Big[\Big(-\frac{1}{\beta}+\frac{z^2}{\beta^2} \Big) 
\left(e^{-\frac1\beta z^2/2} - e^{-N/2}\right)_+
%	\spliteq
%	+ 2G_{\frac1\beta,N}^{-1} e^{-\frac1\beta z^2/2}
%	\left(e^{-\frac1\beta z^2/2} - e^{-N/2}\right)_+
%	\spliteq
+ \frac{z^2}{\beta^2}
\bm{1}_{\bigl[-\sqrt{\beta N},\sqrt{\beta N}\bigr]}(z)\Big] ,
\end{split}
\]
from which
%\[
%\abs{(h^2_\eps)'(z)} \leq 2G_{\frac1\beta,N}^{-1} {\abs{z}}
%	e^{-\frac1\beta z^2} \bm{1}_{\bigl[-\sqrt{\beta N},\sqrt{\beta N}\bigr]}(z) .
%\]
\[
\beta^{1/2} \abs{(h^2_\eps)'(z)} + \beta \abs{(h^2_\eps)''(z)}
\leq 2 G_{\frac1\beta,N}^{-1} \left(\frac{\abs{z}}{\beta^{1/2}}+ 1 + 2\frac{z^2}{\beta} \right)
e^{-\frac1\beta z^2} \bm{1}_{\bigl[-\sqrt{\beta N},\sqrt{\beta N}\bigr]}(z);
\]
integrating this inequality and disregarding the last factor $\bm{1}_{\bigl[-\sqrt{\beta N},\sqrt{\beta N}\bigr]}(z)$ we get \eqref{eqn:deriv h2}.

We adopt the following notation to denote rescaled functions by putting a subscript between curly braces
\[
\phi_{\{\zeta\}}(x) = \frac1\zeta\phi\oleft(\frac x\zeta\right)
\]
and we observe that $\norm{\phi_{\{\zeta\}}}_{L^1} = \norm{\phi}_{L^1}$, $(f*g)_{\{\zeta\}} = f_{\{\zeta\}}*g_{\{\zeta\}}$ and $(a\Id+b)_\#\rho(y) = \rho_{\{a\}}(y-b)$.

Let $\eta_\eps^i = \proj^i_\#\Gamma_{M_{\eps,\beta},N}$ be the two marginals of the Gaussian, namely
\begin{align*}
\eta^1_\eps(x)
%= \int_\setR \Gamma_{M_{\eps,\beta},N}(x,y) \d y
&= \int_\setR h^1_\eps(w) h^2_\eps(z) \d y
= \int_\setR h^1_\eps(\cos\theta\,y-\sin\theta\,x) h^2_\eps(\sin\theta\,y+\cos\theta\,x) \d y, \\
\eta^2_\eps(y)
%= \int_\setR \Gamma_{M_{\eps,\beta},N}(x,y) \d x
&= \int_\setR h^1_\eps(w) h^2_\eps(z) \d x
 = \int_\setR h^1_\eps(\cos\theta\,y-\sin\theta\,x) h^2_\eps(\sin\theta\,y+\cos\theta\,x)  \d x,
\end{align*}
and $\tilde\eta_\eps^2=(a\Id)_\#\eta^1_\eps$. We claim that we can write them as a convolution of rescalings of $h^1_\eps$ and $h^2_\eps$, where the only difference between $\eta_\eps^2$ and $\tilde\eta_\eps^2$ lies in the parameter of rescaling of the first function
\begin{align*}
\eta^2_\eps(y) &=
%\bigl((h^1_\eps)_{\cot\theta} * h^2_\eps\bigr)_{\sin\theta}(y)
(h^1_\eps)_{\{\sin\theta\cot\theta\}} * (h^2_\eps)_{\{\sin\theta\}} , \\
\tilde\eta^2_\eps(y) &=
%\bigl((h^1_\eps)_{\tan\theta} * h^2_\eps\bigr)_{\sin\theta}(y)
(h^1_\eps)_{\{\sin\theta\tan\theta\}} * (h^2_\eps)_{\{\sin\theta\}} .
\end{align*}
Indeed, with the change of variable $ \cos\theta\,y-\sin\theta\,x = \frac{t-x}{\sin\theta}$, which rewrites also as $\sin\theta\,y+\cos\theta\,x = \frac{t}{\cos\theta}$ and $\d y = \frac{\d t}{\sin\theta\cos\theta}$,
%\[
%\left\{
%\begin{aligned}
%w &= \cos\theta\,y-\sin\theta\,x = \frac{t-x}{\sin\theta} , \\
%z &= \sin\theta\,y+\cos\theta\,x = \frac{t}{\cos\theta} , \\
%\d y &= \frac{\d t}{\sin\theta\cos\theta},
%\end{aligned}\right.
%\]
and using the fact that $h^1_\eps$ is symmetric, we can compute %the first marginal as
\begin{equation}
\label{eqn:eta-1-repr}
\begin{split}
\eta^1_\eps(x)
%&= \int_\setR \Gamma_{M_{\eps,\beta},N}(x,y) \d y
 %= \int_\setR h^1_\eps(w) h^2_\eps(z) \d y \\
&= \int_\setR \frac1{\sin\theta} h^1_\eps\oleft(\frac{t-x}{\sin\theta}\right)
	\frac1{\cos\theta} h^2_\eps\oleft(\frac{t}{\cos\theta}\right) \d t \\
&= \int_\setR (h^1_\eps)_{\{\sin\theta\}}(x-t) \cdot (h^2_\eps)_{\{\cos\theta\}}(t) \d t
 = \bigl[(h^1_\eps)_{\{\sin\theta\}} * (h^2_\eps)_{\{\cos\theta\}} \bigr](x) .
\end{split}
\end{equation}
In a similar fashion, with another substitution $\cos\theta\,y-\sin\theta\,x = \frac{t}{\cos\theta}$, which implies $\sin\theta\,y+\cos\theta\,x = \frac{y-t}{\sin\theta}$ and $\d x = \frac{\d t}{\sin\theta\cos\theta}$, we find
%\[
%\left\{
%\begin{aligned}
%w &= \cos\theta\,y-\sin\theta\,x = \frac{t}{\cos\theta} , \\
%z &= \sin\theta\,y+\cos\theta\,x = \frac{y-t}{\sin\theta} , \\
%\d x &= \frac{\d t}{\sin\theta\cos\theta},
%\end{aligned}\right.
%\]
%\alert{(qui ci sarebbe un segno $-$ in $\d x$, che però non conta perché si prende il modulo del Jacobiano)}
%we can compute the second marginal as
\begin{equation}
\label{eqn:eta-2-repr}
\begin{split}
\eta^2_\eps(y)
&= \int_\setR \Gamma_{M_{\eps,\beta},N}(x,y) \d x
 = \int_\setR h^1_\eps(w) h^2_\eps(z) \d x \\
&= \int_\setR \frac1{\cos\theta} h^1_\eps\oleft(\frac{t}{\cos\theta}\right)
	\frac1{\sin\theta} h^2_\eps\oleft(\frac{y-t}{\sin\theta}\right) \d t \\
&= \int_\setR (h^1_\eps)_{\{\cos\theta\}}(t) \cdot (h^2_\eps)_{\{\sin\theta\}}(y-t) \d t
 = \bigl[(h^1_\eps)_{\{\cos\theta\}} * (h^2_\eps)_{\{\sin\theta\}} \bigr](y) .
\end{split}
\end{equation}
Finally, by the properties of the rescalings, we have that
\[
\begin{split}
\tilde\eta^2_\eps
= (\tan\theta\Id)_\#\eta^1_\eps = (\eta^1_\eps)_{\{\tan\theta\}} 
%\\
%&= \bigl[(h^1_\eps)_{\{\sin\theta\}}
%	* (h^2_\eps)_{\{\cos\theta\}} \bigr]_{\{\tan\theta\}}
= (h^1_\eps)_{\{\sin\theta\,\tan\theta\}} * (h^2_\eps)_{\{\sin\theta\}}.
\end{split}
\]

A simple consequence of these representations of $\eta^i_\eps$ is for instance an $L^\infty$ bound which will be useful later: from \eqref{eqn:eta-1-repr}, since $\|(h^1_\eps)_{\{\sin\theta\}}\|_{L^1}= \|h^1_\eps\|_{L^1}= 1$ and by \eqref{eq:Gainfty-GaN} we have
\begin{equation}
\label{eqn:eta-l-infty}
\|\eta^1_\eps\|_{L^\infty}
%&= \int_\setR \Gamma_{M_{\eps,\beta},N}(x,y) \d y
 %= \int_\setR h^1_\eps(w) h^2_\eps(z) \d y \\
 \leq \|(h^1_\eps)_{\{\sin\theta\}}\|_{L^1} \|(h^2_\eps)_{\{\cos\theta\}}\|_{L^\infty} \leq C_H \|h^2_\eps \|_{L^\infty} \leq \frac{C_H}{G_{\frac1\beta,N}} \leq \frac{C_H}{G_{\frac1\beta,\infty}} =  \frac{C_H}{\beta^{1/2}}  .
\end{equation}

%because the rescaling and the convolution commute.

The following proposition considers a plan $\lambda$ supported on the graph of $(a\Id+b)$ and introduces the convolved plan $\lambda \ast \Gamma_{M_{\eps,\beta},N}$, its marginals $\rho^1_\eps$ and $\rho^2_\eps$  and the pushforward $\tilde \rho^2_\eps$  of the first marginal through the linear map $(a\Id+b)$. It uses the representation of the marginals in terms of convolutions of $\eta_\eps^i$ with the marginal of $\lambda$ to provide estimates on their distance and an interpolating curve between $ \rho^2_\eps$  and $\tilde \rho^2_\eps$.

\begin{proposition}\label{prop:conv-prod}
Let $a\in(L^{-1},L)$, $\theta=\arctan(a)$, $b\in\setR$, $A$ a degenerate positive-semidefinite symmetric $2\times2$ matrix with a positive eigenvalue $q\in(L^{-1},L)$, $\ker(A)=\graph(a\Id)$ and $M_{\eps,\beta}=A/\eps^{1/2}+I/\beta$.

Given $\lambda\in\Meas_+(\setR^2)$ with $\supp\lambda\subset\graph(a\Id+b)$, define $\lambda_\eps = \lambda * \Gamma_{M_{\eps,\beta},N}$, let $\rho^i = \proj^i_\#\lambda$, $\rho^i_\eps = \proj^i_\#\lambda_\eps$ for $i=1,2$ and define $\tilde\rho^2_\eps = (a\Id+b)_\#\rho^1_\eps$.
Then we have that
\[
\mu_t = (h^1_\eps)_{\{\sin\theta\,(\tan\theta)^{2t-1}\}}
	* (h^2_\eps)_{\{\sin\theta\}} * \rho^2, \qquad t\in[0,1],
\]
is a curve interpolating between $\rho^2_\eps$ and $\tilde\rho^2_\eps$ satisfying $\partial_t \mu_t = \div(m_t)$ with
\begin{equation}
\label{ts:l-infty-rho}
\norm{\mu_t - \rho^2_\eps}_\infty
\leq CL \frac{\eps^{1/2}N}{\beta} \norm{\rho^2}_\infty,
\end{equation}
$\mu_t$ and $m_t$ supported on the convex hull of $\supp \tilde \rho^2_\eps \cup \supp\rho^2_\eps$ and
\begin{equation}
\label{eqn:mt-bound}
\abs{m_t} \leq C_H \frac{\eps^{1/2} N}{\beta^{1/2}} \norm{\rho^2}_\infty.
\end{equation}

%{\color{gray}
%Suppose that we put some measure $\rho^1$ on a little interval where $T$ is linear and let $\rho^2=T_\#\rho^1$. Then we consider $\gamma=(\Id,T)_\#\rho^1$ and $\gamma_\eps = \gamma * \Gamma_{M_{\eps,\beta},N}$. The marginals of $\gamma_\eps$ are $\rho^1_\eps = \eta^1_\eps*\rho^1$, $\rho^2_\eps = \eta^2_\eps*\rho^2$. Define
%\[
%\tilde\rho^2_\eps = T_\#\rho^1_\eps = (\rho^1_\eps)_{a}
%= (\eta^1_\eps*\rho^1)_{a}=(\eta^1_\eps)_{a}*(\rho^1_\eps)_{a}
%= T_\#\eta^1_\eps * T_\#\rho^1_\eps = \tilde\eta^2_\eps*\rho^2_\eps .
%\]
%Letting $\tilde h^i_\eps = (h^i_\eps)_{\sin\theta}$, we have
%\begin{align*}
%\rho^2_\eps &= (\tilde h^1_\eps)_{\cot\theta} * \tilde h^2_\eps * \rho^2 \\
%\tilde\rho^2_\eps &= (\tilde h^1_\eps)_{\tan\theta} * \tilde h^2_\eps * \rho^2.
%\end{align*}
%Then 
%\[
%\mu_t = (\tilde h^1_\eps)_{(\tan\theta)^{2t-1}} * \tilde h^2_\eps * \rho^2
%\]
%is a curve interpolating between $\rho^2_\eps$ and $\tilde\rho^2_\eps$, satisfying $\partial_t \mu_t = \div(m_t)$ with
%\[
%\norm{\rho^2_\eps - \tilde\rho^2_\eps}_\infty
%%\leq C \frac{\eps^{1/2}N}{\beta^{1/2}} \Lip(\rho)
%\leq C \frac{\eps^{1/2}N}{\beta} \norm{\rho}_\infty
%\]
%and
%\[
%\abs{m_t} \leq \frac{C H \sqrt\eps m_2(h^1)}{\sqrt{\beta}} \norm{\rho^2}_\infty,
%\]
%where $T',1/T'<H$.
%}
\end{proposition}

%{\color{gray}
%In particular we have that $h^{\ep}_1$ solves at some approximation a differenential equation (depending on the point $x,y$) at scale $\sqrt[4]{\ep}$, while we chose $h^{\ep}_2$ having scale $\sqrt{\beta(\ep)}$.
%In order to preserve the fact that $\Gamma_{M_{\eps,\beta},N}$ is a probability density, we can assume that for some fixed probability densities $h_1$ and $h_2$ we have
%\begin{align*}
%h^1_\eps(w) &= \sqrt{\frac{q}{\eps^{1/2}}+\frac1\beta}
%	h_1\oleft(\sqrt{\frac{q}{\eps^{1/2}}+\frac1\beta}w\right) &
%h^2_\eps &= \frac1{\sqrt{\beta}} h_2\oleft( \frac{z}{\sqrt{\beta}}\right).
%\end{align*}
%}

\begin{proof}

Up to a translation, we may assume that $b=0$.  %We now turn to the study of the marginals of the convolution $\lambda_\eps$.
We have  that $\rho^i_\eps = \eta^i_\eps \ast \rho^i$ and, recalling that the kernel $h^2_\eps$ is symmetric, 
\[
\begin{split}
\mu_t(y) - \rho^2_\eps(y)
&= \int_\setR \bigl[(h^1_\eps)_{\{\sin\theta\,(\tan\theta)^{2t-1}\}}
	- (h^1_\eps)_{\{\sin\theta\cot\theta\}}\bigr](y-z)
	\cdot \bigl((h^2_\eps)_{\{\sin\theta\}} * \rho^2\bigr)(z) \d z \\
&= \int_\setR \bigl[(h^1_\eps)_{\{\sin\theta\,(\tan\theta)^{2t-1}\}}
	- (h^1_\eps)_{\{\sin\theta\cot\theta\}}\bigr](y-z)
		\spliteq\qquad
	\cdot \biggl[\bigl((h^2_\eps)_{\{\sin\theta\}} * \rho^2\bigr)(y)
	+ \bigl((h^2_\eps)_{\{\sin\theta\}} * \rho^2\bigr)'(y) \cdot (z-y)
		\spliteq\qquad\qquad
	+ \frac12\bigl((h^2_\eps)_{\{\sin\theta\}} * \rho^2\bigr)''(y_z) \cdot (z-y)^2
	\biggr] \d z \\
&= \int_\setR \bigl[(h^1_\eps)_{\{\sin\theta\,(\tan\theta)^{2t-1}\}}
	- (h^1_\eps)_{\{\sin\theta\cot\theta\}}\bigr](y-z)
		\spliteq\qquad
	\cdot \frac{\bigl((h^2_\eps)_{\{\sin\theta\}} * \rho^2\bigr)''(y_z)}2
	\cdot (z-y)^2 \d z,
\end{split}
\]
%\[\color{gray}
%\begin{split}
%\rho^2_\eps(y) - \tilde\rho^2_\eps(y)
%&= \int \bigl[( h^1_\eps)_{\{\sin\theta\cot\theta\}} - ( h^1_\eps)_{\{\sin\theta\tan\theta\}}\bigr](y-z)
%	\cdot \bigl((h^2_\eps)_{\{\sin\theta\}} * \rho^2\bigr)(z) \d z \\
%&= \int \bigl[( h^1_\eps)_{\{\sin\theta\cot\theta\}} - ( h^1_\eps)_{\{\sin\theta\tan\theta\}}\bigr](y-z)
%%\hspace{2cm}
%	\cdot\left[
%	\bigl((h^2_\eps)_{\{\sin\theta\}} * \rho^2\bigr)(y)\right.
%\\ &\hspace{2em}\left.	+ \bigl((h^2_\eps)_{\{\sin\theta\}} * \rho^2\bigr)'(y) \cdot (z-y)
%	+ \frac12\bigl((h^2_\eps)_{\{\sin\theta\}} * \rho^2\bigr)''(y_z) \cdot (z-y)^2
%	\right] \d z \\
%&= \int \bigl[( h^1_\eps)_{\sin\theta\cot\theta} - ( h^1_\eps)_{\sin\theta\tan\theta}\bigr](y-z)
%	\cdot \frac{\bigl((h^2_\eps)_{\{\sin\theta\}} * \rho^2\bigr)''(y_z)}2 \cdot (z-y)^2 \d z,
%\end{split}
%\]
where $y_z$ is a point in the segment between $y$ and $z$. %for every $z$.
Therefore, since the integrand is nonzero only for $|y-z|^2 \leq C_H \eps^{1/2}N$ and thanks to \eqref{eqn:deriv h2}, we obtain \eqref{ts:l-infty-rho}
\[
\begin{split}
\norm{\mu_t - \rho^2_\eps}_\infty
&\leq C_H\norm*{(h^1_\eps)_{\{\sin\theta\,(\tan\theta)^{2t-1}\}}
	- (h^1_\eps)_{\{\sin\theta\cot\theta\}}}_1
	\cdot\norm*{[(h^2_\eps)_{\{\sin\theta\}}  * \rho^2]''}_\infty \eps^{1/2}N \\
%&\leq \frac12 2 \norm{(\tilde h^2_\eps)'}_1 \cdot \norm{(\rho^2)'}_\infty
%\leq C \frac{\eps^{1/2}N}{\beta^{1/2}} \Lip(\rho^2) .
&\leq C_H \eps^{1/2}N  \norm{[(h^2_\eps)_{\{\sin\theta\}}]''}_1 \cdot \norm{\rho^2}_\infty
\leq C_H \frac{\eps^{1/2}N}{\beta} \norm{\rho^2}_\infty .
\end{split}
\]
%This concludes the $L^\infty$ estimate.

Define $H^1_\eps(x) = \int_{-\infty}^x rh^1_\eps(r)\d r$, so that $(H^1_\eps)_{\{\zeta\}}(x)$ is the primitive of $\frac{x}{\zeta^2}(h^1_\eps)_{\{\zeta\}}(x)$ with the same compact support.
%Notice that when we differentiate with respect to the parameter $\zeta$ we have
This allows us to express the derivative of $(h^1_\eps)_{\{\zeta\}}(x)$ with respect to the parameter $\zeta$ as
\[
\begin{split}
\frac{\partial(h^1_\eps)_{\{\zeta\}}(x)}{\partial\zeta}
&= -\frac1{\zeta^2}h^1_\eps\oleft(\frac x\zeta\right)
	- \frac x{\zeta^3}(h^1_\eps)'\oleft(\frac x\zeta\right)
= -\frac{\partial}{\partial x}\oleft(\frac x{\zeta^2}h^1_\eps\oleft(\frac x\zeta\right)\right) \\
&= -\frac{\partial}{\partial x}\oleft(\zeta \frac x{\zeta^2}(h^1_\eps)_{\{\zeta\}}(x)\right) 
= -\zeta \frac{\partial^2}{\partial x^2}\oleft((H^1_\eps)_{\{\zeta\}}(x)\right) ,
\end{split}
\]
Let now $\zeta(t) = \sin \theta (\tan\theta)^{2t-1}$. We deduce that
\[
\begin{split}
\partial_t\mu_t
&= \frac{\partial( h^1_\eps)_{\{\zeta(t)\}}}{\partial t} * \tilde h^2_\eps * \rho^2 \\
&= \zeta'(t) \frac{\partial( h^1_\eps)_{\{\zeta(t)\}}}{\partial\zeta}
	* \tilde h^2_\eps * \rho^2
= -\zeta'(t) \zeta(t) [(H^1_\eps)_{\{\zeta(t)\}}]''
	* \tilde h^2_\eps * \rho^2 \\
&= -\zeta'(t) \zeta(t) \div\oleft( (H^1_\eps)_{\{\zeta(t)\}}
	* [(h^2_\eps)_{\{\sin\theta\}}]' * \rho^2 \right)
= \div(m_t) .
\end{split}
\]
We can bound by \eqref{eqn:deriv h2}
\[
\norm{m_t}_\infty
\leq C_H \norm{ (H^1_\eps)_{\zeta(t)}
}_1 \cdot \norm{[(h^2_\eps)_{\{\sin\theta\}}]'}_1
	\cdot \norm{\rho^2}_\infty
\leq C_H \norm{ H^1_\eps}_1 \cdot \frac{C}{\beta^{1/2}}
	\cdot \norm{\rho^2}_\infty
%\frac{C H \sqrt\eps m_2(h^1)}{\sqrt{\beta}} \norm{\rho^2}_\infty
\]
and estimate the first factor in the right-hand side in terms of %his scaling in 
powers of $\eps$ using that $h^1_\eps$ is supported in a set of size $C_H  \eps^{1/4} N^{1/2}$
%(we assume $h_1$ to have support in $[-C,C]$)
\[
\begin{split}
\norm{H^1_\eps}_1
&\leq -2\int_{-\infty}^0 \int_{-\infty}^x r h^1_\eps(r)\d r \d x
= -2\int_{-\infty}^0 \int_{r}^0 {r}h^1_\eps(r)\d x \d r \\
&= 2\int_{-\infty}^0 {r^2}h^1_\eps(r) \d r
\leq C_H \eps^{1/2}N \norm{h^1_\eps}_1%H \int_{-\infty}^\infty r^2 (h^1_\eps)_{\{\sin\theta\}} (r) \d r
%= H \sqrt\eps \int_{-\infty}^\infty r^2 (h^1_\eps)_{\{\sin\theta\}} (r) \d r
\leq C_H \eps^{1/2}N.
\end{split}
\]
%and
%\[
%\norm*{\frac{\partial\tilde h^2_\eps}{\partial r}}_1
%\leq \frac{C}{\sqrt{\beta}} .
%\]
%so
%$$|m_t(x)| \leq |1-a| \cdot (t+(1-t)a) \cdot \sqrt{\ep}m_2(h_1) \cdot \sup_{ |y-x| \leq C \sqrt[4]{\ep} } \{ \frac d {dx} ( \beta \star \rho_2)\}.$$
Hence we obtain \eqref{eqn:mt-bound}.
%
%{\color{gray}In particular if we are in a region where $|T'| \leq H$ for some $H>1$, and $\rho_2$ is Lipschitz, then
%
%$$|m_t(x)| \leq \sqrt{\ep} H^2 m_2(h_1) \cdot Lip(\rho_2).$$
%
%Another softer estimate comes from the decrasingness of the $2$-norm with respect to the convolution with respect to some probability kernel; normalizing we get
%
%\[
%\int m_t(x)^2  \, dx \leq (1-a)^2 \cdot (t+(1-t)a)^2 \cdot \ep \cdot m_2(h_1)^2 \cdot  \int  d {dx} ( \beta \star \rho_2)^2 \,d x
%\]}
\end{proof}

\subsection{Wasserstein estimate of \autoref{prop:pot-cost-rem-mass}}\label{sec:pot-rem}

\begin{proof}[Proof of \autoref{prop:pot-cost-rem-mass}]
Define the map $T_t=(1-t)T+tT_\delta$ and define $T_\delta^i$ to be the affine extension of ${T_\delta}\rvert_{I^1_i}$. Define also the marginals $\rho^1_{\eps,i}$ and $\rho^2_{\eps,i}$ of $\bar\gamma_{\eps,i}$.
We introduce the following three curves of measures:
\begin{itemize}
\item $\rho_t = {T_t}_\#\rho$ from $T_\#\rho=\rho$ to ${T_\delta}_\#\rho$;
\item $\xi^i_t$ from $\rho^2_{\eps,i}$ to ${T_\delta^i}_\#\rho^1_{\eps,i}$ using the construction presented in \autoref{prop:conv-prod};
\item $\zeta^i_t$ from ${T_\delta^i}_\#\rho^1_{\eps,i}$ to ${T_\delta}_\#\rho^1_{\eps,i}$ defined by linearly stretching $\rho^1_{\eps,i}\res{I^1_{i+1}}$ from $T^i_\delta$ to $T^{i+1}_\delta$ and similarly for $\rho^1_{\eps,i}\res{I^1_{i-1}}$ with \autoref{lem:stretching}.
\end{itemize}
We transform the remaining mass in three steps:
\[
\rho - \rho^2_\eps
	\xto{\rho_t-\rho^2_\eps}
{T_\delta}_\#\rho - \rho^2_\eps
	\xto{{T_\delta}_\#\rho - \sum_i \xi^i_t}
{T_\delta}_\#\rho - \sum_i {T^i_\delta}_\#\rho^1_{\eps,i}
	\xto{{T_\delta}_\#\rho - \sum_i \zeta^i_t}
{T_\delta}_\#\rho - \sum_i {T_\delta}_\#\rho^1_{\eps,i} .
\]
We estimate the $W_2$ distance of each step with Benamou-Brenier.

\textit{Step 1.} For every $x\in \R$ we observe that the map $t \to T_t(x)$ solves the ODE 
\[
\partial_t X(t) = (T_\delta-T) \circ T_t^{-1} \circ X(t)
\]
with initial datum $X(0)= x$. Notice that $(T_\delta-T) \circ T_t^{-1}$ is a locally Lipschitz vector field, therefore $T_t(x)$ represents its flow. As a consequence $\rho_t = (T_t)_\# \rho$ solves
\[
\partial_t\rho_t=-\div\bigl((T_\delta-T)\circ T_t^{-1}\rho_t\bigr).
\]

 From the definition of $\rho_t$ we get that $\rho_t(T_tx) = \rho(x)/T_t'(x)$, whereas $\rho(Tx) = \rho(x)/T'(x)$, since $T_\#\rho=\rho$. These can be rewritten as
\begin{align*}
\rho(y) &= \frac{\rho\bigl(T^{-1}(y)\bigr)}{T'\bigl(T^{-1}(y)\bigr)}, &
\rho_t(y) &= \frac{\rho\bigl(T_t^{-1}(y)\bigr)}{T_t'\bigl(T_t^{-1}(y)\bigr)}.
\end{align*}
Our goal is to show that these two densities are close to one another for $y\in \Omega_H$. We start by comparing the denominators.
 Let $x\in I^1_i$ and $y=T(x)$; we have $T_t^{-1}(y)\in I^1_i$. Then $T'\bigl(T^{-1}(y)\bigr) = T'(x) \in T'(I^1_i)$; but also $T_t'\bigl(T_t^{-1}(y)\bigr) = (1-t)T'\bigl(T_t^{-1}(y)\bigr)+tT_\delta'\bigl(T_t^{-1}(y)\bigr) \in \mathop{\mathrm{conv}} \bigl(T'(I^1_i)\bigr) = T'(I^1_i)$. Therefore we have
\[
\abs*{T_t'\bigl(T_t^{-1}(y)\bigr) - T'\bigl(T^{-1}(y)\bigr)}
\leq \Lip(T') \diam(I^1_i) \leq L\delta,
\]
\begin{equation}
\label{eqn:min}
\min\{T_t'\bigl(T_t^{-1}(y)\bigr) : y \in \Omega_{H}\} \geq \min\{T'\bigl(T_t^{-1}(y)\bigr) : y \in \Omega_{H}\} \geq \frac 1 L.
\end{equation}
Let us now estimate the numerators of the densities.
 Differentiating with respect to $t$ the identity $y = T_t\bigl(T_t^{-1}(y)\bigr)$ and using the fact that $\d T_t/\d t = T_\delta-T$ we get
\[
0 = \frac{\d T_t}{\d t}\bigl(T_t^{-1}(y)\bigr)
	+ T_t'\bigl(T_t^{-1}(y)\bigr) \frac{\d T_t^{-1}}{\d t}(y),
\]
from which we deduce by  \eqref{eqn:Tdelta-Tclose} and \eqref{eqn:min} that for $y \in \Omega_H$
%\[
%\frac{\d T_t^{-1}}{\d t}(y)
%= \frac{(T-T_\delta)\bigl(T_t^{-1}(y)\bigr)}{T_t'\bigl(T_t^{-1}(y)\bigr)}
%\]
%and the estimate
\[
\abs*{\frac{\d T_t^{-1}}{\d t}(y)} = \abs*{\frac{(T-T_\delta)\bigl(T_t^{-1}(y)\bigr)}{T_t'\bigl(T_t^{-1}(y)\bigr)}
}
\leq L^2 \delta^2.
\]
This allows us to compute for $y \in \Omega_H$ (notice that for $y \in (\Omega_H)^c$ the quantity below is $0$)
\[
\abs*{T_t^{-1}(y) - T^{-1}(y)} = \abs*{T_t^{-1}(y) - T_0^{-1}(y)}
\leq \int_0^t \abs*{\frac{\d T_t^{-1}}{\d t}(y)} \d t
\leq L^2 \delta^2,
\]
therefore
\[
\abs*{\rho\bigl(T_t^{-1}(y)\bigr) - \rho\bigl(T^{-1}(y)\bigr)}
\leq \Lip(\rho) L^2 \delta^2 \leq L^3 \delta^2 .
\]
From this we conclude that
\[
\begin{split}
\abs{\rho_t(x)- {\rho(x)}}
&\leq \frac{\abs{\rho\bigl(T^{-1}(y)\bigr)-\rho\bigl(T_t^{-1}(y)\bigr)}}
		{T'\bigl(T^{-1}(y)\bigr)}
	+ \rho\bigl(T_t^{-1}(y)\bigr)
		\frac{\abs{T_t'\bigl(T_t^{-1}(y)\bigr)-T'\bigl(T^{-1}(y)\bigr)}}
		{T_t'\bigl(T_t^{-1}(y)\bigr)T'\bigl(T^{-1}(y)\bigr)} \\
&\leq L^4\delta^2 + L^4\delta
\leq C_H \delta.
\end{split}
\]
%\[
%1-CL^3\delta \leq \frac{\rho_t(x)}{\rho(x)} \leq 1+CL^3\delta
%\]
%for some constant $C$.
On the other hand, by \eqref{eq:lower-marginal} of \autoref{lem:approx-plan} we have that $\rho^2_\eps \leq \rho - c_H\tau$. % Combining this with 
%\[
%\rho^1_\eps(x)
%= \sum_i \Bigl(\bigl[\bigl(\rho-f(\eps)\bigr)\chi_{I^1_i}\bigr]*k^1_i\Bigr)(x),
%\]
%where the $k^1_i$ are convolution kernels of size $\sqrt{\beta(\eps)N(\eps)}$; therefore, by the H\"olderianity of $\rho$, we have (if $\rho>f(\eps)$)
%\[
%\rho^1_\eps(x) \leq \rho(x)-\tau
%	+ \norm{\rho}_{C^{0,1/2}}\bigl(\beta(\eps)N(\eps)\bigr)^{1/4}.
%\]
Combined together, these two density estimates imply that
\[
C_H \geq \rho_t(x) \geq \rho^2_\eps(x) + c_H\tau/2
\qquad \text{for every $x \in \Omega_H$}.
\]
Therefore, with Benamou-Brenier % and \eqref{eqn:Tdelta-Tclose}
and noticing that  $T_\delta(y)-T(y) =0$ for $y \in (\Omega_H)^c$, we can compute
\[
\begin{split}
W_2^2(\rho-\rho^2_\eps, {T_\delta}_\#\rho-\rho^2_\eps)
&\leq \int_0^1 \int_{\Omega_H'} \frac1{\rho_t(x) - \rho^2_\eps(x)}
	\abs*{(T_\delta-T)\circ T_t^{-1}(x) \rho_t(x)}^2 \d x \d t \\
&\leq \int_0^1 \int_{\Omega_H'} \frac{C_H}{\tau}
	\norm{T_\delta-T}_\infty^2 \norm{\rho}_\infty^2 \d x \d t
\leq \frac{C_H\delta^4}{\tau},
\end{split}
\]
where in the last step we use the estimate \eqref{eqn:Tdelta-Tclose}.

\textit{Step 2.}  By \autoref{prop:conv-prod} applied with  $\lambda = (\Id,T_\delta)_\#\bigl((\rho-\tau)\bm{1}_{I^1_i}\bigr)$ (recall that by \eqref{eq:bar-gamma-eps-i} with this definition $\lambda_\eps$ has marginals $\rho^1_{\eps,i}$ and $\rho^2_{\eps,i}$)
for every $i$ there exists a curve $\xi^i_t$ connecting  $\rho^1_{\eps,i}$ and $\rho^2_{\eps,i}$ such that $\partial_t \xi^i_t = \div (m^i_t)$ and
\[
\norm{\xi^i_t - \rho^2_{\eps,i}}_\infty = \norm{\xi^i_t - \xi^i_0}_\infty
\leq C_H \frac{\eps^{1/2}N}{\beta},
\]
\[
\abs{m^i_t} \leq C_H \frac{\eps^{1/2} N }{\beta^{1/2}}.
\]
The second curve to consider is ${T_\delta}_\#\rho - \sum_i \xi^i_t$, which solves $\partial_t({T_\delta}_\#\rho-\sum_i\xi^i_t) = \div(\sum_i m^i_t)$ and connects ${T_\delta}_\#\rho - \rho^2_\eps $ and $
{T_\delta}_\#\rho - \sum_i {T^i_\delta}_\#\rho^1_{\eps,i} $. Since both $\xi^i_t$ and  $m^i_t$ are supported in a small neighborhood of $I^i_\delta$, at each point $x$ at most two of each of these objects can overlap. Therefore for a.e. $x\in \Omega_H$
\[
\abs*{\sum_i\xi^i_t - \rho^2_\eps} + \abs*{\sum_i m^i_t}
\leq 4 C_H \frac{\eps^{1/2}N}{\beta^{1/2}}
\]
because the terms overlap at most twice. This means again that the density of ${T_\delta}_\#\rho - \sum_i \xi^i_t$ is at least $C_H\tau/2$ and 
\[
\begin{split}
W_2^2\oleft({T_\delta}_\#\rho-\rho^2_\eps,
	{T_\delta}_\#\rho-\sum_i{T^i_\delta}_\#\rho^1_{\eps,i}\right)
&\leq \int_0^1 \int_{\Omega_H'} \frac1{{T_\delta}_\#\rho - \sum_i \xi^i_t}
	\abs*{\sum_i m^i_t}^2 \d x \d t \\
&\leq \int_0^1 \int_{\Omega_H'}  C_H \frac{\eps N^2}{\beta\tau} \d x \d t \leq C_H\frac{\eps N^2}{\beta\tau} .
\end{split}
\]

\textit{Step 3}.
Finally, we deal with the curve ${T_\delta}_\#\rho - \sum_i \zeta^i_t$ which connects ${T_\delta}_\#\rho - \sum_i {T^i_\delta}_\#\rho^1_{\eps,i}$ to ${T_\delta}_\#\rho - \sum_i {T_\delta}_\#\rho^1_{\eps,i}$.
We split each term $\rho^1_{\eps,i}$ as the sum of three contributions
\[
\rho^1_{\eps,i} = \rho^1_{\eps,i}\bm{1}_{I^1_{i-1}} + \rho^1_{\eps,i}\bm{1}_{I^1_i} + \rho^1_{\eps,i}\bm{1}_{I^1_{i+1}} .
\]
On the interval $I^1_i$ the maps $T^i_\delta$ and $T_\delta$ coincide, therefore we leave the central mass $\rho^1_{\eps,i}\bm{1}_{I^1_i}$ still. On the interval $I^1_{i+1}$, on the other hand, the map $T_\delta$ coincides with $T^{i+1}_\delta$, so the mass $\rho^1_{\eps,i}\bm{1}_{I^1_{i+1}}$ must be transformed from ${T^i_\delta}_\#(\rho^1_{\eps,i}\bm{1}_{I^1_{i+1}})$ to ${T^{i+1}_\delta}_\#(\rho^1_{\eps,i}\bm{1}_{I^1_{i+1}})$.
We perform this transformation by applying \autoref{lem:stretching} to $\mu = \rho^1_{\eps,i}\bm{1}_{I^1_{i+1}}$, $x_0$ equal to the point separating $I^1_i$ and $I^1_{i+1}$, $\ell=\sqrt{\beta N}$, $S_{\lambda_0} = T^i_\delta$ and $S_{\lambda_1} = T^{i+1}_\delta$. Similarly, we need to transform the mass in $I^1_{i-1}$ from ${T^i_\delta}_\#(\rho^1_{\eps,i}\bm{1}_{I^1_{i-1}})$ to ${T^{i-1}_\delta}_\#(\rho^1_{\eps,i}\bm{1}_{I^1_{i-1}})$, and we do so with a mirrored version of \autoref{lem:stretching}. In summary, the curve $\zeta^i_t$ is composed of three terms: a central mass $\rho^1_{\eps,i}\bm{1}_{I^1_i}$ which does not move, and two masses one on each side which move according to the construction in \autoref{lem:stretching}.
Notice that, when summing over $i$, all the contributions in $\sum_i\zeta^i_t$ deriving from the application of \autoref{lem:stretching} are disjoint, because there are two terms in every $I^1_i$ but they don't overlap (since $\sqrt{\beta N}\ll\delta$).

The lower bound on the density of ${T_\delta}_\#\rho-\sum_i\zeta^i_t$ follows from the same estimate as in the first step and the estimate of the density in \autoref{lem:stretching}. Moreover, \autoref{lem:stretching} also gives us the estimate on the momentum
\[
\abs{m^i_t} \leq \beta^{1/2}\delta \frac{\max T'}{(\min T')^2} \leq L^3 \beta^{1/2}\delta,
\]
which again do not overlap. Therefore we can compute
\[
\begin{split}
W_2^2\oleft({T_\delta}_\#\rho-\sum_i{T^i_\delta}_\#\rho^1_{\eps,i},
	{T_\delta}_\#\rho-\sum_i{T_\delta}_\#\rho^1_{\eps,i}\right)
&\leq \int_0^1 \int_{\Omega_H'} \frac1{{T_\delta}_\#\rho-\sum_i\zeta^i_t}
	\abs*{\sum_i m^i_t}^2 \d x \d t \\
&\leq \int_0^1 \int_{\Omega_H'} \frac{C_H}{\tau} \beta \delta^2 \d x \d t
\leq C_H \frac{\beta\delta^2}{\tau} .
\end{split}
\]

%Therefore the thesis is correct if we require the first three inequalities in \eqref{eqn:choice-tau}, namely
%\begin{align*}
%\tau &\gg \frac{\delta^4}{\eps^{1/2}} , &
%\tau &\gg \frac{\eps^{1/2}N}{\beta} , &
%\tau &\gg \frac{\beta\delta^2}{\eps^{1/2}} . \qedhere
%\end{align*}

\textit{Step 4: proof of \eqref{eqn:deconv-hp-verif2}.}
We apply the main estimate in \cite[Proposition 1.2]{BJR07}: if $\mu_0, \mu_1$ are probability measures on a bounded interval $I_0$, with  $\mu_0$  absolutely continuous and with density bounded below by $\tau_0$, then the optimal transport map $S$ between $\mu_0$ and $\mu_1$ satisfies the estimate
\begin{equation}\label{eqn:BJR}
\norm{S(x) - x}_{L^\infty}^3 \leq \frac{C_\ell}{\tau_0} W_2^2(\mu_0, \mu_1).
\end{equation}

Let $S_\eps$ be the optimal map from
\[
\rho - \rho^2_\eps
= \rho\bm1_{\Omega_H'} - \rho^2_\eps + \rho\bm1_{\Omega_H''\setminus\Omega_H}
	+ \rho\bm1_{\setR\setminus\Omega_H''}
\]
to
\[
{T_\delta}_\#(\rho-\rho^1_\eps)
= {T_\delta}_\#(\rho\bm1_{\Omega_H'}-\rho^1_\eps) + \rho\bm1_{\Omega_H''\setminus\Omega_H}
	+ \rho\bm1_{\setR\setminus\Omega_H''} .
\]
Since the tails outside $\Omega_H''$ of both the marginals coincide with $\rho\bm1_{\setR\setminus\Omega_H''}$, we claim that $S_\eps(x)=x$ for $x\in\setR\setminus\Omega_H''$ and that $S_\eps$ is also an optimal map from
\[
\mu_0 = (\rho - \rho^2_\eps)\bm1_{\Omega_H''}
= (\rho\bm1_{\Omega_H'}-\rho^2_\eps) + \rho\bm1_{\Omega_H''\setminus\Omega_H}
\]
to
\[
\mu_1 = {T_\delta}_\#(\rho-\rho^1_\eps)\bm1_{\Omega_H''}
= {T_\delta}_\#(\rho\bm1_{\Omega_H'}-\rho^1_\eps) + \rho\bm1_{\Omega_H''\setminus\Omega_H}.
\]
Applying \eqref{eqn:BJR} with $I_0=\Omega_H''$ and $\tau_0 = c_H\tau$ (notice that in the gap $\Omega_H''\setminus\Omega_H'$ we have $\mu_0=\rho\gg\tau$) we deduce that
\[
\begin{split}
\norm{S_\eps(x) - x}_{L^\infty(\Omega_H'')}^3
&\leq \frac{C_H}{\tau} W_2^2\oleft((\rho - \rho^2_\eps)\bm1_{\Omega_H''},
	{T_\delta}_\#(\rho-\rho^1_\eps)\bm1_{\Omega_H''}\right) \\
&= \frac{C_H}{\tau} \int_{\Omega_H''} \abs{S_\eps(x)-x}^2 (\rho(x)-\rho^2_\eps(x)) \d x \\
&\leq \frac{C_H}{\tau} \int_\setR \abs{S_\eps(x)-x}^2 (\rho(x)-\rho^2_\eps(x)) \d x \\
&= \frac{C_H}{\tau} W_2^2\oleft(\rho-\rho^2_\eps, {T_\delta}_\#(\rho-\rho^1_\eps)\right)
\leq C_H\left( \frac{\delta^4}{\tau^2} + \frac{\eps N^2}{\beta\tau^2}
	+ \frac{\eps^{1/2}N}{\beta \tau}\right) .
\end{split}
\]

As a final remark, notice that even without invoking \eqref{eqn:BJR} one could directly verify that the map between  $\rho-\rho^2_\eps$ and ${T_\delta}_\#(\rho-\rho^1_\eps)$  implicitly constructed above (via the curves of measures proposed in Steps 1,2,3) satisfies an $L^\infty$ estimate of the type \eqref{eqn:deconv-hp-verif2}. The proof would meet more technical difficulties than this Step 4, but would however be sufficient.
\end{proof}

\begin{lemma}\label{lem:stretching}
Let $\mu\in\Meas_+([x_0,x_0+\ell])$, $\mu\ll\leb^1$, $\lambda_0,\lambda_1\in\setR_+$ and $S_\lambda:\setR\to\setR:x\mapsto x_0+\lambda(x-x_0)$. Then
$\mu_t=(S_{(1-t)\lambda_0+t\lambda_1})_\#\mu$ solves $\partial_t\mu_t=\div(m_t)$ with
\[
\abs{m_t} \leq \ell\frac{\abs{\lambda_1-\lambda_0}\max\{\lambda_0,\lambda_1\}}{\min\{\lambda_0,\lambda_1\}^2}\norm{\mu}_\infty .
\]

Moreover,
\[
\norm{\mu_t-\mu_s}_{L^\infty(0,\infty)} \leq C \abs{\lambda_t-\lambda_s}+ \ell \Lip(\mu)
	\max\left\{ \frac{\lambda_t}{\lambda_s}, \frac{\lambda_s}{\lambda_t} \right\}.
\]
\end{lemma}

\begin{proof}
Without loss of generality, up to translation we may assume $x_0=0$.
Let $\lambda_t=(1-t)\lambda_0+t\lambda_1$. Then $\mu_t(x)=\frac1{\lambda_t}\mu\oleft(\frac x{\lambda_t}\right)$. The measure $\mu_t$ is advected by the vector field $v_t(x)=\frac{\lambda'_t}{\lambda_t}x=\frac{\lambda_1-\lambda_0}{\lambda_t}x$: in fact
\[
\partial_t\mu_t + \div(v_t\mu_t)
= \frac{\partial}{\partial t}\left[\frac1{\lambda_t}\mu\oleft(\frac x{\lambda_t}\right)\right]
+ \frac{\partial}{\partial x} \left[\frac{\lambda'_t}{\lambda_t^2}x\mu\left(\frac x{\lambda_t}\right)\right] = 0.
\]
Therefore $m_t=-v_t\mu_t$ and
\[
\abs{m_t} \leq \ell\frac{\abs{\lambda_1-\lambda_0}\max\{\lambda_0,\lambda_1\}}{\min\{\lambda_0,\lambda_1\}^2}\norm{\mu}_\infty
\]
because $\abs{x} < \max{\lambda_0,\lambda_1}\ell$ where $\mu(x/\lambda_t)>0$.

The $L^\infty$ estimate comes from the explicit formula for the density of $\mu_t$.
\end{proof}

%{\color{gray}
%\begin{proposition}
%Let $\sigma^1=\rho-\rho^1_\eps$ and $\sigma^1=\rho-\rho^2_\eps$. Let $\pi$ be given by \autoref{cor:pot-cost-rem-mass}. Define $\pi_\eps=\pi*\gamma_\eps$, where $\gamma_\eps$ is the rescaling of the standard Gaussian in $\setR^2$ with variance $\eps^{1/2}$. Then
%\[
%\int V\d\pi_\eps \leq \sqrt\eps\norm{\pi}_1.
%\]
%\end{proposition}
%
%\begin{proof}
%\end{proof}
%}

\subsection{Kinetic energy of the remaining mass}\label{sec:kin-rem}

\begin{proof}[Proof of \autoref{prop:cost-rem-mass}(ii)]
We start by showing the bound for $\KE(\rho-\rho^1_\eps)$. We split the kinetic energy in the domain $S=\supp\rho^1_\eps\subset\Omega_H'$ and the complement, making use of \eqref{eq:lower-marginal} of \autoref{lem:approx-plan} and the fact that in $S^c=(\supp\rho^1_\eps)^c$ we have $\rho-\rho^1_\eps=\rho$. What we obtain is
%\[
%\begin{split}
%\KE(\rho-\rho^1_\eps)
%&= \int_{\supp\rho^1_\eps} \frac{\abs{\partial_x[\rho(x)-\rho^1_\eps(x)]}^2}
%	{4[\rho(x)-\rho^1_\eps(x)]} \d x \\
%&= \int_{\supp\rho^1_\eps} \frac{\abs{\partial_x[\rho(x)-\rho^1_\eps(x)]}^2}
%	{4[\rho(x)-\rho^1_\eps(x)]} \d x
%	+ \int_{(\supp\rho^1_\eps)^c} \frac{\abs{\partial_x[\rho(x)-\rho^1_\eps(x)]}^2}
%	{4[\rho(x)-\rho^1_\eps(x)]} \d x \\
%&\leq \frac{C_H}{\tau} \int_{\supp\rho^1_\eps} \oleft( \abs{\partial_x\rho(x)}^2
%		+ \abs{\partial_x\rho^1_\eps(x)}^2 \right) \d x
%	+ \int_{(\supp\rho^1_\eps)^c} \frac{\abs{\partial_x \rho(x)}^2}{4\rho(x)} \d x \\
%&\leq \frac{C_H\norm{\rho}_\infty}{\tau} \int_{\supp\rho^1_\eps}
%		\frac{\abs{\partial_x\rho(x)}^2}{4\rho(x)} \d x
%	+ \frac{C_H\norm{\rho}_\infty}{\tau} \int_{\supp\rho^1_\eps}
%		\abs{\partial_x\rho^1_\eps(x)}^2 \d x
%%	\spliteq
%	+ \int_{(\supp\rho^1_\eps)^c} \frac{\abs{\partial_x \rho(x)}^2}{4\rho(x)} \d x \\
%&\leq 
%\end{split}
%\]
\[
\begin{split}
\KE(\rho-\rho^1_\eps)
&= \int_{S} \frac{\abs{\partial_x[\rho(x)-\rho^1_\eps(x)]}^2}
	{8[\rho(x)-\rho^1_\eps(x)]} \d x \\
&= \int_{S} \frac{\abs{\partial_x[\rho(x)-\rho^1_\eps(x)]}^2}
	{8[\rho(x)-\rho^1_\eps(x)]} \d x
	+ \int_{S^c} \frac{\abs{\partial_x[\rho(x)-\rho^1_\eps(x)]}^2}
	{8[\rho(x)-\rho^1_\eps(x)]} \d x \\
&\leq \frac{C_H}{\tau} \int_{S} \oleft( \abs{\partial_x\rho(x)}^2
		+ \abs{\partial_x\rho^1_\eps(x)}^2 \right) \d x
	+ \int_{S^c} \frac{\abs{\partial_x \rho(x)}^2}{4\rho(x)} \d x \\
&\leq \frac{C_H\norm{\rho}_{L^\infty(\Omega_H')}}{\tau} \int_{S}
		\frac{\abs{\partial_x\rho(x)}^2}{4\rho(x)} \d x
	+ \frac{C_H}{\tau} \int_{S} \abs{\partial_x\rho^1_\eps(x)}^2 \d x
	+ \int_{S^c} \frac{\abs{\partial_x \rho(x)}^2}{4\rho(x)} \d x \\
&\leq \frac{C_HL}{\tau} \KE(\rho)
	+ \frac{C_H}{\tau} \int_{S} \abs{\partial_x\rho^1_\eps(x)}^2 \d x .
\end{split}
\]
Let us recall that the domain $\Omega_H$ defined in \eqref{eqn:def-omegaH} has two connected components which are both covered by an essentially disjoint family of intervals $I^1_i=[a_i,b_i]$ of length $\delta/2<b_i-a_i<\delta$ with
\begin{gather*}
T(r_H) = a_1 < b_1 = a_2 < \dots < b_{k_0} = T(H) , \\
r_H = a_{k_0+1} < b_{k_0+1} = a_{k_0+2} < \dots < b_{k_1} = H .
\end{gather*}
Notice that the cardinality of the intervals is $k_1 \leq C_H/\delta$.
As already mentioned in \eqref{eq:rho1eps-sum}, we have that
\[
\rho^1_\eps(x)
= \sum_{i=1}^{k_1} \Bigl(\bigl[(\rho-\tau)\bm1_{I^1_i}\bigr]*\eta^1_{\eps,i}\Bigr)(x),
\]
where $\eta^1_{\eps,i} = \proj^1_\# \Gamma_{M_{\eps,\beta}(x_i),N}$ were first introduced in the proof of \autoref{lem:approx-plan}. As a consequence, using the fact that the terms in the sum overlap at most twice, we have
\[
\begin{split}
\norm{\partial_x\rho^1_\eps}_2^2
&= \norm*{\sum_{i=1}^{k_1} \bigl[(\rho-\tau)\bm1_{I^1_i}\bigr]
	* \partial_x\eta^1_{\eps,i}}_2^2 \\
&\leq 2\sum_{i=1}^{k_1} \norm*{\bigl[(\rho-\tau)\bm1_{I^1_i}\bigr]
	* \partial_x\eta^1_{\eps,i}}_2^2
\leq 2\sum_{i=1}^{k_1} \norm*{(\rho-\tau)\bm1_{I^1_i}}_2^2
	\norm{\partial_x\eta^1_{\eps,i}}_1^2 .
\end{split}
\]
Using the representation \eqref{eqn:eta-1-repr} we see that the function $\eta^1_{\eps,i}$ is increasing in $(-\infty,0]$, decreasing in $[0,\infty)$ and vanishes at infinity, hence
\[
\begin{split}
\norm{\partial_x\eta^1_{\eps,i}}_1
&= 2\norm{\eta^1_{\eps,i}}_\infty
\leq 2\norm{(h^1_\eps)_{\{\sin\theta\}}}_1 \norm{(h^2_\eps)_{\{\cos\theta\}}}_\infty \\
&= 2(\cos\theta)^{-1} \norm{h^2_\eps}_\infty
= 2(\cos\theta)^{-1} G_{\frac1\beta,N}^{-1}
\leq C_H \beta^{-1/2}.
\end{split}
\]
As a consequence,
\[
\norm{\partial_x\rho^1_\eps}_2^2
\leq \frac{C_H}{\beta} \sum_{i=1}^{k_1} \norm*{(\rho-\tau)\bm1_{I^1_i}}_2^2
= \frac{C_H}{\beta} \norm{\rho-\tau}_{L^2(\Omega_H)}^2
\leq \frac{C_H}{\beta} \norm{\rho}_1 \norm{\rho}_{L^\infty(\Omega_H')}
\leq \frac{C_H}{\beta} .
\]
%On the other hand,
%\[
%\norm{(\rho-\tau)\bm1_{I^1_i}}_2^2
%= \norm{\rho-\tau}_{L^2(I^1_i)}^2
%\leq \norm{\rho}_{L^1(I^1_i)} \norm{\rho}_{L^\infty(I^1_i)}
%\leq \norm{\rho}_1 \norm{\rho}_{L^\infty(\Omega_H')}
%\leq L .
%\]
%Using these two inequalities to estimate the norms, we get
%\[
%\norm{\partial_x\rho^1_\eps}_2^2
%\leq C_H \frac{\beta}{\delta}
%\]
Using this inequality in the original estimate of the kinetic energy we finally get
\[
\KE(\rho-\rho^1_\eps)
\leq \frac{C_H}{\tau} \KE(\rho) + \frac{C_H}{\tau\beta} .
\]

Let us now move on to the analogous estimate for $\KE(\rho-\rho^2_\eps)$. With the same argument as we did at the beginning, splitting the kinetic energy in the domain $S=\supp\rho^2_\eps\subset\Omega_H'$ and the complement we get
\[
\KE(\rho-\rho^2_\eps)
\leq \frac{C_H}{\tau} \KE(\rho)
	+ \frac{C_H}{\tau} \int_\setR \abs{\partial_x\rho^2_\eps(x)}^2 \d x .
\]
As already mentioned in the proof of \eqref{eq:lower-marginal} of \autoref{lem:approx-plan}, we have that
\[
\rho^2_\eps(y) = \sum_{i=1}^{k_1}
	\Bigl(\bigl[{T_\delta}_\#(\rho-\tau)\bm1_{I^2_i}\bigr]*\eta^2_{\eps,i}\Bigr)(y) ,
\]
where $\eta^2_{\eps,i} = \proj^2_\#\Gamma_{M_{\eps,\beta}(x_i),N}$ are the second marginals of the Gaussians. The intervals $[c_i,d_i]=I^2_i=T_\delta(I^1_i)=T(I^1_i)=[T(a_i),T(b_i)]$ have lengths $\delta/(2L)\leq\diam(I^2_i)\leq L\delta$. If $y\in I^2_i$ we have that
\[
{T_\delta}_\#(\rho-\tau)(y) = [\rho(T_\delta^{-1}(y))-\tau]/T'(x_i)
\leq L \norm{\rho}_{L^\infty(\Omega_H')} ;
\]
on the other hand, from \eqref{eqn:eta-2-repr} we have
\[
\begin{split}
\norm{\eta^2_{\eps,i}}_\infty
\leq \norm{(h^1_\eps)_{\{\cos\theta\}}}_1 \norm{(h^2_\eps)_{\{\sin\theta\}}}_\infty
\leq (\sin\theta)^{-1}\norm{h^2_\eps}_\infty
\leq C_H \beta^{-1/2},
\end{split}
\]
therefore
\[
\norm{\partial_y \eta^2_{\eps,i}}_1
= 2\norm{\eta^2_{\eps,i}}_\infty
\leq C_H \beta^{-1/2} .
\]
This implies that
\[
\begin{split}
\norm{\partial_y \rho^2_\eps}_2^2
&= \norm*{\sum_{i=1}^{k_1}
	\bigl[{T_\delta}_\#(\rho-\tau)\bm1_{I^2_i}\bigr]
	* \partial_y \eta^2_{\eps,i}}_2^2
\leq 2\sum_{i=1}^{k_1} \norm*{\bigl[{T_\delta}_\#(\rho-\tau)\bm1_{I^2_i}\bigr]
	* \partial_y \eta^2_{\eps,i}}_2^2 \\
&\leq 2\sum_{i=1}^{k_1} \norm*{{T_\delta}_\#(\rho-\tau)\bm1_{I^2_i}}_2^2
	\norm{\partial_y\eta^2_{\eps,i}}_1^2
\leq \frac{C_H}{\beta} \sum_{i=1}^{k_1} \norm*{{T_\delta}_\#\rho\bm1_{I^2_i}}_1
	\norm*{{T_\delta}_\#\rho\bm1_{I^2_i}}_\infty \\
&\leq \frac{C_H}{\beta} L \norm{\rho}_{L^\infty(\Omega_H')}
	\sum_{i=1}^{k_1} \norm*{\rho\bm1_{I^2_i}}_1
\leq \frac{C_H}{\beta} L \norm{\rho}_{L^\infty(\Omega_H')} \norm*{\rho}_{L^1(\Omega_H')}
\leq \frac{C_H}{\beta} ,
\end{split}
\]
%\[
%\begin{split}
%\norm{{T_\delta}_\#(\rho-\tau)\bm1_{I^2_i}}_2^2
%&= \norm{{T_\delta}_\#(\rho-\tau)}_{L^2(I^2_i)}^2
%\leq \norm{{T_\delta}_\#(\rho-\tau)}_{L^1(I^2_i)}
%	\norm{{T_\delta}_\#(\rho-\tau)}_{L^\infty(I^2_i)} \\
%&\leq \norm{\rho-\tau}_{L^1(I^1_i)} \norm{\rho-\tau}_{L^\infty(I^1_i)} / T'(x_i)
%\leq L^2 .
%\end{split}
%\]
which inserted in the original estimate of the kinetic energy leads to
\[
\KE(\rho-\rho^2_\eps)
\leq \frac{C_H}{\tau} \KE(\rho) + \frac{C_H}{\tau\beta}. \qedhere
\]
\end{proof}

\section{Proof of \autoref{thm:main_intro} and \autoref{thm:rec}}\label{sec:proof}

%By translating we may assume that the origin $0$ is the median of the probability $\rho$, i.e.\ $\rho\bigl((-\infty,0]\bigr)=\rho\bigl([0,\infty)\bigr)=1/2$. %(CONTROLLARE SE LO USIAMO)

\begin{proof}[Proof of \autoref{thm:rec}]
Given $H>1$ and $L=L(H)$ defined by \eqref{eq:L}, let $\bar\gamma_\eps$ be as defined in \eqref{eq:bar-gamma-eps}, built with admissible parameters $N,\beta,\delta,\tau$ satisfying \eqref{eqn:choice-N}-\eqref{eqn:choice-tau}, for instance the explicit choice \eqref{eqn:choice}. By \autoref{lem:approx-plan} we have that
\[
\limsup_{\eps\to0} E_\eps(\sqrt{\bar \gamma_\eps})
\leq \frac12 \int_{\Omega_H} \tr\oleft(\sqrt{\nabla^2V(x, T(x))}\right) \rho(x) \d x
%\limsup_{\eps\to0} E_\eps(\sqrt{ \gamma_\eps})
\leq \frac12 \int_{\setR^2} \tr\bigl(\sqrt{D^2V}\bigr) \d\gamma .
\]
The problem is that $\bar \gamma_\eps$ is not a recovery sequence because it has the wrong mass and marginals, so we have to use the remaining mass to fix the marginals.

For $i=1,2$, let $\upsilon^i_\eps = \proj^i_\#(\gamma) - \proj^i_\#(\bar \gamma_\eps) = \rho - \rho^i_\eps$, which is a positive measure as remarked in \eqref{eq:lower-marginal} of \autoref{lem:approx-plan}.% with total mass computed in \eqref{eq:remaining-mass} and given by $ \rho(\Omega_H^c) + \tau\abs{\Omega_H} $.

Let $\pi_\eps\in\Pi(\upsilon^1_\eps,\upsilon^2_\eps)$ be the plan given by \autoref{prop:cost-rem-mass}(i). Notice that, by \eqref{eq:remaining-mass}, $\norm{\upsilon^i_\eps}_1=\norm{\pi_\eps}_1\leq C_H\tau + \rho(\Omega_H^c)$. With $H$ fixed the first term goes to $0$ when $\eps\to0$, whereas the second term goes to zero as $H\to\infty$.

Let $\tilde\pi_\eps \in \Pi(\upsilon^1_\eps,\upsilon^2_\eps)$ be the deconvolved plan defined in \autoref{sec:deconvolution} starting from $\Pi_0=\pi_\eps$ (whose marginals are $\sigma^1=\upsilon^1_\eps$ and $\sigma^2=\upsilon^2_\eps$)%, which clearly has the same total mass as $\pi_\eps$
. By \autoref{prop:cost-deconvolved-plan} (notice that the assumption \eqref{eqn:deconv-hp} is satisfied because of \eqref{eqn:deconv-hp-verif}) and the bound on the mass of $\pi_\eps$ we have
\[
\begin{split}
E_\eps(\sqrt{\tilde\pi_\eps})
&= \eps^{1/2}\KE(\tilde\pi_\eps) + \eps^{-1/2} \int_{\setR^2} V\d\tilde\pi_\eps \\
&\leq \eps^{1/2}[\KE(\upsilon^1_\eps)+\KE(\upsilon^2_\eps)]
	+ C_H\eps^{-1/2} \int_{\setR^2} V\d\pi_\eps + C\norm{\pi_\eps}_1 \\
&\leq \eps^{1/2}[\KE(\upsilon^1_\eps)+\KE(\upsilon^2_\eps)]
	+ C_H\eps^{-1/2}\int_{\setR^2} V\d\pi_\eps + C_H\tau + C\rho(\Omega_H^c).
\end{split}
\]
By \autoref{prop:cost-rem-mass}, we have  $\lim_{\eps\to 0}\eps^{1/2}\int_{\setR^2} V\d\pi_\eps=0= \lim_{\eps\to 0}\eps^{-1/2}\KE(\upsilon^i_\eps) $, therefore
\begin{equation}\label{eq:energia-rimanente-sparisce}
\lim_{\eps\to 0}E_\eps(\sqrt{\tilde\pi_\eps}) \leq C\rho(\Omega_H^c).
\end{equation}

Finally, let $\psi_\eps = \sqrt{\bar\gamma_\eps+\tilde\pi_\eps}$; then $\psi_\eps$ has the correct marginals to be a recovery sequence for $\gamma$.
Moreover, by the subadditivity of $E_\eps(\sqrt\gamma)$ with respect to $\gamma$ and \eqref{eq:energia-rimanente-sparisce}, we have
\[\begin{split}
\limsup_{\eps\to0} E_\eps(\psi_\eps)
&\leq \limsup_{\eps\to0} [E_\eps(\sqrt{\bar\gamma_\eps}) + E_\eps(\sqrt{\tilde\pi_\eps})] \\
&\leq \limsup_{\eps\to0} E_\eps(\sqrt{\bar\gamma_\eps}) + C\rho(\Omega_H^c)
\leq \frac12 \int_{\setR^2} \tr\bigl(\sqrt{D^2V}\bigr) \d\gamma + C\rho(\Omega_H^c).
\end{split}
\]

We now use a diagonal argument to conclude. For every $n\in\setN_+$ we can find $H_n$ large enough such that $C\rho(\Omega_{H_n}^c)<1/n$. We can then find $\eps_n$ such that the $\psi_{\eps_n}$ constructed as above with parameters $H_n,\eps_n,\beta_n,\delta_n,\tau_n,N_n$ satisfies
\[
E_{\eps_n}(\psi_{\eps_n})
\leq \frac12 \int_{\setR^2} \tr\bigl(\sqrt{D^2V}\bigr) \d\gamma
	+ C\rho(\Omega_{H_n}^c) + 1/n
\leq \frac12 \int_{\setR^2} \tr\bigl(\sqrt{D^2V}\bigr) \d\gamma + 2/n,
\]
which is precisely the thesis of the theorem.
\end{proof}
\begin{proof}[Proof of \autoref{thm:main_intro}]

To prove \eqref{eq:Eeps-liminf}, we take a sequence $\psi_\eps$ which almost realizes the minimum in \eqref{eq:Feps-FOT} up to a vanishing error. Up to subsequence, $\abs{\psi_\eps}^2\weakto\gamma$ where $\gamma\in\Pi_0(\rho)$ with $\gamma(\{V>0\})=0$. Applying \autoref{thm:nonfixed-marginals} to this sequence we get
\[
\liminf_{\eps\to0} \frac{F^\eps(\rho)-F_{OT}(\rho)}{\sqrt\eps}
\geq \liminf_{\eps\to0} E_\eps(\psi_\eps)
\geq \frac12 \int_{\setR^{dN}} \tr\oleft(\sqrt{D^2V}\right) \d\gamma
= F_{ZPO}(\rho) .
\]

To prove \eqref{eq:Eeps-limsup}, we take $\gamma\in\Pi_0(\rho)$ with $\gamma(\{V>0\})=0$ and $F_{ZPO}(\rho) = \frac12 \int_{\setR^{dN}} \tr\oleft(\sqrt{D^2V}\right) \d\gamma$. If $\psi_\eps$ is a recovery sequence given by \autoref{thm:rec}, we get
\[
\limsup_{\eps\to0} \frac{F^\eps(\rho)-F_{OT}(\rho)}{\sqrt\eps}
\leq \limsup_{\eps\to0} E_\eps(\psi_\eps) \leq F_{ZPO}(\rho). \qedhere
\]

%Let's start with the first claim, which is a direct consequence of \autoref{thm:nonfixed-marginals}. Let $\gamma\in\Pi_0(\rho)$ be such that
%\[
%F_{ZPO}(\rho) = \frac12 \int_{\setR^{dN}} \tr\oleft(\sqrt{D^2V}\right) \d\gamma
%\]
%and let $\psi_\eps\mapsto\rho$ be a sequence such that
%\[
%\frac{F^\eps(\rho)-F_{OT}(\rho)}{\sqrt\eps} \geq E_\eps(\psi_\eps) - \eps.
%\]
%From \autoref{thm:nonfixed-marginals}~(a) applied to the sequence $\psi_\eps=\rho$ we get that
%\[
%\liminf_{\eps\to0} \frac{F^\eps(\rho)-F_{OT}(\rho)}{\sqrt\eps}
%\geq \liminf_{\eps\to0} E_\eps(\psi_\eps)
%\geq \frac12 \int_{\setR^{dN}} \tr\oleft(\sqrt{D^2V}\right) \d\gamma
%= F_{ZPO}(\rho),
%\]
%which is the thesis.
%
%Let us now deal with the second part of the statement, which will follow from \autoref{thm:rec} instead.
\end{proof}
\section*{Appendix: proof of \autoref{lem:norm-diff-gauss} and \autoref{lem:diff-energies}}\label{sec:proof-lemmas}

\begin{proof}[Proof of \autoref{lem:norm-diff-gauss}]
Notice that $G_{M,N}<G_{M,\infty}$ by definition, so we only need to prove the second inequality of \eqref{eq:GN-infty}.
We have $G_{\alpha,N}<G_{\alpha,\infty}$ and we claim that
\begin{equation}\label{eq:Gainfty-GaN}
G_{\alpha,\infty} - G_{\alpha,N} \leq 3e^{-N/2} G_{\alpha,\infty} .
\end{equation}
By a change of variable we reduce to prove it only for $\alpha=1$. We apply the pointwise inequality 
\begin{equation}
\label{eqn:f-n}f^2\leq (f-e^{-N/2})^2 + 2fe^{-N/2}%\qquad \mbox{for every } f\geq 0, N\geq 0
\end{equation}
 to $f= e^{-t^2/2}$ for $t< \sqrt N$ and use $N\geq1$ to get
\begin{equation*}
\begin{split}
G_{1,\infty}
&= \int_{\{\abs{t}\leq\sqrt N\}} e^{-t^2} \d t
	+ \int_{\{\abs{t}>\sqrt N\}} e^{-t^2} \d t \\
&\leq \int_{\{\abs{t}\leq\sqrt N\}} \left[
		\left(e^{- t^2/2}-e^{-N/2}\right)^2 + 2e^{-t^2/2}e^{-N/2} \right] \d t
	+ \int_{\{\abs{t}> \sqrt N\}} t e^{- t^2}  \d t \\
&\leq G_{1,N} + 2 \sqrt{2\pi}e^{-N/2} + e^{-N}
= G_{1,N} + \left(2\sqrt{2} + \frac{e^{-N/2}}{\sqrt\pi}\right) e^{-N/2} G_{1,\infty} ,
\end{split}
\end{equation*}
which proves the claim \eqref{eq:Gainfty-GaN} for $N\geq3$, because the factor inside parentheses is less than $3$. From this, \eqref{eq:GN-infty} follows because
\begin{equation*}%\label{eq:GMinfty-GMN}
\begin{split}
G_{M,\infty} - G_{M,N}
&= G_{a,\infty} (G_{b,\infty}-G_{b,N}) + G_{b,N} (G_{a,\infty} - G_{a,N}) \\
&\leq G_{a,\infty}3e^{-N/2}G_{b,\infty} + G_{b,N}3e^{-N/2}G_{a,\infty} \\
&\leq 6e^{-N/2}G_{a,\infty}G_{b,\infty} = 6e^{-N/2}G_{M,\infty} .
\end{split}
\end{equation*}
As a consequence, we get also
\begin{equation}\label{eq:1_GMinfty-1_GMN}
\frac1{G_{M,N}}-\frac1{G_{M,\infty}}
= \frac{G_{M,\infty}-G_{M,N}}{G_{M,\infty}G_{M,N}}
\leq \frac{6e^{-N/2}}{G_{M,N}}
\leq \frac{6e^{-N/2}}{(1-6e^{-N/2})G_{M,\infty}}
\leq \frac{7e^{-N/2}}{G_{M,\infty}} .
\end{equation}

Let now $f=e^{-aw^2/2}$ and $g=e^{-bz^2/2}$. Notice that $f,g\leq1$ %, $f^2\leq \bigl(f-e^{-N/2}\bigr)_+^2 + 2fe^{-N/2}$ and similarly for $g$
and that we can apply \eqref{eqn:f-n} to $f$ and $g$.
 Therefore
\[
\begin{split}
\abs*{\tilde\Gamma_{M,\infty}(\x)-\tilde\Gamma_{M,N}(\x)}
&= f^2g^2 - \bigl(f-e^{-N/2}\bigr)_+^2 \bigl(g-e^{-N/2}\bigr)_+^2 \\
&= f^2 \left(g^2 - \bigl(g-e^{-N/2}\bigr)_+^2\right)
%		\spliteq
	+ \left(f^2 - \bigl(f-e^{-N/2}\bigr)_+^2\right) \bigl(g-e^{-N/2}\bigr)_+^2 \\
&\leq 2 (f^2g+fg^2) e^{-N/2} \leq 4 e^{-N/2}.
\end{split}
\]

%{\color{gray}
%Let $a=e^{-\left(\frac{q}{\eps^{1/2}}+\frac1\beta\right)w^2/2}$ and $b=e^{-\frac1\beta z^2/2}$. Notice that $a,b\leq 1$. If either $a\leq e^{-N/2}$ or $b\leq e^{-N/2}$, then
%\[
%\abs*{\tilde\Gamma_{M,\infty}(x)-\tilde\Gamma_{M,N}(x)}
%= a^2b^2 - \bigl(a-e^{-N/2}\bigr)_+^2\bigl(b-e^{-N/2}\bigr)_+^2
%= a^2b^2-0 \leq e^{-N}.
%\]
%Otherwise, if $\min\{a,b\}>e^{-N/2}$, then
%\[
%\begin{split}
%\abs*{\tilde\Gamma_{M,\infty}(x)-\tilde\Gamma_{M,N}(x)}
%&= a^2b^2 - \bigl(a-e^{-N/2}\bigr)^2\bigl(b-e^{-N/2}\bigr)^2 \\
%&= a^2b^2 - \bigl(a^2-2ae^{-N/2}+e^N\bigr) \bigl(b^2-2be^{-N/2}+e^N\bigr) \\
%&= 2(a^2b+ab^2)e^{-N/2} - (a^2+4ab+b^2)e^{-N} + 2(a+b)e^{-3N/2} - e^{-2N} \\
%&\leq 5 e^{-N/2}
%\end{split}
%\]
%for $N$ sufficiently large. In any case, $\norm{\tilde\Gamma_{M,\infty}(x)-\tilde\Gamma_{M,N}(x)}_\infty \leq 5e^{-N/2}$.
%}

We can now use this information and \eqref{eq:1_GMinfty-1_GMN} to estimate the difference of the two normalized Gaussians as
\[
\begin{split}
\norm{\Gamma_{M,\infty}-\Gamma_{M,N}}_\infty
%&\mathrel{\phantom{=}}
%&= \norm*{ \frac{\tilde\Gamma_{M,\infty}}{G_{M,\infty}}
%	- \frac{\tilde\Gamma_{M,N}}{G_{M,N}} }_\infty \\
&= \norm*{ \frac{\tilde\Gamma_{M,\infty} - \tilde\Gamma_{M,N}}{G_{M,\infty}}
	- \tilde\Gamma_{M,N} \left(\frac1{G_{M,N}} - \frac1{G_{M,\infty}}\right)}_\infty \\
&= \frac{\norm{\tilde\Gamma_{M,\infty}(x)-\tilde\Gamma_{M,N}(x)}_\infty}{G_{M,\infty}}
	+ \norm{\tilde\Gamma_{M,N}}_\infty
		\left(\frac1{G_{M,N}} - \frac1{G_{M,\infty}}\right) \\
&\leq \frac{4e^{-N/2}}{G_{M,\infty}} + \frac{7e^{-N/2}}{G_{M,\infty}}
= \frac{11e^{-N/2}}{G_{M,\infty}}
= \frac{11}\pi \sqrt{\det M} e^{-N/2} ,
\end{split}
\]
proving \eqref{eq:Linfty-norm-trunc}. 
To prove the $L^1$ estimate \eqref{eq:L1-norm-trunc}, we divide the integral on $\supp\Gamma_{M,N}$ and the complement. In the first, thanks to \eqref{eq:Linfty-norm-trunc} we obtain
\[
\begin{split}
\norm{\Gamma_{M,N}-\Gamma_{M,\infty}}_{L^1(\supp\Gamma_{M,N})} &\leq 
 \norm{\Gamma_{M,N}-\Gamma_{M,\infty}}_\infty
	\cdot \abs{\supp\Gamma_{M,N}}
%	+ \left( 1 - \int_{\supp\Gamma_{M,N}} \Gamma_{M,\infty}(\x) \d \x \right)
\\&\leq C\sqrt{\det M}e^{-N/2} \frac{4N}{\sqrt{ab}} \leq C e^{-N/2} N
%	+ \left( 1 - \frac{G_{M,N}}{G_{M,\infty}} \right) \\
%&\leq C e^{-N/2} N + 6 e^{-N/2}
%\leq C e^{-N/2}N 
.
\end{split}
\]
In the complement, only $\Gamma_{M,\infty}$ is nonzero and by \eqref{eq:GN-infty} we have
\[
\begin{split}
\norm{\Gamma_{M,N}-\Gamma_{M,\infty}}&_{L^1((\supp\Gamma_{M,N})^c)} \\
&= \int_{(\supp\Gamma_{M,N})^c} \Gamma_{M,\infty}(\x) \d \x
= 1 - \int_{\supp\Gamma_{M,N}} \Gamma_{M,\infty}(\x) \d \x \\
&\leq  1 - \int_{\supp\Gamma_{M,N}} \frac{\tilde\Gamma_{M,N}(\x)}{G_{M,\infty}} \d \x
= 1 - \frac{G_{M,N}}{G_{M,\infty}} \leq C e^{-N/2}N .
\end{split}
\]

Let us now prove \eqref{eq:Linfty-norm-trunc-marg}. Define $R_x$ to be the section above $x$ of the rectangle $\supp\Gamma_{M,N}$, i.e.
\[
R_x = \set{y}{(x,y)\in\supp\Gamma_{M,N}}
= \set{y}{aw^2\leq N,\ bz^2\leq N}.
\]
Since $b y^2 \leq b(w^2+z^2) \leq aw^2+bz^2 \leq a(w^2+z^2) = a(x^2+y^2)$, we have that
\[
\set{(x,y)\in\setR^2}{a(x^2+y^2)\leq N} \subset \supp\Gamma_{M,N}
\]
and
\[
R_x \subseteq \set{y}{by^2\leq 2N} ,\qquad
R_x^c \subseteq \set{y}{ay^2 > N-ax^2}.
\]
Therefore, using \eqref{eq:Linfty-norm-trunc}, we can estimate
\[
\begin{split}
\abs{\eta_{M,\infty}(x) - \eta_{M,N}(x)}
&\leq \int_\setR \abs{\Gamma_{M,\infty}(x,y)-\Gamma_{M,N}(x,y)} \d y \\
&= \int_{R_x} \abs{\Gamma_{M,\infty}(x,y)-\Gamma_{M,N}(x,y)} \d y
	+ \int_{R_x^c} \Gamma_{M,\infty}(x,y) \d y . \\
%&\leq \abs{R_x} \cdot \norm{\Gamma_{M,\infty}-\Gamma_{M,N}}_\infty
%	+ G_{M,\infty}^{-1} \int_{R_x^c} e^{-(aw^2+bz^2)} \d y \\
%&\leq 2\sqrt2\frac{\sqrt N}{\sqrt{b}} C\sqrt{ab}e^{-N/2}
%	+ \frac{\sqrt{ab}}{\pi} \int_{\{a\abs{y}^2>N-ax^2\}} e^{-ax^2-ay^2} \d y \\
%&\leq 2\sqrt2\sqrt{N}C\sqrt{a}e^{-N/2}
%	+ \frac{\sqrt{ab}}{\pi} \frac{\sqrt\pi}{\sqrt a} e^{-ax^2}
%		\erfc\oleft(\sqrt{N-ax^2}\right) \\
%&\leq 2\sqrt2\sqrt{N}C\sqrt{a}e^{-N/2}
%	+ \frac{\sqrt{b}}{\sqrt\pi} e^{-ax^2} e^{-N+ax^2} \\
%&\leq C\sqrt a \sqrt{N} e^{-N/2} .
\end{split}
\]

The first integral can be estimated with \eqref{eq:Linfty-norm-trunc} as
\[
\begin{split}
\int_{R_x} \abs{\Gamma_{M,\infty}(x,y)-\Gamma_{M,N}(x,y)} \d y
&\leq \abs{R_x} \cdot \norm{\Gamma_{M,\infty}-\Gamma_{M,N}}_\infty \\
&\leq \left(2\sqrt2\frac{\sqrt N}{\sqrt{b}}\right) \left(C\sqrt{ab}e^{-N/2}\right)
\leq C\sqrt a \sqrt{N} e^{-N/2} .
\end{split}
\]

For the second integral we proceed as follows.
Let $B=B(0,\sqrt{N/b}) = \{b(w^2+z^2)<N\} = \{b(x^2+y^2)<N\}$. Note that in $R^c_x$ we have $e^{-aw^2-bz^2}<e^{-N}$ because at least one among $aw^2>N$ and $bz^2>N$ is true.
As a consequence, we have that
\[
\int_{\R^c_x\cap B} e^{-aw^2-bz^2}\d y
\leq e^{-N} \diam(B) = \frac{2}{\sqrt b}e^{-N} \sqrt N
\]
while on the other hand
\[
\begin{split}
\int_{\R^c_x\setminus B} e^{-aw^2-bz^2}\d y
&\leq \int_{\R^c_x\setminus B} e^{-bw^2-bz^2} \d y
= \int_{\{by^2>N-bx^2\}} e^{-bx^2-by^2} \d y \\
&= \frac{\sqrt\pi}{\sqrt b} e^{-bx^2} \erfc\oleft(\sqrt{N-bx^2}\right)
\leq \frac{\sqrt\pi}{\sqrt b} e^{-bx^2} e^{-N+bx^2}
= \frac{\sqrt\pi}{\sqrt b} e^{-N} ,
\end{split}
\]
therefore
\[
\begin{split}
\int_{\R^c_x} e^{-aw^2-bz^2}\d y
= \int_{\R^c_x\cap B} e^{-aw^2-bz^2}\d y + \int_{\R^c_x\setminus B} e^{-aw^2-bz^2}\d y
\leq \frac{C}{\sqrt b} e^{-N}\sqrt N ,
\end{split}
\]
hence
\[
\int_{R_x^c} \Gamma_{M,\infty}(x,y) \d y
= G_{M,\infty}^{-1} \int_{R_x^c} e^{-(aw^2+bz^2)} \d y
\leq \frac{\sqrt{ab}}{\pi} \frac{C}{\sqrt b} e^{-N}\sqrt N
\leq C \sqrt{a} e^{-N}\sqrt N.
\]

In conclusion, putting the two estimates together, we have
\[
\abs{\eta_{M,\infty}(x) - \eta_{M,N}(x)}
\leq C\sqrt a \sqrt{N} e^{-N/2} . \qedhere
\]
\end{proof}

\begin{proof}[Proof of \autoref{lem:diff-energies}]
In some computations of this proof we use the standard error function and the complementary error function
\begin{align*}
\erf(z) &= \frac2{\sqrt\pi} \int_0^z e^{-t^2}\d t, &
\erfc(z) &= \frac2{\sqrt\pi} \int_z^\infty e^{-t^2}\d t = 1-\erf(z).
\end{align*}
The crucial property that we will need is the fast decay at infinity of $\erfc$ implied by the bound
\[
\erfc(z) \leq \frac2{\sqrt\pi} \int_z^\infty \frac{t}{z}e^{-t^2}\d t
\leq \frac{e^{-z^2}}{\sqrt\pi z} .
\]
%and it is well known that for $z\to\infty$
%\[
%\erfc(z) = e^{-z^2}\left(\frac1{\sqrt\pi z} + o(z^{-1})\right),
%\]
%which can be readily proved with de l'Hôpital.

Since
\[
(\supp\Gamma_{M,N})^c \subset \{\abs{w}>\sqrt{N/a}\} \cup \{\abs{z}>\sqrt{N/b}\},
\]
with the aid of \eqref{eq:L1-norm-trunc} we can estimate
\[
\begin{split}
%\biggl\lvert \int_{\setR^2} &\abs{M\x}^2 \Gamma_{M,N}(\x) \d\x
%	- \int_{\setR^2} \abs{M\x}^2 \Gamma_{M,\infty}(\x) \d\x \biggr\rvert \\
\int_{\setR^2} & \abs{M\x}^2 \abs{\Gamma_{M,N}(\x)-\Gamma_{M,\infty}(\x)} \d\x \\
&\leq \int_{\supp\Gamma_{M,N}} \abs{M\x}^2
		\abs{\Gamma_{M,N}(\x)-\Gamma_{M,\infty}(\x)} \d\x
%	\spliteq
	+ \int_{(\supp\Gamma_{M,N})^c} \abs{M\x}^2 \Gamma_{M,\infty}(\x) \d\x \\
&\leq \left( \sup_{\supp\Gamma_{M,N}} \abs{M\x}^2 \right)
		\norm{\Gamma_{M,N}-\Gamma_{M,\infty}}_1
	\spliteq
	+ G_{M,\infty}^{-1} 2\int_{\sqrt{N/a}}^\infty \int_\setR
		(a^2w^2+b^2z^2) e^{-aw^2-bz^2} \d z \d w
	\spliteq
	+ G_{M,\infty}^{-1} 2\int_{\sqrt{N/b}}^\infty \int_\setR
		(a^2w^2+b^2z^2) e^{-aw^2-bz^2} \d w \d z \\
&\leq (a+b)CNe^{-N/2} + (a+b)\left(
	\erfc(\sqrt\pi) + \frac{\sqrt{N}}{\sqrt\pi} e^{-N} \right)
\leq C\tr(M) N e^{-N/2} .
\end{split}
\]

%We can compute
%\[
%\nabla\sqrt{\Gamma_{M,\infty}}
%= G_{M,\infty}^{-1/2} \nabla\bigl(e^{-aw^2/2}e^{-bz^2/2}\bigr)
%= -G_{M,\infty}^{-1/2} \begin{pmatrix} aw \\ bz \end{pmatrix}
%	\sqrt{\tilde\Gamma_{M,\infty}},
%\]
%therefore the kinetic energy of the whole Gaussian is
%\[
%\begin{split}
%\KE(\Gamma_{M,\infty})
%%&= G_{M,\infty}^{-1} \KE(\tilde\Gamma_{M,\infty})
%%&= \int_{\setR^2} \abs*{\nabla\sqrt{\Gamma_{M,\infty}}}^2 \d w \d z
%&= G_{M,\infty}^{-1} \int_{\setR^2} (a^2w^2+b^2z^2)e^{-aw^2}e^{-bz^2} \d w \d z \\
%&= G_{M,\infty}^{-1} \left(
%	\frac{\sqrt\pi}{\sqrt b} \int_\setR a^2 w^2 e^{-aw^2} \d w
%	+ \frac{\sqrt\pi}{\sqrt a} \int_\setR b^2 z^2 e^{-bz^2} \d z \right) \\
%&= G_{M,\infty}^{-1} \left( \frac{\sqrt\pi}{\sqrt b} \cdot \frac{\sqrt a \sqrt\pi}{2}
%	+ \frac{\sqrt\pi}{\sqrt a} \cdot \frac{\sqrt b \sqrt\pi}{2} \right) \\
%&= \frac{\sqrt{ab}}{\pi} \cdot \frac{\pi}{2} \left( \frac{a+b}{\sqrt{ab}} \right)
%= \frac{a+b}{2} = \frac{\tr M}{2} .
%\end{split}
%\]

For the kinetic energy, we can compute
\begin{equation}\label{eq:gradient-sqrt-gamma}
\begin{split}
\nabla_{w,z}\sqrt{\Gamma_{M,N}(\x)}
&= G_{M,N}^{-1/2} \begin{pmatrix}
-aw e^{-aw^2/2} \left(e^{-bz^2/2} - e^{-N/2}\right)_+ \\
-bz e^{-bz^2/2} \left(e^{-aw^2/2} - e^{-N/2}\right)_+
\end{pmatrix} \bm{1}_{\supp\Gamma_{M,N}},
\end{split}
\end{equation}
hence
\[
\begin{split}
\abs*{\nabla\sqrt{\Gamma_{M,N}(\x)}}^2
&= G_{M,N}^{-1} \biggl[
a^2 w^2 e^{-aw^2} \left(e^{-bz^2/2} - e^{-N/2}\right)_+^2
	\spliteq\hspace{2cm}
+ b^2 z^2 e^{-bz^2} \left(e^{-aw^2/2} - e^{-N/2}\right)_+^2
\biggr] \bm{1}_{aw^2\leq N, bz^2\leq N}
\end{split}
\]
and the kinetic energy of the truncated Gaussian ends up being
\[
\begin{split}
\KE&(\Gamma_{M,N}) \\
%&= \int_{\setR^2} \abs*{\nabla\sqrt{\Gamma_{M,N}(\x)}}^2 \\
&= \frac12 G_{M,N}^{-1} \left(\int_{\{aw^2\leq N\}} a^2 w^2 e^{-aw^2} \d w\right)G_{b,N}
%	\spliteq
	+ \frac12 G_{M,N}^{-1} \left(\int_{\{bz^2\leq N\}} b^2 z^2 e^{-bz^2} \d z\right)G_{a,N} \\
&= \frac12 G_{a,N}^{-1} \sqrt{a} \left(\frac{\sqrt\pi}{2}\erf(\sqrt{N})
		- e^{-N}\sqrt{N}\right)
%	\spliteq
	+ \frac12 G_{b,N}^{-1} \sqrt{b} \left(\frac{\sqrt\pi}{2}\erf(\sqrt{N})
		- e^{-N}\sqrt{N}\right) \\
&= (G_{a,N}^{-1}\sqrt{a} + G_{b,N}^{-1}\sqrt{b})
	\frac{\sqrt\pi}{4} \left(\erf(\sqrt N) - \frac{2}{\sqrt\pi}e^{-N}\sqrt{N}\right) ,
\end{split}
\]
where to pass from the first to the second line we integrated by parts and changed variables.
Recalling that $G_{\alpha,\infty}^{-1}=\sqrt{\alpha}/\sqrt\pi$, for $N=\infty$ this gives
\[
\KE(\Gamma_{M,\infty})
= \left(\frac{\sqrt{a}}{\sqrt\pi}\sqrt{a} + \frac{\sqrt{b}}{\sqrt\pi}\sqrt{b}\right)
	\frac{\sqrt\pi}{4}
= \frac{a+b}{4} = \frac{\tr M}{4},
\]
as already stated in \eqref{eqn:daqualcheparte}.
In order to compare the two energies, we estimate separately the difference of the terms involving $a$ and $b$. Using the fact that, thanks to \eqref{eq:Gainfty-GaN}, %\eqref{eq:1_GMinfty-1_GMN},
\[
\frac1{G_{\alpha,N}}-\frac1{G_{\alpha,\infty}}
= \frac{G_{\alpha,\infty}-G_{\alpha,N}}{G_{\alpha,\infty}G_{\alpha ,N}}
\leq \frac{3e^{-N/2}}{G_{\alpha,N}}
\leq \frac{3e^{-N/2}}{(1-3e^{-N/2})G_{\alpha,\infty}}
\leq \frac{4e^{-N/2}}{G_{\alpha,\infty}} ,
\]
we can estimate
\[
\begin{split}
\biggl\lvert \frac{a}{4} &- G_{a,N}^{-1} \sqrt{a}
	\frac{\sqrt\pi}{4} \left(\erf(\sqrt N) - \frac{2}{\sqrt\pi}e^{-N}\sqrt{N}\right)
	\biggr\rvert \\
&\leq \abs*{\frac{a}{4} - G_{a,\infty}^{-1} \sqrt{a}
	\frac{\sqrt\pi}{4} \left(\erf(\sqrt N) - \frac{2}{\sqrt\pi}e^{-N}\sqrt{N}\right)}
		\spliteq
	+ \abs*{(G_{a,N}^{-1}-G_{a,\infty}^{-1}) \sqrt{a}
	\frac{\sqrt\pi}{4} \left(\erf(\sqrt N) - \frac{2}{\sqrt\pi}e^{-N}\sqrt{N}\right)} \\
&\leq \frac{a}{4} \abs*{\erfc(\sqrt N) + \frac{2}{\sqrt\pi}e^{-N}\sqrt{N}}
	+ 4 e^{-N/2} G_{a,\infty}^{-1} \sqrt{a} \frac{\sqrt\pi}{4}
		\abs*{\erf(\sqrt N) - \frac{2}{\sqrt\pi}e^{-N}\sqrt{N}} \\
&\leq C \frac{a}{4} e^{-N/2},
\end{split}
\]
and similarly for $b$, therefore
\[
\begin{split}
\abs*{\KE(\Gamma_{M,\infty}) - \KE(\Gamma_{M,N})}
&\leq C \frac{a+b}{4} e^{-N/2} = C \KE(\Gamma_{M,\infty}) e^{-N/2}. \qedhere
\end{split}
\]
\end{proof}

\phantomsection
\addcontentsline{toc}{section}{\refname}
\printbibliography

\end{document}